# AUXILIARY SDES FOR HOMOGENIZATION OF QUASILINEAR PDES WITH PERIODIC COEFFICIENTS

By François Delarue

*Université Paris VII*

We study the homogenization property of systems of quasi-linear PDEs of parabolic type with periodic coefficients, highly oscillating drift and highly oscillating nonlinear term. To this end, we propose a probabilistic approach based on the theory of forward–backward stochastic differential equations and introduce the new concept of "auxiliary SDEs."

## 1. Introduction and assumptions.

1.1. *Objective and structure of the paper.* The aim of the following work is to present a probabilistic approach to the homogenization of systems of quasi-linear parabolic PDEs with periodic coefficients. More precisely, for a given arbitrary positive real $T$, we wish to study the asymptotic behavior of the solutions of the following family of equations as $\varepsilon \to 0$:

For $(t,x) \in [0,T[ \times \mathbb{R}^P$ and $\ell \in \{1, \ldots, Q\}$,

$$\frac{\partial(\theta_\varepsilon)_\ell}{\partial t}(t,x) + \frac{1}{2}\sum_{i,j=1}^{P} a_{i,j}(\varepsilon^{-1}x, \theta_\varepsilon(t,x))\frac{\partial^2(\theta_\varepsilon)_\ell}{\partial x_i \, \partial x_j}(t,x)$$

$$\mathcal{E}(\varepsilon) \qquad + \sum_{i=1}^{P}[\varepsilon^{-1}b_i(\varepsilon^{-1}x, \theta_\varepsilon(t,x)) + c_i(\varepsilon^{-1}x, \theta_\varepsilon(t,x), \nabla_x\theta_\varepsilon(t,x))]\frac{\partial(\theta_\varepsilon)_\ell}{\partial x_i}(t,x)$$

$$+ \varepsilon^{-1}e_\ell(\varepsilon^{-1}x, \theta_\varepsilon(t,x)) + f_\ell(\varepsilon^{-1}x, \theta_\varepsilon(t,x), \nabla_x\theta_\varepsilon(t,x)) = 0,$$

$$\theta_\varepsilon(T,x) = H(x),$$

where the coefficients $a$, $b$, $c$, $e$ and $f$ are periodic in $x$ of period one in each direction of the space $\mathbb{R}^P$ (pay attention to the fact that $u$ is $\mathbb{R}^Q$-valued and that $\nabla_x u$ is therefore $\mathbb{R}^{Q \times P}$-valued). Equations such as $(\mathcal{E}(\varepsilon))_{\varepsilon > 0}$ are called quasi-linear since their coefficients depend on the solution or on its gradient (recall that they are said to be "semi-linear" when the differential parts are purely linear, i.e., when $a$, $b$ and $c$ just depend on the variable $x$).

Actually, several results of homogenization of second-order partial differential equations have been already established both by analytical methods [we refer to the monographs of Bensoussan, Lions and Papanicolaou (1978) and Jikov, Kozlov and Oleinik (1994) on the subject] and by probabilistic ones. In particular, Freidlin (1964) and Bensoussan, Lions and Papanicolaou (1978) have established earlier homogenization results in the case of linear equations by means of the well-known Feynman–Kac formula.

The theory of backward stochastic differential equations (in short, BSDEs), whose study is now well known [see, e.g., Pardoux and Peng (1990) and Pardoux (1999a)], provides a probabilistic representation of systems of semi-linear PDEs of second order. Thanks to this deep connection between BSDEs and semi-linear PDEs, two different probabilistic schemes have been developed so far to treat homogenization of such PDEs. On the one hand, Buckdahn, Hu and Peng (1999) base their approach on a stability property of BSDEs. They first deal with the case of classical solutions and then manage from a regularization procedure to weaken the sense given to these solutions. On the other hand, Pardoux (1999a, b) uses weak convergence techniques. He achieves to relax the assumptions made on the coefficients in Buckdahn, Hu and Peng (1999) when the nonlinear term does not depend on the gradient of the solution. Inspired by this latter scheme, Gaudron and Pardoux (2001) consider PDEs whose nonlinearity has a quadratic growth in the gradient, whereas Lejay (2002) deals with the case of PDE operators in divergence form.

As already mentioned, our article treats, again from a probabilistic point of view, the quasi-linear case. It requires work with forward–backward SDEs [see, e.g., Delarue (2002b), Ma, Protter and Yong (1994) and Pardoux and Tang (1999)], which provide a stochastic representation of systems of quasi-linear parabolic PDEs: given a probability space $(\Omega, \mathcal{F}, \mathbf{P})$ and a $P$-dimensional Brownian motion $(B_t)_{t \geq 0}$, we can connect, for every $\varepsilon > 0$, the system $\mathcal{E}(\varepsilon)$ to a family of FBSDEs $(\mathrm{E}(\varepsilon, t, x))_{(t,x) \in [0,T] \times \mathbb{R}^P}$. For every initial condition $(t, x) \in [0, T] \times \mathbb{R}^P$, $\mathrm{E}(\varepsilon, t, x)$ is then given by

$$
\begin{aligned}
&\forall s \in [t, T], \\
&X_s(\varepsilon, t, x) = x + \int_t^s \left(\frac{b}{\varepsilon} + c\right) \left(\frac{X_r(\varepsilon, t, x)}{\varepsilon}, Y_r(\varepsilon, t, x), Z_r(\varepsilon, t, x)\right) dr \\
&\qquad\qquad + \int_t^s \sigma\left(\frac{X_r(\varepsilon, t, x)}{\varepsilon}, Y_r(\varepsilon, t, x)\right) dB_r, \\
&Y_s(\varepsilon, t, x) = H(X_T(\varepsilon, t, x))
\end{aligned}
$$

$\mathrm{E}(\varepsilon, t, x)$

$$+ \int_s^T \left( \frac{e}{\varepsilon} + f \right) \left( \frac{X_r(\varepsilon, t, x)}{\varepsilon}, Y_r(\varepsilon, t, x), Z_r(\varepsilon, t, x) \right) dr$$

$$- \int_s^T Z_r(\varepsilon, t, x) \sigma \left( \frac{X_r(\varepsilon, t, x)}{\varepsilon}, Y_r(\varepsilon, t, x) \right) dB_r,$$

$$\mathbf{E} \int_t^T (|X_s(\varepsilon, t, x)|^2 + |Y_s(\varepsilon, t, x)|^2 + |Z_s(\varepsilon, t, x)|^2) \, ds < \infty,$$

where, for all $\varepsilon > 0$ and $(x, y, z) \in \mathbb{R}^P \times \mathbb{R}^Q \times \mathbb{R}^{Q \times P}$, $(\varepsilon^{-1}b(\mathrm{resp.}\ e) + c(\mathrm{resp.}\ f))(x, y, z) = \varepsilon^{-1}b(\mathrm{resp.}\ e)(x, y) + c(\mathrm{resp.}\ f)(x, y, z)$. Note that, for every $(x, y) \in \mathbb{R}^P \times \mathbb{R}^Q$, the matrix $a(x, y)$ can be assumed to be symmetric: $\sigma$ is then chosen to satisfy, for all $(x, y) \in \mathbb{R}^P \times \mathbb{R}^Q$, $a(x, y) = \sigma\sigma^*(x, y)$. Recall also that the connection between $\mathcal{E}(\varepsilon)$ and $\mathrm{E}(\varepsilon, t, x)$ can be summarized in the following rough way: $\forall s \in [t, T]$, $Y_s(\varepsilon, t, x) = \theta_\varepsilon(s, X_s(\varepsilon, t, x))$.

In order to adapt the stability property approach to this new framework, [Buckdahn and Hu (1998)](#) have restricted themselves to the smooth coefficients case. This permits them to avoid any regularization procedure as done in [Buckdahn, Hu and Peng (1999)](#) and to use crucial estimates of the solutions of the homogenized system. Basically, we manage in this paper to improve this approach and to relax in a very significant way these regularity assumptions. We finally establish the homogenization property despite the weaker regularity of the solution of the limit equation. The reader may find this result quite anecdotal. However, we feel that techniques that we develop to reach our objective permit us to work under assumptions that seem natural when compared to the usual quasi-linear PDEs theory. We guess that they also lead to a better understanding of the FBSDEs machinery and we finally think that our strategy might be applied to similar nonlinear asymptotic problems.

Our approach is largely based on our previous paper, [Delarue (2002b)](#), on FBSDEs. In fact, instead of dealing in a first time with the smooth coefficients case and applying then a regularization procedure as proposed by [Buckdahn, Hu and Peng (1999)](#), we show that it may be more powerful to couple homogenization and regularization in one unique step, as soon as this regularization procedure is efficiently controlled and does not compete with homogenization. Roughly speaking, this strategy can be divided in four steps:

1. We first modify the processes $X(\varepsilon, t, x)$ and $Y(\varepsilon, t, x)$ in order to get rid of the highly oscillating terms $\varepsilon^{-1}b$ and $\varepsilon^{-1}e$. This step is well known in the probabilistic literature devoted to homogenization and leads to the study of the so-called "auxiliary problems."

2. Denoting by $\hat{X}(\varepsilon, t, x)$ and $\hat{Y}(\varepsilon, t, x)$ the modified processes, we then aim to compare the quantity $\theta(\cdot, \hat{X}_\cdot(\varepsilon, t, x))$ with the process $\hat{Y}_\cdot(\varepsilon, t, x)$, where $\theta$ denotes the solution of the presumed limit system. Basically, we aim to write from Itô's formula $\theta(\cdot, \hat{X}.(\varepsilon, t, x))$ as the solution of a BSDE and then to apply a stability property of BSDEs to estimate the distance from this process to $\hat{Y}.(\varepsilon, t, x)$.

3. Since $\theta$ is not regular enough under our assumptions, we introduce a well-chosen regularization sequence of $\theta$, denoted by $(\zeta_n)_{n \in \mathbb{N}^*}$, and then estimate the distance between $\zeta_n(\cdot, \hat{X}.(\varepsilon, t, x))$ and $\hat{Y}.(\varepsilon, t, x)$. Unfortunately, we will see that our assumptions are still too weak to estimate these quantities in a good way. In a nutshell, to apply this strategy, we should be able to control in an appropriate sense the partial derivatives $\nabla^2_{x,x}\zeta_n$ uniformly in $n$, but we are actually not.

4. To face this difficulty, we introduce so-called "auxiliary SDEs" (as far as we know, this point is completely new in the probabilistic literature devoted to homogenization). Basically, we aim to compare $\hat{Y}.(\varepsilon, t, x)$ with a quantity of the form $\zeta_n(\cdot, U.(\varepsilon, n, t, x))$, where $n$ belongs to $\mathbb{N}^*$ and $U(\varepsilon, n, t, x)$ denotes the solution of a well-chosen SDE. In this approach, the distance between the processes $U(\varepsilon, n, t, x)$ and $X(\varepsilon, t, x)$ has to tend toward 0 as $\varepsilon \to 0$ and $n \to +\infty$. We will see in this article how to perform a good choice for these auxiliary SDEs.

From an adapted ergodic property, we finally establish, for every $(t, x) \in [0, T] \times \mathbb{R}^P$, the following kind of convergence:

$$\lim_{\varepsilon \to 0} \Big( \mathbf{E} \sup_{t \le s \le T} |Y_s(\varepsilon, t, x) - \theta(s, X_s(\varepsilon, t, x))|^2$$

$$+ \mathbf{E} \int_t^T \Big| Z_s(\varepsilon, t, x)$$

$$- \nabla_x \theta(s, X_s(\varepsilon, t, x)) \chi_1(\varepsilon^{-1} X_s(\varepsilon, t, x), Y_s(\varepsilon, t, x))$$

$$- \chi_2(\varepsilon^{-1} X_s(\varepsilon, t, x), Y_s(\varepsilon, t, x)) \Big|^2 ds \Big) = 0,$$

where $\chi_1$ and $\chi_2$ denote appropriate corrector terms. It is readily seen from this convergence property that, for every $(t, x) \in [0, T] \times \mathbb{R}^P$, $\theta_\varepsilon(t, x) \to \theta(t, x)$ as $\varepsilon \to 0$. The reader can find a precise statement of the homogenization property (in particular, a precise description of $\chi_1$ and $\chi_2$ and of the limit system) at the end of this section. Note that we also manage to deduce the asymptotic behavior of the distribution of the processes:

$$\Big( X.(\varepsilon, t, x), Y.(\varepsilon, t, x), \int_t^{\cdot} Z_r(\varepsilon, t, x) \sigma(\varepsilon^{-1} X_r(\varepsilon, t, x), Y_r(\varepsilon, t, x)) dB_r,$$

$$\int_t^{\cdot} \varphi(r, \varepsilon^{-1} X_r(\varepsilon, t, x), X_r(\varepsilon, t, x), Y_r(\varepsilon, t, x), Z_r(\varepsilon, t, x)) dr \Big)_{\varepsilon > 0},$$

on the space of continuous functions on $[0, T]$ endowed with the topology of the uniform convergence, where $\varphi$ denotes an arbitrary continuous $\mathbb{R}^d$-valued function, $d \in \mathbb{N}^*$. From our point of view, this result improves the approach proposed by Pardoux (1999a, b).

In fact, the reader can find another proof of the homogenization property established in this paper in our Ph.D. thesis [Delarue (2002a), Chapter 3]. In this latter work, we adapt to our case the compactness techniques introduced by Pardoux (1999a, b) (in particular, we propose a method to handle nonlinearities in gradient in this approach), but we just manage to describe the limit distribution of the processes $(X.(\varepsilon, t, x), Y.(\varepsilon, t, x), \int_t^{\cdot} Z_r(\varepsilon, t, x) \sigma(\varepsilon^{-1} X_r(\varepsilon, t, x), Y_r(\varepsilon, t, x)) \, dB_r)_{\varepsilon > 0}$ in a weaker sense than above: roughly speaking, $Y.(\varepsilon, t, x)$ and $\int_t^{\cdot} Z_r \times \sigma(\varepsilon^{-1} X_r(\varepsilon, t, x), Y_r(\varepsilon, t, x)) \, dB_r$ are considered as elements of the space of càdlàg functions on $[0, T]$ endowed with $S$-topology [$S$-topology is a topology weaker than the Skorohod topology whose compact sets are described in a simple way from the supremum norm and the $\eta$-variations of a càdlàg function. See Jakubowski (1997)]. Nevertheless, to obtain such a result, we also use the concept of "auxiliary SDEs," and therefore, confirm the deep interest of these equations in homogenization of quasi-linear PDEs.

Note finally that we do not know any published analytical proofs of the result established in our paper. In fact, most of the analytical results devoted to homogenization of nonlinear equations deal with divergence form operators: we refer the reader to the monograph of Pankov (1997) for a detailed survey of available results. Nevertheless, in a different framework than ours, several articles deal with the nondivergent case [see, e.g., Avellaneda and Lin (1989) or Evans (1989)].

Our paper is organized as follows. In Section 2 we give several basic preliminary results that are crucial to establish the homogenization property. Most of them are quite well known in the literature devoted to this subject. We detail in Section 3 our strategy to establish the homogenization property. The proof is then given in Section 4. In Appendices A–D, we prove some of the results given in Section 2.

### 1.2. *General notation and assumption.*

**Notation.**

1. For every $N \in \mathbb{N}^*$, let $\langle \cdot, \cdot \rangle$ and $| \cdot |$ denote, respectively, the Euclidean scalar product and the Euclidean norm on $\mathbb{R}^N$; for every $x \in \mathbb{R}^N$, let $x_i$ denote the $i$th coordinate of the vector $x$; and for every $R \geq 0$, let $B_N(x, R)$ [resp. $\bar{B}_N(x, R)$] denote the open (resp. closed) Euclidean ball of dimension $N$, of center $x$ and of radius $R$.

2. For every $N \in \mathbb{N}^*$, let $\rho_N$ denote a nonnegative and smooth function from $\mathbb{R}^N$ into $\mathbb{R}$, null outside $\{u \in \mathbb{R}^N, |u| \leq 1\}$, such that $\int_{\mathbb{R}^N} \rho_N(u) \, du = 1$.

3. Recall that $(\Omega, \mathcal{F}, \mathbf{P})$ denotes a probability space and $(B_t)_{t\geq 0}$ a $P$-dimensional Brownian motion. The usual augmentation of the natural filtration of $(B_t)_{t\geq 0}$ is denoted by $\{\mathcal{F}_t\}_{t\geq 0}$.

4. Let $\mathbb{T}^P$ denote the torus of dimension $P$ (i.e., the quotient space $\mathbb{R}^P/\mathbb{Z}^P$). Then, for every $q \geq 1$, $\|\cdot\|_q$, $\|\cdot\|_{1,q}$ and $\|\cdot\|_{2,q}$ denote, respectively, the $L^q$, $W^{1,q}$ and $W^{2,q}$ norms on $\mathbb{T}^P$. Note that we do not make any differences between a periodic function $g\colon \mathbb{R}^P \to \mathbb{R}^d$ of period one in each direction and the associated function $\tilde{g}\colon \tilde{x} \in \mathbb{T}^P \mapsto g(x)$ with $x \in \tilde{x}$.

5. Recall that $T$ is an arbitrarily prescribed positive real and that

$$b : \mathbb{T}^P \times \mathbb{R}^Q \to \mathbb{R}^P,$$

$$c : \mathbb{T}^P \times \mathbb{R}^Q \times \mathbb{R}^{Q\times P} \to \mathbb{R}^P,$$

$$e : \mathbb{T}^P \times \mathbb{R}^Q \to \mathbb{R}^Q,$$

$$f : \mathbb{T}^P \times \mathbb{R}^Q \times \mathbb{R}^{Q\times P} \to \mathbb{R}^Q,$$

$$\sigma, a : \mathbb{T}^P \times \mathbb{R}^Q \to \mathbb{R}^{P\times P},$$

$$H : \mathbb{R}^P \to \mathbb{R}^Q$$

are measurable functions satisfying, for all $(x,y) \in \mathbb{T}^P \times \mathbb{R}^Q$, $a(x,y) = (\sigma\sigma^*)(x,y)$.

6. Recall finally that for all $q \geq 2$ and $R > 0$ [see also Ladyzenskaja, Solonnikov and Ural'ceva (1968), Chapter I, Section 1],

$$W^{1,2,q}(]0,T[\times B_P(0,R), \mathbb{R}^Q)$$
$$= \{h : ]0,T[\times B_P(0,R) \to \mathbb{R}^Q,$$
$$\quad |h|, |\nabla_t h|, |\nabla_x h|, |\nabla^2_{x,x} h| \in L^q(]0,T[\times B_P(0,R))\}.$$

We then assume that the coefficients $b$, $c$, $e$, $f$, $H$ and $\sigma$ satisfy in the whole paper the following properties:

ASSUMPTION $(\mathcal{H})$.   We say that the functions $b$, $c$, $e$, $f$, $H$ and $\sigma$ satisfy Assumption $(\mathcal{H})$ if there exist three constants $k$, $\lambda > 0$ and $\Lambda$ and an increasing function $K : \mathbb{R}_+ \to \mathbb{R}_+$ such that:

$(\mathcal{H}.1)$. The functions $b$, $c$, $e$, $f$, $H$ and $\sigma$ are continuous on their definition set.

$(\mathcal{H}.2)$. $\forall (x,y) \in \mathbb{R}^P \times \mathbb{R}^Q$, $\forall (x',y') \in \mathbb{R}^P \times \mathbb{R}^Q$, $\forall z \in \mathbb{R}^{Q\times P}$,

$$\langle b(x',y) - b(x,y), x'-x \rangle \leq K(0)|x'-x|^2,$$

$$\langle c(x',y,z) - c(x,y,z), x'-x \rangle \leq K(|y|+|z|)|x'-x|^2,$$

$$\langle f(x,y',z) - f(x,y,z), y'-y \rangle \leq K(|y|+|y'|+|z|)|y'-y|^2.$$

(ℋ.3). $\forall (x,y,z) \in \mathbb{R}^P \times \mathbb{R}^Q \times \mathbb{R}^{Q \times P}$, $\forall (x',y',z') \in \mathbb{R}^P \times \mathbb{R}^Q \times \mathbb{R}^{Q \times P}$,

$$|\sigma(x',y') - \sigma(x,y)| + |H(x') - H(x)| \le k(|x'-x| + |y'-y|),$$

$$|b(x,y') - b(x,y)| + |e(x',y') - e(x,y)| \le K(0)(|x'-x| + |y'-y|),$$

$$|f(x',y,z') - f(x,y,z)| \le K(|y| + |z|)$$
$$\times (|x'-x| + |z'-z|),$$

$$|c(x,y',z') - c(x,y,z)| \le K(|y| + |y'| + |z|)$$
$$\times (|y'-y| + |z'-z|).$$

(ℋ.4). $\forall (x,y,z) \in \mathbb{R}^P \times \mathbb{R}^Q$,

$$|b(x,y)| + |e(x,y)| + |\sigma(x,y)| + |H(x)| \le \Lambda,$$

$$|c(x,y,z)| + |f(x,y,z)| \le \Lambda(1 + |y| + |z|).$$

(ℋ.5). $\forall (x,y) \in \mathbb{R}^P \times \mathbb{R}^Q$, $\forall \xi \in \mathbb{R}^P$, $\langle \xi, a(x,y)\xi \rangle \ge \lambda |\xi|^2$.

(ℋ.6). $\sigma \in \mathcal{C}^{0,2}(\mathbb{T}^P \times \mathbb{R}^Q, \mathbb{R}^{P \times P})$, $b \in \mathcal{C}^{0,2}(\mathbb{T}^P \times \mathbb{R}^Q, \mathbb{R}^P)$ and $e \in \mathcal{C}^{0,2}(\mathbb{T}^P \times \mathbb{R}^Q, \mathbb{R}^Q)$.

(ℋ.7). $\forall (x,y) \in \mathbb{R}^P \times \mathbb{R}^Q$,

$$|\nabla_y b(x,y)| + |\nabla_{y,y}^2 b(x,y)| + |\nabla_y e(x,y)|$$
$$+ |\nabla_{y,y}^2 e(x,y)| + |\nabla_y \sigma(x,y)| + |\nabla_{y,y}^2 \sigma(x,y)| \le K(0).$$

1.3. *Precise statement of the homogenization property.* Here is the main result of the paper:

Theorem (HP). *Suppose that Assumption (ℋ) is in force. Then, set*

$$\mathcal{V} = \Bigg\{ h \in \mathcal{C}([0,T] \times \mathbb{R}^P, \mathbb{R}^Q) \cap \bigcap_{q \ge 2} \bigcap_{R > 0} W^{1,2,q}(]0,T[ \times B_P(0,R), \mathbb{R}^Q),$$

$$\sup_{t \in [0,T]} \Bigg[ \sup_{x \in \mathbb{R}^P} |h(t,x)| + \sup_{x,x' \in \mathbb{R}^P, x \ne x'} \frac{|h(t,x') - h(t,x)|}{|x' - x|} \Bigg] < \infty \Bigg\}.$$

1. *For every $\varepsilon > 0$, the system of PDEs $\mathcal{E}(\varepsilon)$ admits a unique solution in the space $\mathcal{V}$. It is denoted by $\theta_\varepsilon$.*

2. *For every $y \in \mathbb{R}^Q$, there exists a unique density on the torus, denoted by $p(\cdot, y)$, belonging to the space $W^{1,2}(\mathbb{T}^P)$ and satisfying the following equation on $\mathbb{T}^P$:*

Inv(y) $$L_y^*(p(\cdot,y)) = 0, \qquad \int_{\mathbb{T}^P} p(x,y)\,dx = 1,$$

where

$\mathcal{G}$en$(y)$         $$L_y = \frac{1}{2} \sum_{i,j=1}^{P} a_{i,j}(\cdot, y) \frac{\partial^2}{\partial x_i \partial x_j} + \sum_{i=1}^{P} b_i(\cdot, y) \frac{\partial}{\partial x_i}.$$

3. *If, for every $y \in \mathbb{R}^Q$, the functions $(b_i(\cdot, y))_{1 \le i \le P}$ and $(e_j(\cdot, y))_{1 \le j \le Q}$ satisfy*

($\mathcal{H}.8$).

$$\forall i \in \{1, \ldots, P\}, \qquad \int_{\mathbb{T}^P} b_i(x, y) p(x, y) \, dx = 0,$$

$$\forall j \in \{1, \ldots, Q\}, \qquad \int_{\mathbb{T}^P} e_j(x, y) p(x, y) \, dx = 0,$$

*then, $\forall (t, x) \in [0, T] \times \mathbb{R}^P$, $\lim_{\varepsilon \to 0} \theta_\varepsilon(t, x) = \theta(t, x)$, where $\theta$ denotes the unique solution in the space $\mathcal{V}$ of the system of PDEs:*

*For $(t, x) \in [0, T[ \times \mathbb{R}^P$ and $\ell \in \{1, \ldots, Q\}$,*

$\mathcal{E}$(lim)

$$\frac{\partial \theta_\ell}{\partial t}(t, x) + \frac{1}{2} \sum_{i,j=1}^{P} \bar{\alpha}_{i,j}(\theta(t, x)) \frac{\partial^2 \theta_\ell}{\partial x_i \partial x_j}(t, x)$$

$$+ \sum_{i=1}^{P} \bar{u}_i(\theta(t, x), \nabla_x \theta(t, x)) \frac{\partial \theta_\ell}{\partial x_i}(t, x)$$

$$+ \bar{v}_\ell(\theta(t, x), \nabla_x \theta(t, x)) = 0,$$

$$\theta(T, x) = H(x).$$

4. *The so-called "homogenized" coefficients $\bar{\alpha}$, $\bar{u}$ and $\bar{v}$ are given by the relationships*

$\forall (y, z) \in \mathbb{R}^Q \times \mathbb{R}^{Q \times P}, \qquad \bar{u}(y, z) = \int_{\mathbb{T}^P} u(x, y, z(I + \nabla_x \hat{b})(x, y)) p(x, y) \, dx,$

$$\bar{v}(y, z) = \int_{\mathbb{T}^P} v(x, y, z(I + \nabla_x \hat{b})(x, y)) p(x, y) \, dx,$$

$$\bar{\alpha}(y) = \int_{\mathbb{T}^P} \alpha(x, y) p(x, y) \, dx,$$

*where for all $(x, y, z) \in \mathbb{T}^P \times \mathbb{R}^Q \times \mathbb{R}^{Q \times P}$:*

$$u(x, y, z) = (I + \nabla_x \hat{b})(x, y) c(x, y, z + \nabla_x \hat{e}(x, y)) - (\nabla_y \hat{b} e)(x, y)$$

$$+ \nabla_{x,y}^2 \hat{b}(x, y) [a(x, y)(z + \nabla_x \hat{e}(x, y))^*],$$

$$v(x, y, z) = f(x, y, z + \nabla_x \hat{e}(x, y))$$

$$+ \nabla_x \hat{e}(x, y) c(x, y, z + \nabla_x \hat{e}(x, y)) - (\nabla_y \hat{e} e)(x, y)$$

$$+ \nabla_{x,y}^2 \hat{e}(x, y) [a(x, y)(z + \nabla_x \hat{e}(x, y))^*],$$

$$\alpha(x, y) = (I + \nabla_x \hat{b})(x, y) a(x, y)(I + \nabla_x \hat{b})^*(x, y),$$

*with the convention* $\nabla_{x,y}\hat{b}(x,y)z^* = (\sum_{i,j} \frac{\partial \hat{b}_\ell}{\partial x_i \partial y_j} z_{j,i})_{1\leq \ell \leq P}$, $\nabla_{x,y}\hat{e}(x,y)z^* = (\sum_{i,j} \frac{\partial \hat{e}_\ell}{\partial x_i \partial y_j} z_{j,i})_{1\leq \ell \leq Q}$.

The functions $\hat{b}$ and $\hat{e}$ appearing in the definitions of these coefficients are given as the solutions of an auxiliary family of PDEs. More precisely, for every $y \in \mathbb{R}^Q$, the functions $\hat{b}(\cdot,y)$ and $\hat{e}(\cdot,y)$ are the unique solutions in the spaces $(\bigcap_{q\geq 2} W^{2,q}(\mathbb{T}^P))^P$ and $(\bigcap_{q\geq 2} W^{2,q}(\mathbb{T}^P))^Q$ of the systems of equations on $\mathbb{T}^P$:

$$\forall i \in \{1,\ldots,P\}, \qquad L_y \hat{b}_i(\cdot,y) + b_i(\cdot,y) = 0 \quad \text{and}$$

$$\int_{\mathbb{T}^P} \hat{b}_i(x,y)\,dx = 0,$$

$\mathcal{A}\mathrm{ux}(y)$

$$\forall j \in \{1,\ldots,Q\}, \qquad L_y \hat{e}_j(\cdot,y) + e_j(\cdot,y) = 0 \quad \text{and}$$

$$\int_{\mathbb{T}^P} \hat{e}_j(x,y)\,dx = 0.$$

Note that the existence of the integrals and derivatives that may appear in the definitions of the homogenized coefficients will be justified in the paper.

We will also establish the probabilistic counterpart of the latter analytical result:

**Theorem** (PHP). *Suppose that Assumption ($\mathcal{H}$) is in force and that the additional hypothesis ($\mathcal{H}.8$) holds. Then:*

1. *For every $\varepsilon > 0$, for every $(t,x) \in [0,T] \times \mathbb{R}^P$, there exists a unique progressively measurable solution to the FBSDE $\mathrm{E}(\varepsilon,t,x)$, which is denoted by $(X_s(\varepsilon,t,x),Y_s(\varepsilon,t,x),Z_s(\varepsilon,t,x))_{t\leq s\leq T}$.*

2. *Let $(t,x) \in [0,T] \times \mathbb{R}^P$ and set for all $\varepsilon > 0$ and $s \in [t,T]$, $\hat{Z}_s(\varepsilon,t,x) = Z_s(\varepsilon,t,x) - \nabla_x \hat{e}(\varepsilon^{-1} X_s(\varepsilon,t,x), Y_s(\varepsilon,t,x))$. Then, the following convergence holds as $\varepsilon \to 0$:*

$$\lim_{\varepsilon \to 0} \mathbf{E} \sup_{t\leq s\leq T} |Y_s(\varepsilon,t,x) - \theta(s,X_s(\varepsilon,t,x))|^2$$

$$+ \mathbf{E} \int_t^T |\hat{Z}_s(\varepsilon,t,x)$$

$$- \nabla_x \theta(s,X_s(\varepsilon,t,x))(I + \nabla_x \hat{b})(\varepsilon^{-1} X_s(\varepsilon,t,x), Y_s(\varepsilon,t,x))|^2 \, ds$$

$$= 0.$$

3. *For every $(t,x) \in [0,T] \times \mathbb{R}^P$, we associate to the system $\mathcal{E}(\mathrm{lim})$ and to the initial condition $(t,x)$ the following FBSDE:*

$$\forall s \in [t,T],$$

$$\text{E}(t,x) \qquad \begin{aligned} X_s &= x + \int_t^s \bar{u}(Y_r, Z_r)\, dr + \int_t^s \bar{\alpha}^{1/2}(Y_r)\, dB_r, \\[2mm] Y_s &= H(X_T) + \int_s^T \bar{v}(Y_r, Z_r)\, dr - \int_s^T Z_r \bar{\alpha}^{1/2}(Y_r)\, dB_r, \\[2mm] &\text{E}\int_t^T (|X_r|^2 + |Y_r|^2 + |Z_r|^2)\, dr < \infty. \end{aligned}$$

It admits a unique progressively measurable solution denoted by $(X_s(t,x), Y_s(t,x), Z_s(t,x))_{t \le s \le T}$. Then, for a given bounded function, $\varphi \colon [0,T] \times \mathbb{T}^P \times \mathbb{R}^P \times \mathbb{R}^Q \times \mathbb{R}^{Q \times P} \to \mathbb{R}^N$ satisfying

$$\forall R > 0, \ \forall u \in \mathbb{T}^P,$$

$$\forall (t,x,y,z), (t',x',y',z') \in [0,T] \times \bar{B}_P(0,R) \times \bar{B}_Q(0,R) \times \bar{B}_{Q \times P}(0,R),$$

$$|\varphi(t',u,x',y',z') - \varphi(t,u,x,y,z)|$$

$$\le w_R(|t'-t| + |x'-x| + |y'-y| + |z'-z|),$$

where, for every $R > 0$, $\lim_{\eta \to 0} w_R(\eta) = 0$, we claim that for every $(t,x) \in [0,T] \times \mathbb{R}^P$,

$$\left( X_s(\varepsilon,t,x), Y_s(\varepsilon,t,x), \int_t^s \hat{Z}_r(\varepsilon,t,x) \sigma(\varepsilon^{-1} X_r(\varepsilon,t,x), Y_r(\varepsilon,t,x))\, dB_r, \right.$$

$$\left. \int_t^s \varphi(r, \varepsilon^{-1} X_r(\varepsilon,t,x), X_r(\varepsilon,t,x), Y_r(\varepsilon,t,x), \hat{Z}_r(\varepsilon,t,x))\, dr \right)_{t \le s \le T}$$

$$\implies \left( X_s(t,x), Y_s(t,x), \int_t^s Z_r(t,x) \bar{\alpha}^{1/2}(Y_r(t,x))\, dB_r, \right.$$

$$\left. \int_t^s \bar{\varphi}(r, X_r(t,x), Y_r(t,x), Z_r(t,x))\, dr \right)_{t \le s \le T},$$

where $\implies$ stands for the convergence in law on the space $\mathcal{C}([t,T], \mathbb{R}^P \times \mathbb{R}^Q \times \mathbb{R}^Q \times \mathbb{R}^N)$ endowed with the topology of the uniform convergence and where we have set

$$\forall (t,x,y,z) \in [0,T] \times \mathbb{R}^P \times \mathbb{R}^Q \times \mathbb{R}^{Q \times P},$$

$$\bar{\varphi}(t,x,y,z) = \int_{\mathbb{T}^P} \varphi(t,r,x,y,z(I + \nabla_x \hat{b}(r,y))) p(r,y)\, dr.$$

**2. Preliminary results.** This section aims to present basic materials needed to solve the homogenization problem. As most of the tools introduced in the following lines are quite well known, we have decided to summarize the main results without detailing their proofs. Nevertheless, for the sake of completeness, the reader can find the demonstrations (or at least sketches of them) at the end of the paper.

This section is organized in the following way:

1. We first give basic solvability results of equations $(\mathcal{E}(\varepsilon))_{\varepsilon>0}$ and $(\mathrm{E}(\varepsilon, t, x))_{\varepsilon>0, t\in[0,T], x\in\mathbb{R}^P}$. We then show that it is sufficient to investigate the convergence, as $\varepsilon$ tends toward 0, of the processes $(X(\varepsilon, 0, x_0), Y(\varepsilon, 0, x_0), Z(\varepsilon, 0, x_0))$ for a given $x_0 \in \mathbb{R}^P$.

2. We solve the so-called "auxiliary problems" and give useful estimates of the associated solutions.

3. In order to get rid of the highly oscillating terms, we modify the forward and backward processes involved in the representation of the nonlinear systems of PDEs. As a by-product, we obtain basic estimates of these processes.

4. Thanks to the estimates of the solutions of the auxiliary problems, we successfully adapt the usual "ergodic theorem" to our nonlinear framework.

5. We finally give estimates of the solution of the limit system that are crucial to establish the homogenization property.

2.1. *Solvability results and choice of the initial condition of $X(\varepsilon)$.* According to Delarue (2002b), we know that, for all $\varepsilon > 0$ and $(t, x) \in [0, T] \times \mathbb{R}^P$, the FBSDE $\mathrm{E}(\varepsilon, t, x)$ admits a unique $\{\mathcal{F}_s\}_{t \le s \le T}$-progressively measurable solution. It is denoted by $(X_s(\varepsilon, t, x), Y_s(\varepsilon, t, x), Z_s(\varepsilon, t, x))_{t \le s \le T}$. We then define

(2.1) $$\theta_\varepsilon : (t, x) \in [0, T] \times \mathbb{R}^P \mapsto Y_t(\varepsilon, t, x).$$

According to Delarue (2002b), $\theta_\varepsilon$ is bounded, 1/2-Hölderian in $t$ (uniformly in $x$) and Lipschtizan in $x$ (uniformly in $t$). Moreover, for every $(t, x) \in [0, T] \times \mathbb{R}^P$,

(2.2) $$\forall s \in [t, T], \qquad Y_s(\varepsilon, t, x) = \theta_\varepsilon(s, X_s(\varepsilon, t, x)).$$

In fact, we will prove in Appendix D that, for every $\varepsilon > 0$, $\theta_\varepsilon$ is the unique solution in the space $\mathcal{V}$ of the system of PDEs $\mathcal{E}(\varepsilon)$. In particular, relationships (2.1) and (2.2) will be crucial to pass from Theorem (PHP) to Theorem (HP).

Since we aim to establish the pointwise convergence of the functions $(\theta_\varepsilon)_{\varepsilon>0}$ as $\varepsilon \to 0$, note that we just have to prove the convergence of the sequence $(\theta_\varepsilon(0, x_0))_{\varepsilon>0}$ for a given $x_0 \in \mathbb{R}^P$.

Hence, our strategy is clear: fix once for all $x_0$ in $\mathbb{R}^P$ and denote for the sake of simplicity, for every $\varepsilon > 0$, the process $(X_s(\varepsilon, 0, x_0), Y_s(\varepsilon, 0, x_0), Z_s(\varepsilon, 0, x_0))_{0 \le s \le T}$ by $(X_s(\varepsilon), Y_s(\varepsilon), Z_s(\varepsilon))_{0 \le s \le T}$. We then aim to establish, in the sequel of the paper, the convergence, in a suitable sense, of the family $(X(\varepsilon), Y(\varepsilon), Z(\varepsilon))_{\varepsilon>0}$ as $\varepsilon$ tends toward 0.

2.2. *Auxiliary problems.* In this section we state the main solvability results of $(\mathfrak{Inv}(y))_{y \in \mathbb{R}^Q}$ and $(\mathcal{A}\mathrm{ux}(y))_{y \in \mathbb{R}^Q}$ that we establish in Appendix A. We also detail the regularity properties of the solutions with respect to $y$.

THEOREM 2.1 (Estimates of the invariant measures).

1. *For every $y \in \mathbb{R}^Q$, the $\mathbb{T}^P$-valued Markov process associated to the operator $L_y$ admits a unique invariant probability measure, which is denoted by $m(\cdot, y)$. $m(\cdot, y)$ is absolutely continuous with respect to the Lebesgue measure on the torus $\mathbb{T}^P$, and its density, denoted by $p(\cdot, y)$, is the unique solution of $\mathfrak{Inv}(y)$ in $W^{1,2}(\mathbb{T}^P)$.*

    *Moreover, there exists a constant $C_{2.1} > 0$, only depending on $k$, $\lambda$, $\Lambda$ and $P$, such that*

$$\forall (x, y) \in \mathbb{T}^P \times \mathbb{R}^Q, \qquad \frac{1}{C_{2.1}} \le p(x, y) \le C_{2.1}.$$

2. *The function $p : y \in \mathbb{R}^Q \mapsto p(\cdot, y) \in L^2(\mathbb{T}^P)$ is twice continuously differentiable [the partial derivatives are denoted by $(\frac{\partial p}{\partial y_i})_{1 \le i \le Q}$ and $(\frac{\partial^2 p}{\partial y_i \partial y_j})_{1 \le i,j \le Q}$]. Moreover, there exists a constant $C_{2.2}$, only depending on $k$, $K$, $\lambda$, $\Lambda$, $P$ and $Q$, such that*

$$\forall (i, j) \in \{1, \dots, Q\}^2, \ \forall y \in \mathbb{R}^Q, \qquad \left\| \frac{\partial p}{\partial y_i}(\cdot, y) \right\|_2 + \left\| \frac{\partial^2 p}{\partial y_i \partial y_j}(\cdot, y) \right\|_2 \le C_{2.2}.$$

THEOREM 2.2 (Auxiliary problems). *Suppose that the functions $b$ and $e$ satisfy the additional assumption $(\mathcal{H}.8)$. Then:*

1. *For every $y \in \mathbb{R}^Q$, we can define the vector valued functions $\hat{b}(\cdot, y)$ and $\hat{e}(\cdot, y)$ as the unique solutions of the equations $\mathcal{A}\mathrm{ux}(y)$ in the spaces $(\bigcap_{q \ge 2} W^{2,q}(\mathbb{T}^P))^P$ and $(\bigcap_{q \ge 2} W^{2,q}(\mathbb{T}^P))^Q$.*

2. *For every $q \ge 2$, the functions $y \in \mathbb{R}^Q \mapsto \hat{b}(\cdot, y) \in (W^{2,q}(\mathbb{T}^P))^P$ and $y \in \mathbb{R}^Q \mapsto \hat{e}(\cdot, y) \in (W^{2,q}(\mathbb{T}^P))^Q$ are twice continuously differentiable. Moreover, there exists a constant $C_{2.3}^{(q)}$, only depending on $k$, $K$, $\lambda$, $\Lambda$, $P$, $q$ and $Q$, such that for every $y \in \mathbb{R}^Q$,*

$$\|\hat{b}(\cdot, y)\|_{2,q} + \|\hat{e}(\cdot, y)\|_{2,q} + \|\nabla_y \hat{b}(\cdot, y)\|_{2,q}$$
$$+ \|\nabla_y \hat{e}(\cdot, y)\|_{2,q} + \|\nabla_{y,y}^2 \hat{b}(\cdot, y)\|_{2,q} + \|\nabla_{y,y}^2 \hat{e}(\cdot, y)\|_{2,q} \le C_{2.3}^{(q)}.$$

3. *In particular, from the Sobolev imbedding theorems, the functions $(x, y) \in \mathbb{T}^P \times \mathbb{R}^Q \mapsto \hat{b}(x, y) \in \mathbb{R}^P$ and $(x, y) \in \mathbb{T}^P \times \mathbb{R}^Q \mapsto \hat{e}(x, y) \in \mathbb{R}^Q$ are continuously differentiable with respect to $x$ and twice continuously differentiable with respect to $y$ [note that, for every $q \ge 2$, the usual partial derivatives in $y$ of $\hat{b}$ and $\hat{e}$ coincide with the derivatives in $y$ of $\hat{b}$ and $\hat{e}$ seen as functions with values in $(W^{2,q}(\mathbb{T}^P))^P$ and $(W^{2,q}(\mathbb{T}^P))^Q$]. In the same way,* $(x, y) \in \mathbb{T}^P \mapsto \nabla_x \hat{b}(x, y) \in \mathbb{R}^{P \times P}$ and $(x, y) \in \mathbb{T}^P \mapsto \nabla_x \hat{e}(x, y) \in \mathbb{R}^{Q \times P}$ are twice continuously differentiable with respect to $y$. All the associated derivatives are uniformly bounded on $\mathbb{T}^P \times \mathbb{R}^Q$.

2.3. *Modification of the processes* $X(\varepsilon)$, $Y(\varepsilon)$ *and* $Z(\varepsilon)$. Suppose that the additional assumption $(\mathcal{H}.8)$ holds. Then, thanks to the definitions of $\hat{b}$ and $\hat{e}$, we are able to get rid of the highly oscillating terms $\varepsilon^{-1} b(\varepsilon^{-1} X(\varepsilon), Y(\varepsilon))$ and $\varepsilon^{-1} e(\varepsilon^{-1} X(\varepsilon), Y(\varepsilon))$ by setting for all $\varepsilon > 0$ and $t \in [0, T]$,

$$
\begin{aligned}
\bar{X}_t(\varepsilon) &= \varepsilon^{-1} X_t(\varepsilon), \\
\hat{X}_t(\varepsilon) &= X_t(\varepsilon) + \varepsilon \hat{b}(\bar{X}_t(\varepsilon), Y_t(\varepsilon)), \\
\hat{Y}_t(\varepsilon) &= Y_t(\varepsilon) - \varepsilon \hat{e}(\bar{X}_t(\varepsilon), Y_t(\varepsilon)), \\
\hat{Z}_t(\varepsilon) &= Z_t(\varepsilon) - \nabla_x \hat{e}(\bar{X}_t(\varepsilon), Y_t(\varepsilon)).
\end{aligned}
\tag{2.3}
$$

Referring to the notation introduced in point 4 of Theorem (HP) and according to the Itô–Krylov formula [see Krylov (1980), Theorem 1, Paragraph 10, Chapter II, as well as Pardoux and Veretennikov (2002)] and to Theorem 2.2, we then have that for all $\varepsilon > 0$ and $t \in [0, T]$,

$$
\begin{aligned}
\hat{X}_t(\varepsilon) = {}& \hat{X}_0(\varepsilon) + \int_0^t u(\bar{X}_s(\varepsilon), Y_s(\varepsilon), \hat{Z}_s(\varepsilon)) \, ds \\
& - \varepsilon \int_0^t (\nabla_y \hat{b} f)(\bar{X}_s(\varepsilon), Y_s(\varepsilon), Z_s(\varepsilon)) \, ds \\
& + \tfrac{1}{2} \varepsilon \int_0^t \nabla_{y,y}^2 \hat{b}(\bar{X}_s(\varepsilon), Y_s(\varepsilon)) [Z_s(\varepsilon) a(\bar{X}_s(\varepsilon), Y_s(\varepsilon)) Z_s^*(\varepsilon)] \, ds \\
& + \int_0^t [(I + \nabla_x \hat{b}) \sigma](\bar{X}_s(\varepsilon), Y_s(\varepsilon)) \, dB_s \\
& + \varepsilon \int_0^t \nabla_y \hat{b}(\bar{X}_s(\varepsilon), Y_s(\varepsilon)) Z_s(\varepsilon) \sigma(\bar{X}_s(\varepsilon), Y_s(\varepsilon)) \, dB_s,
\end{aligned}
\tag{2.4}
$$

$$
\begin{aligned}
\hat{Y}_t(\varepsilon) = {}& \hat{Y}_T(\varepsilon) + \int_t^T v(\bar{X}_s(\varepsilon), Y_s(\varepsilon), \hat{Z}_s(\varepsilon)) \, ds \\
& - \varepsilon \int_t^T (\nabla_y \hat{e} f)(\bar{X}_s(\varepsilon), Y_s(\varepsilon), Z_s(\varepsilon)) \, ds \\
& + \frac{1}{2} \varepsilon \int_t^T \nabla_{y,y}^2 \hat{e}(\bar{X}_s(\varepsilon), Y_s(\varepsilon)) [Z_s(\varepsilon) a(\bar{X}_s(\varepsilon), Y_s(\varepsilon)) Z_s^*(\varepsilon)] \, ds \\
& - \int_t^T \hat{Z}_s(\varepsilon) \sigma(\bar{X}_s(\varepsilon), Y_s(\varepsilon)) \, dB_s \\
& + \varepsilon \int_t^T \nabla_y \hat{e}(\bar{X}_s(\varepsilon), Y_s(\varepsilon)) Z_s(\varepsilon) \sigma(\bar{X}_s(\varepsilon), Y_s(\varepsilon)) \, dB_s.
\end{aligned}
\tag{2.5}
$$

Using Itô's calculus, we can establish basic estimates of processes $((X_t(\varepsilon), Y_t(\varepsilon), Z_t(\varepsilon))_{0 \leq t \leq T})_{\varepsilon > 0}$. In particular, the following theorem is proved in Appendix B.

THEOREM 2.3 (Estimates of representation processes). *There exists a constant $C_{2.4}$, only depending on $k$, $K$, $\lambda$, $\Lambda$, $P$, $Q$ and $T$, such that*

$$\forall \varepsilon > 0, \qquad \mathbf{P}\left\{\sup_{t \in [0,T]} |Y_t(\varepsilon)| \leq C_{2.4}\right\} = 1,$$

$$\mathbf{E} \sup_{t \in [0,T]} |X_t(\varepsilon)|^2 + \mathbf{E}\left(\int_0^T |Z_s(\varepsilon)|^2 \, ds\right)^2 \leq C_{2.4}.$$

2.4. *Ergodic theorem.* Suppose that the additional assumption $(\mathcal{H}.8)$ holds and consider, for a given $d \in \mathbb{N}^*$, the following family of $\mathbb{R}^d$-valued semimartingales:

(2.6)
$$\forall \varepsilon > 0, \ \forall t \in [0,T],$$
$$G_t(\varepsilon) = G_0(\varepsilon) + \int_0^t g_s(1,\varepsilon) \, ds + \int_0^t g_s(2,\varepsilon) \, dB_s,$$

where, for every $\varepsilon > 0$, $G_0(\varepsilon)$ is $\mathcal{F}_0$-measurable and $g(1,\varepsilon)$ and $g(2,\varepsilon)$ are two progressively measurable processes satisfying

(2.7)
$$\sup_{0 < \varepsilon < 1} \left(\mathbf{E}|G_0(\varepsilon)| + \mathbf{E}\int_0^T (|g_t(1,\varepsilon)| + |g_t(2,\varepsilon)|^2) \, dt\right) < \infty.$$

Then, we will establish in Appendix C the following ergodic property.

THEOREM 2.4 (Ergodic theorem).

1. *Let $\varphi : [0,T] \times \mathbb{T}^P \times \mathbb{R}^Q \times \mathbb{R}^d \to \mathbb{R}$ be a bounded and measurable function, such that for every compact set $\kappa \subset [0,T] \times \mathbb{R}^Q \times \mathbb{R}^d$, the family of functions $((t,y,g) \in \kappa \mapsto \varphi(t,x,y,g))_{x \in \mathbb{T}^P}$ is equicontinuous. Then*

$$\lim_{\varepsilon \to 0} \mathbf{E} \sup_{0 \leq t \leq T} \left| \int_0^t [\varphi(r, \varepsilon^{-1} X_r(\varepsilon), Y_r(\varepsilon), G_r(\varepsilon)) - \bar{\varphi}(r, Y_r(\varepsilon), G_r(\varepsilon))] \, dr \right| = 0,$$

*with*

$$\forall (t,y,g) \in [0,T] \times \mathbb{R}^Q \times \mathbb{R}^d, \qquad \bar{\varphi}(t,y,g) = \int_{\mathbb{T}^P} \varphi(t,x,y,g) p(x,y) \, dx.$$

2. *Let $\varphi : [0,T] \times \mathbb{T}^P \times \mathbb{R}^Q \times \mathbb{R}^d \to \mathbb{R}$ satisfy*

$$\forall R > 1, \ \forall (t,t') \in [0, (1-R^{-1})T]^2,$$
$$\forall (y,y') \in (\bar{B}_Q(0,R))^2, \forall (g,g') \in (\bar{B}_d(0,R))^2:$$
$$|\varphi(t',y',g') - \varphi(t,y,g)| \leq \nu_R(|t'-t| + |y'-y| + |g'-g|),$$

where, for every $R > 0$, $\lim_{\delta \to 0} \nu_R(\delta) = 0$. If, for every $0 \leq \delta \leq T$,

$$\mathbf{E} \int_{T-\delta}^{T} |\varphi(s, \bar{X}_s(\varepsilon), Y_s(\varepsilon), G_s(\varepsilon)) - \bar{\varphi}(s, Y_s(\varepsilon), G_s(\varepsilon))| \, ds \leq w(\delta),$$

with $\lim_{\delta \to 0} w(\delta) = 0$, then

$$\forall \eta > 0, \qquad \lim_{\varepsilon \to 0} \mathbf{P} \left\{ \sup_{0 \leq t \leq T} \left| \int_0^t [\varphi(s, \bar{X}_s(\varepsilon), Y_s(\varepsilon), G_s(\varepsilon)) \right.\right.$$
$$\left.\left. - \bar{\varphi}(s, Y_s(\varepsilon), G_s(\varepsilon))] \, ds \right| \geq \eta \right\} = 0.$$

2.5. *Limit system.* Suppose that the additional assumption $(\mathcal{H}.8)$ holds. We then claim that the limit system of PDEs $\mathcal{E}(\lim)$ is uniquely solvable in the space $\mathcal{V}$, and that the solution satisfies relevant estimates.

We first note that limit coefficients, as defined in Section 1.3, satisfy usual assumptions: there exist three constants $\bar{k}$, $\bar{\lambda}$ and $\bar{\Lambda}$ and a nondecreasing function $\bar{K}$ from $\mathbb{R}_+$ into itself, only depending on $\lambda$, $\Lambda$, $k$, $K$, $P$ and $Q$, such that

(2.5.i) *For every $x \in \mathbb{R}^P$, $u(x, \cdot, \cdot)$ and $v(x, \cdot, \cdot)$ satisfy with respect to $\bar{\Lambda}$ and $\bar{K}$ the same properties as $c$ and $f$ in $(\mathcal{H}.1)$–$(\mathcal{H}.4)$.*

(2.5.ii) $\bar{u}$, $\bar{v}$, $H$ and $\bar{\alpha}^{1/2}$ *satisfy $(\mathcal{H}.1)$–$(\mathcal{H}.5)$ with respect to $\bar{k}$, $\bar{\lambda}$, $\bar{\Lambda}$ and $\bar{K}$ ($\bar{u}$, $\bar{v}$ and $\bar{\alpha}^{1/2}$ playing the role of $c$, $f$ and $\sigma$).*

PROOF OF (2.5.i) AND (2.5.ii). Let us first prove nondegeneracy of $\bar{\alpha}$. For $(\xi, y) \in \mathbb{R}^P \times \mathbb{R}^Q$,

$$(2.8) \qquad \langle \xi, \bar{\alpha}(y)\xi \rangle \geq \frac{\lambda}{C_{2.1}} \int_{\mathbb{T}^P} \langle \xi, (I + \nabla_x \hat{b})(I + \nabla_x \hat{b})^*(x, y)\xi \rangle \, dx.$$

From the relation

$$(2.9) \qquad \int_{\mathbb{T}^P} \langle \xi, (\nabla_x \hat{b} + (\nabla_x \hat{b})^*)(x, y)\xi \rangle \, dx = 0,$$

we deduce that $\bar{\alpha}$ satisfies $(\mathcal{H}.5)$ with respect to $\lambda/C_{2.1}$.

The other properties are easily proved from Assumption $(\mathcal{H})$ and Theorems 2.1 and 2.2. $\square$

Referring to Delarue (2002b), we then claim that, for every $(t, x) \in [0, T] \times \mathbb{R}^P$, the limit FBSDE $\mathrm{E}(t, x)$ admits a unique solution, denoted by $(X(t, x), Y(t, x), Z(t, x))$. Hence, we can define

$$(2.10) \qquad \theta : (t, x) \in [0, T] \times \mathbb{R}^P \mapsto Y_t(t, x).$$

Thanks to Delarue (2002b) [see also Ladyzenskaja, Solonnikov and Ural'ceva (1968) and Delarue (2003)], there exist a constant $C_{2.5}$, only depending on $\bar{\Lambda}$ and $T$, and a constant $C_{2.6}$, only depending on $\bar{k}$, $\bar{K}$, $\bar{\lambda}$, $\bar{\Lambda}$, $P$, $Q$ and $T$, such that

$$(2.11) \qquad \forall\, (t,x) \in [0,T] \times \mathbb{R}^P, \qquad |\theta(t,x)| \leq C_{2.5},$$

$$\forall\, (t,t') \in [0,T]^2,\ \forall\, (x,x') \in (\mathbb{R}^P)^2,$$
$$(2.12)$$
$$|\theta(t',x') - \theta(t,x)| \leq C_{2.6}(|t'-t|^{1/2} + |x'-x|).$$

Actually, we will prove in Appendix D that $\theta$ is the unique solution of the limit system of PDEs $\mathcal{E}(\lim)$ in the space $\mathcal{V}$. Moreover, we will also show (see Appendix D) that $\theta$ is continuously differentiable on $[0,T[\,\times\mathbb{R}^P$ with respect to $x$, and that there exists a constant $0 < \beta \leq 1$, only depending on $\bar{\lambda}$, $\bar{\Lambda}$ and $P$, such that for every $0 < \eta < T$,

$$(2.13) \qquad \|\nabla_x \theta\|_{\mathcal{C}^{(\beta/2,\beta)}([0,T-\eta]\times\mathbb{R}^P,\mathbb{R}^Q)} \leq C_{2.7}^{(\eta)},$$

where $C_{2.7}^{(\eta)}$ only depends on $\bar{k}$, $\bar{\lambda}$, $\bar{\Lambda}$, $P$, $Q$, $T$ and $\eta$.

Note that estimates (2.11)–(2.13) will be crucial to establish the homogenization property.

**3. Strategy to prove the homogenization property.** In this section we wish to expose our strategy to prove the homogenization property. This strategy is then applied in Section 4. In the whole section, the additional hypothesis $(\mathcal{H}.8)$ is assumed to be in force.

**3.1.** *A few words on modified processes.* Referring to the notation introduced in (2.3) and according to (2.4) and (2.5), for every $\varepsilon > 0$, the forward and backward equations of FBSDE $\mathrm{E}(\varepsilon,0,x_0)$ can be put into the following form:

$$\forall\, t \in [0,T], \qquad \hat{X}_t(\varepsilon) = X_t(\varepsilon) + \varepsilon\hat{b}\big(\bar{X}_t(\varepsilon), Y_t(\varepsilon)\big)$$
$$= x_0 + \varepsilon\hat{b}\left(\frac{x_0}{\varepsilon}, Y_0(\varepsilon)\right)$$
$$(3.1)$$
$$+ \int_0^t u\big(\bar{X}_s(\varepsilon), Y_s(\varepsilon), \hat{Z}_s(\varepsilon)\big)\, ds$$
$$+ \int_0^t [(I + \nabla_x\hat{b})\sigma]\big(\bar{X}_s(\varepsilon), Y_s(\varepsilon)\big)\, dB_s + R_t(\varepsilon),$$

$$\forall\, t \in [0,T], \qquad \hat{Y}_t(\varepsilon) = Y_t(\varepsilon) - \varepsilon\hat{e}\big(\bar{X}_t(\varepsilon), Y_t(\varepsilon)\big)$$
$$= H(X_T(\varepsilon)) - \varepsilon\hat{e}\big(\bar{X}_T(\varepsilon), Y_T(\varepsilon)\big)$$
$$(3.2)$$
$$+ \int_t^T v\big(\bar{X}_s(\varepsilon), Y_s(\varepsilon), \hat{Z}_s(\varepsilon)\big)\, ds$$
$$- \int_t^T \hat{Z}_s(\varepsilon)\sigma\big(\bar{X}_s(\varepsilon), Y_s(\varepsilon)\big)\, dB_s + S_T(\varepsilon) - S_t(\varepsilon),$$

where, for every $\varepsilon > 0$ and every $t \in [0, T]$,

$$
\begin{aligned}
R_t(\varepsilon) = {} & -\varepsilon \int_0^t (\nabla_y \hat{b} f)(\bar{X}_s(\varepsilon), Y_s(\varepsilon), Z_s(\varepsilon))\, ds \\
& + \tfrac{1}{2}\varepsilon \int_0^t \nabla_{y,y}^2 \hat{b}(\bar{X}_s(\varepsilon), Y_s(\varepsilon))[Z_s(\varepsilon) a(\bar{X}_s(\varepsilon), Y_s(\varepsilon)) Z_s^*(\varepsilon)]\, ds \\
& + \varepsilon \int_0^t \nabla_y \hat{b}(\bar{X}_s(\varepsilon), Y_s(\varepsilon)) Z_s(\varepsilon) \sigma(\bar{X}_s(\varepsilon), Y_s(\varepsilon))\, dB_s \\
= {} & \int_0^t R_s(1, \varepsilon)\, ds + \int_0^t R_s(2, \varepsilon)\, dB_s,
\end{aligned}
$$

$$
\begin{aligned}
S_t(\varepsilon) = {} & -\varepsilon \int_0^t (\nabla_y \hat{e} f)(\bar{X}_s(\varepsilon), Y_s(\varepsilon), Z_s(\varepsilon))\, ds \\
& + \tfrac{1}{2}\varepsilon \int_0^t \nabla_{y,y}^2 \hat{e}(\bar{X}_s(\varepsilon), Y_s(\varepsilon))[Z_s(\varepsilon) a(\bar{X}_s(\varepsilon), Y_s(\varepsilon)) Z_s^*(\varepsilon)]\, ds \\
& + \varepsilon \int_0^t \nabla_y \hat{e}(\bar{X}_s(\varepsilon), Y_s(\varepsilon)) Z_s(\varepsilon) \sigma(\bar{X}_s(\varepsilon), Y_s(\varepsilon))\, dB_s \\
= {} & \int_0^t S_s(1, \varepsilon)\, ds + \int_0^t S_s(2, \varepsilon)\, dB_s.
\end{aligned}
$$

(3.3)

(3.4)

Thanks to Theorem 2.2 (which gives bounds of the solutions of the auxiliary problems) and to Theorem 2.3 [which provides estimates of $(X(\varepsilon), Y(\varepsilon), Z(\varepsilon))_{\varepsilon > 0}$], it is readily seen that extra terms $(R(\varepsilon))_{\varepsilon > 0}$ and $(S(\varepsilon))_{\varepsilon > 0}$ are negligible as $\varepsilon$ tends to 0:

$$
\begin{aligned}
\lim_{\varepsilon \to 0} \mathbf{E} \Big[ & \Big( \int_0^T |R_s(1, \varepsilon)|\, ds \Big)^2 + \Big( \int_0^T |S_s(1, \varepsilon)|\, ds \Big)^2 \\
& + \Big( \int_0^T |R_s(2, \varepsilon)|^2\, ds \Big)^2 + \Big( \int_0^T |S_s(2, \varepsilon)|^2\, ds \Big)^2 \Big] = 0.
\end{aligned}
$$

(3.5)

In particular, $\mathbf{E} \sup_{0 \le t \le T}(|R_t(\varepsilon)|^2 + |S_t(\varepsilon)|^2) \to 0$ as $\varepsilon \to 0$.

3.2. *A few words on the stability property approach.* As mentioned in the Introduction, we wish to apply the approach due to Buckdahn, Hu and Peng (1999). Roughly speaking, this method aims to compare $\theta(\cdot, \hat{X}.(\varepsilon))$ with $\hat{Y}.(\varepsilon)$. To reach such an objective, the authors apply the Itô formula to the quantity $\theta(\cdot, \hat{X}.(\varepsilon))$ in order to write, thanks to the system of PDEs $\mathcal{E}(\lim)$, this process as the solution of a BSDE. Using stability properties of BSDEs, they are then able to compare $\theta(\cdot, \hat{X}.(\varepsilon))$ with $\hat{Y}.(\varepsilon)$.

Of course, this strategy fails if the map $\theta$ is not regular enough. As written in Section 2.5, $\theta$ belongs to the space $\mathcal{V}$. Unfortunately, since we do not know if $(I + \nabla_x \hat{b})\sigma$ is elliptic, we cannot apply the so-called Itô–Krylov formula [see Krylov (1980), Chapter 2, Paragraph 10, Theorem 1] to the term $\theta(\cdot, \hat{X}.(\varepsilon))$. In fact, to get round this difficulty, we still follow the paper of Buckdahn, Hu and Peng (1999) and introduce a regularization sequence of $\theta$.

### 3.3. *Regularization of $\theta$ and related systems of PDEs.*

We first introduce the following regularization sequences of $H$ and $\bar{v}$. Using the notation given in the Introduction, we set for every $n \in \mathbb{N}^*$,

$$\forall x \in \mathbb{R}^P, \qquad H_n(x) = n^P \int H(x') \rho_P(n(x - x')) \, dx',$$

(3.6)  $\forall (y, z) \in \mathbb{R}^Q \times \mathbb{R}^{Q \times P},$

$$\bar{v}_n(y, z) = n^{Q+QP} \int \bar{v}(y', z') \rho_{Q+QP}(n(y - y', z - z')) \, dy' \, dz'.$$

It is readily seen that for every $n \in \mathbb{N}^*$, $H_n$ and $\bar{v}_n$ are infinitely differentiable. Moreover, for every $q \in \mathbb{N}^*$, the derivatives of order $q$ of $H_n$ are bounded.

Following the proof of Proposition 2.2 in Delarue (2002b) and modifying $\bar{\Lambda}$ and $\bar{K}$ if necessary, we can assume without loss of generality that, for every $n \in \mathbb{N}^*$, $\bar{u}$, $\bar{v}_n$, $H_n$ and $\bar{\alpha}^{1/2}$ satisfy (2.5.ii) with respect to $\bar{k}$, $\bar{\lambda}$, $\bar{\Lambda}$ and $\bar{K}$. Hence, thanks to Appendix D, we claim that, for every $n \in \mathbb{N}^*$, the system of quasilinear PDEs.

For $(t, x) \in [0, T[ \times \mathbb{R}^P$ and $\ell \in \{1, \dots, Q\}$,

$$\mathcal{E}_{\mathrm{reg}}(n)$$

$$\frac{\partial (\zeta_n)_\ell}{\partial t}(t, x) + \frac{1}{2} \sum_{i,j=1}^P \bar{\alpha}_{i,j}(\zeta_n(t, x)) \frac{\partial^2 (\zeta_n)_\ell}{\partial x_i \partial x_j}(t, x)$$

$$+ \sum_{i=1}^P \bar{u}_i(\zeta_n(t, x), \nabla_x \zeta_n(t, x)) \frac{\partial (\zeta_n)_\ell}{\partial x_i}(t, x)$$

$$+ (\bar{v}_n)_\ell(\zeta_n(t, x), \nabla_x \zeta_n(t, x)) = 0,$$

$$\zeta_n(T, x) = H_n(x),$$

admits a unique bounded classical solution $\zeta_n$. $\forall 0 < \eta < 1$, $\zeta_n \in \mathcal{C}^{(1+\eta/2, 2+\eta)}([0, T] \times \mathbb{R}^P, \mathbb{R}^Q)$.

Referring to Delarue (2002b) [see also Ladyzenskaja, Solonnikov and Ural'ceva (1968) and Delarue (2003)], there exists a constant $C_{3.1}$, only depending on $\bar{k}$, $\bar{\lambda}$, $\bar{\Lambda}$, $P$, $Q$ and $T$, such that

(3.7)  $\forall n \in \mathbb{N}^*, \ \forall (t, x) \in [0, T] \times \mathbb{R}^P, \qquad |\zeta_n(t, x)| + |\nabla_x \zeta_n(t, x)| \leq C_{3.1}.$

According to Appendix D, for every $n \in \mathbb{N}^*$, there exists a constant $C_{3.2}(n)$ such that

(3.8)  $\forall (t, x) \in [0, T] \times \mathbb{R}^P, \qquad |\nabla^2_{x,x} \zeta_n(t, x)| \leq C_{3.2}(n).$

Finally, thanks again to Appendix D, $(\zeta_n)_{n \in \mathbb{N}^*}$ uniformly converges on every compact subset of $[0, T] \times \mathbb{R}^P$ toward $\theta$ and $(\nabla_x \zeta_n)_{n \in \mathbb{N}^*}$ uniformly converges on every compact subset of $[0, T[ \times \mathbb{R}^P$ toward $\nabla_x \theta$.

3.4. *Regularized solutions and modified processes.* Following our strategy, we apply, for every $n \in \mathbb{N}^*$, Itô's formula to the quantity $\zeta_n(\cdot, \hat{X}.(\varepsilon))$. We then obtain for every $t \in [0, T]$,

$$
\begin{aligned}
\zeta_n(t, \hat{X}_t(\varepsilon)) = {}& H_n(\hat{X}_T(\varepsilon)) + \int_t^T \bar{v}_n(\zeta_n(s, \hat{X}_s(\varepsilon)), \nabla_x \zeta_n(s, \hat{X}_s(\varepsilon))) \, ds \\
& + \int_t^T \nabla_x \zeta_n(s, \hat{X}_s(\varepsilon)) \\
& \qquad \times [\bar{u}(\zeta_n(s, \hat{X}_s(\varepsilon)), \nabla_x \zeta_n(s, \hat{X}_s(\varepsilon))) \\
& \qquad\qquad - u(\bar{X}_s(\varepsilon), Y_s(\varepsilon), \hat{Z}_s(\varepsilon))] \, ds \\
& + \tfrac{1}{2} \int_t^T \nabla_{x,x}^2 \zeta_n(s, \hat{X}_s(\varepsilon)) [\bar{\alpha}(\zeta_n(s, \hat{X}_s(\varepsilon))) - \alpha(\bar{X}_s(\varepsilon), Y_s(\varepsilon))] \, ds \\
& - \int_t^T \nabla_x \zeta_n(s, \hat{X}_s(\varepsilon)) [(I + \nabla_x \hat{b}) \sigma](\bar{X}_s(\varepsilon), Y_s(\varepsilon)) \, dB_s \\
& - \int_t^T [\tfrac{1}{2} \nabla_{x,x}^2 \zeta_n(s, \hat{X}_s(\varepsilon)) \\
& \qquad \times ([(I + \nabla_x \hat{b}) \sigma] R_s^*(2, \varepsilon) \\
& \qquad\qquad + R_s(2, \varepsilon) [(I + \nabla_x \hat{b}) \sigma]^* \\
& \qquad\qquad\qquad + R_s(2, \varepsilon) R_s^*(2, \varepsilon)) (\bar{X}_s(\varepsilon), Y_s(\varepsilon)) \, ds \\
& \qquad\qquad\qquad\qquad + \nabla_x \zeta_n(s, \hat{X}_s(\varepsilon)) \, dR_s(\varepsilon)].
\end{aligned}
$$

(3.9)

Therefore, from (3.2), we deduce that for every $t \in [0, T]$,

$$
\begin{aligned}
\zeta_n(t, & \hat{X}_t(\varepsilon)) - \hat{Y}_t(\varepsilon) \\
= {}& H_n(\hat{X}_T(\varepsilon)) - \hat{Y}_T(\varepsilon) \\
& + \int_t^T [\bar{v}_n(\zeta_n(s, \hat{X}_s(\varepsilon)), \nabla_x \zeta_n(s, \hat{X}_s(\varepsilon))) \\
& \qquad - v(\bar{X}_s(\varepsilon), Y_s(\varepsilon), \hat{Z}_s(\varepsilon))] \, ds \\
& + \int_t^T \nabla_x \zeta_n(s, \hat{X}_s(\varepsilon)) [\bar{u}(\zeta_n(s, \hat{X}_s(\varepsilon)), \nabla_x \zeta_n(s, \hat{X}_s(\varepsilon))) \\
& \qquad\qquad - u(\bar{X}_s(\varepsilon), Y_s(\varepsilon), \hat{Z}_s(\varepsilon))] \, ds \\
& + \tfrac{1}{2} \int_t^T \nabla_{x,x}^2 \zeta_n(s, \hat{X}_s(\varepsilon)) \\
& \qquad \times [\bar{\alpha}(\zeta_n(s, \hat{X}_s(\varepsilon))) - \alpha(\bar{X}_s(\varepsilon), Y_s(\varepsilon))] \, ds \\
& - \int_t^T [\nabla_x \zeta_n(s, \hat{X}_s(\varepsilon)) (I + \nabla_x \hat{b})(\bar{X}_s(\varepsilon), Y_s(\varepsilon)) - \hat{Z}_s(\varepsilon)]
\end{aligned}
$$

(3.10)

$$\times \sigma(\bar{X}_s(\varepsilon), Y_s(\varepsilon)) \, dB_s$$

$$- \int_t^T [\tfrac{1}{2} \nabla_{x,x}^2 \zeta_n(s, \hat{X}_s(\varepsilon))$$

$$\times ([(I + \nabla_x \hat{b}) \sigma] R_s^*(2, \varepsilon)$$

$$+ R_s(2, \varepsilon) [(I + \nabla_x \hat{b}) \sigma]^*$$

$$+ R_s(2, \varepsilon) R_s^*(2, \varepsilon))(\bar{X}_s(\varepsilon), Y_s(\varepsilon)) \, ds$$

$$+ \nabla_x \zeta_n(s, \hat{X}_s(\varepsilon)) \, dR_s(\varepsilon) + dS_s(\varepsilon)].$$

Hence, we can write (with obvious notation)

$$(3.11) \quad \begin{aligned} &\zeta_n(t, \hat{X}_t(\varepsilon)) - \hat{Y}_t(\varepsilon) \\ &\quad = [H_n(\hat{X}_T(\varepsilon)) - H(X_T(\varepsilon)) + \varepsilon \hat{e}(\bar{X}_T(\varepsilon), Y_T(\varepsilon))] \\ &\qquad + \Delta_{3.1}(t, \varepsilon, n) + \Delta_{3.2}(t, \varepsilon, n) + \Delta_{3.3}(t, \varepsilon, n) \\ &\qquad + \Delta_{3.4}(t, \varepsilon, n) + \Delta_{3.5}(t, \varepsilon, n). \end{aligned}$$

Recall at this step that we aim to estimate the distance between $\zeta_n(\cdot, \hat{X}.(\varepsilon))$ and $\hat{Y}.(\varepsilon)$. Recall also that the usual strategy to achieve such an objective consists in writing the process $|\zeta_n(\cdot, \hat{X}.(\varepsilon)) - \hat{Y}.(\varepsilon)|^2$ as a semimartingale. Note, for the sake of simplicity, that we won't detail this operation right now, but keep in mind that we should actually perform it. In particular, pay attention to the fact that this operation usually permits one to apply so-called "monotonicity assumption" [as written in $(\mathcal{H}.2)$].

Note now from ergodic Theorem 2.4 that the right-hand side in (3.11) can be approximated in the following way. For all $n \in \mathbb{N}^*$ and $t \in [0, T]$,

$$\Delta_{3.1}(t, \varepsilon, n)$$

$$(3.12) \quad \begin{aligned} &\underset{\varepsilon \to 0}{\approx} \int_t^T [\bar{v}_n(\zeta_n(s, \hat{X}_s(\varepsilon)), \nabla_x \zeta_n(s, \hat{X}_s(\varepsilon))) \\ &\qquad\qquad - \bar{v}(Y_s(\varepsilon), \nabla_x \zeta_n(s, \hat{X}_s(\varepsilon)))] \, ds \\ &\quad - \int_t^T [v(\bar{X}_s(\varepsilon), Y_s(\varepsilon), \\ &\qquad\qquad \nabla_x \zeta_n(s, \hat{X}_s(\varepsilon))(I + \nabla_x \hat{b})(\bar{X}_s(\varepsilon), Y_s(\varepsilon))) \\ &\qquad\qquad\qquad - v(\bar{X}_s(\varepsilon), Y_s(\varepsilon), \hat{Z}_s(\varepsilon))] \, ds, \end{aligned}$$

$$\Delta_{3.2}(t, \varepsilon, n)$$

$$\underset{\varepsilon\to 0}{\approx} \int_t^T ds \nabla_x \zeta_n(s, \hat{X}_s(\varepsilon))$$

(3.13)
$$\times [\bar{u}(\zeta_n(s, \hat{X}_s(\varepsilon)), \nabla_x \zeta_n(s, \hat{X}_s(\varepsilon))) - \bar{u}(Y_s(\varepsilon), \nabla_x \zeta_n(s, \hat{X}_s(\varepsilon)))$$

$$+ u(\bar{X}_s(\varepsilon), Y_s(\varepsilon), \nabla_x \zeta_n(s, \hat{X}_s(\varepsilon))(I + \nabla_x \hat{b})(\bar{X}_s(\varepsilon), Y_s(\varepsilon)))$$

$$- u(\bar{X}_s(\varepsilon), Y_s(\varepsilon), \hat{Z}_s(\varepsilon))].$$

The sense given to the symbol " $\underset{\varepsilon\to 0}{\approx}$ " will be clearly detailed in Section 4.

Recall now from Section 2.5 that $u$ and $\bar{u}$ are locally Lipschitzian in $y$ and $z$, and that $v$ and $\bar{v}$ are locally monotonous in $y$ and locally Lipschitzian in $z$. Recall also from Section 3.3 that the functions $(\nabla_x \zeta_n)_{n\in\mathbb{N}^*}$ are bounded on $[0,T]\times\mathbb{R}^P$, uniformly in $n$. We are then able to estimate $\Delta_{3.1}(\cdot,\varepsilon,n)$ and $\Delta_{3.2}(\cdot,\varepsilon,n)$ by a negligible term, the distance between $\bar{v}_n$ and $\bar{v}$ and the distance between $(\zeta_n(\cdot,\hat{X}.(\varepsilon)),\nabla_x \zeta_n(\cdot,\hat{X}.(\varepsilon))(I + \nabla_x \hat{b})(\bar{X}.(\varepsilon),Y.(\varepsilon)))$ and $(Y.(\varepsilon),\hat{Z}.(\varepsilon))$. $\Delta_{3.4}(\cdot,\varepsilon,n)$ can be also estimated by the distance between $\nabla_x \zeta_n(\cdot,\hat{X}.(\varepsilon))(I + \nabla_x \hat{b})(\bar{X}.(\varepsilon),Y.(\varepsilon))$ and $\hat{Z}.(\varepsilon)$. Thanks to Section 3.1, $\Delta_{3.5}(\cdot,\varepsilon,n)$ is a negligible term.

Note, however, that the latter approximations are not uniform in $n$. To face this point, recall the conclusion of the former section: $\zeta_n(t,\hat{X}_t(\varepsilon)) - \hat{Y}_t(\varepsilon)$ can be estimated by $\Delta_{3.3}(\cdot,\varepsilon,n)$, a term $\delta(\varepsilon,n)$ that tends to 0 as $\varepsilon\to 0$ for every $n\in\mathbb{N}^*$, and the distance between the triples $(H_n(\hat{X}_T(\varepsilon)),\zeta_n(\cdot,\hat{X}.(\varepsilon)),\nabla_x \zeta_n(\cdot,\hat{X}.(\varepsilon))(I + \nabla_x \hat{b})(\bar{X}.(\varepsilon),Y.(\varepsilon)))$ and $(H(X_T(\varepsilon)),Y.(\varepsilon),\hat{Z}.(\varepsilon))$. Hence, if we manage to treat $\Delta_{3.3}$ as $\Delta_{3.1}$, $\Delta_{3.2}$ and $\Delta_{3.4}$, we can apply a Gronwall argument as often done in BSDEs theory: $\zeta_n(t,\hat{X}_t(\varepsilon)) - \hat{Y}_t(\varepsilon)$ can be then estimated by $\delta(\varepsilon,n)$ and the distance between $H_n(\hat{X}_T(\varepsilon))$ and $H(X_T(\varepsilon))$. Recall finally that $(H_n)_{n\in\mathbb{N}^*}$ uniformly converges on $\mathbb{R}^P$ toward $H$ and that $(\zeta_n)_{n\in\mathbb{N}^*}$ uniformly converges on every compact subset of $[0,T]\times\mathbb{R}^P$ toward $\theta$: fixing first a large enough $n\in\mathbb{N}^*$, we then deduce that $\theta(t,\hat{X}_t(\varepsilon)) - \hat{Y}_t(\varepsilon)$ tends to 0 as $\varepsilon\to 0$. Choosing $t = 0$, the homogenization property is then established.

In fact, we are not able to treat $\Delta_{3.3}(\cdot,\varepsilon,n)$ as the other terms. Indeed, we deduce from Theorem 2.4 that

$$\Delta_{3.3}(t,\varepsilon,n)$$

(3.14)
$$\underset{\varepsilon\to 0}{\approx} \frac{1}{2} \int_t^T \nabla^2_{x,x} \zeta_n(s, \hat{X}_s(\varepsilon))[\bar{\alpha}(\zeta_n(s, \hat{X}_s(\varepsilon))) - \bar{\alpha}(Y_s(\varepsilon))]\, ds.$$

Here, the situation is rather different from (3.12) and (3.13): $\bar{\alpha}$ is locally Lipschitzian with respect to $y$, but we do not control the supremum norm of $\nabla^2_{x,x} \zeta_n$ over $[0,T]\times\mathbb{R}^P$ uniformly in $n$. Actually, referring once again to Gronwall's lemma, we can see that we should then be able to control such a quantity as $\int_0^T \sup_{x\in\mathbb{R}^P} |\nabla^2_{x,x} \zeta_n(s,x)|\, ds$ in order to estimate in a good way $\Delta_{3.3}(\cdot, \varepsilon, n)$. Unfortunately, our assumptions seem too weak to let us establish such an estimate [at the opposite, such an estimate holds in the case treated by Buckdahn and Hu (1998)]. Note, moreover, that such a difficulty does not appear in the semilinear case.

3.5. *Rough definition of auxiliary SDEs.* We now explain how to face the latter difficulty—this part seems to be completely new in the probabilistic literature devoted to homogenization. In fact, instead of estimating the distance between $\zeta_n(\cdot, \hat{X}.(\varepsilon))$ and $\hat{Y}.(\varepsilon)$, we study the distance between $\zeta_n(\cdot, U.(\varepsilon, n))$ and $\hat{Y}.(\varepsilon)$, where $U(\varepsilon, n)$ denotes the solution of a "well-chosen" auxiliary SDE. Roughly speaking, this SDE is given by

$$\forall t \in [0, T], \qquad U_t(\varepsilon, n) = x_0 + \int_0^t u(\bar{X}_s(\varepsilon), Y_s(\varepsilon), \hat{Z}_s(\varepsilon)) \, ds$$

$$(3.15) \qquad\qquad + \int_0^t \frac{p^{1/2}(\bar{X}_s(\varepsilon), \zeta_n(s, U_s(\varepsilon, n)))}{p^{1/2}(\bar{X}_s(\varepsilon), Y_s(\varepsilon))}$$

$$\times [(I + \nabla_x \hat{b})\sigma](\bar{X}_s(\varepsilon), \zeta_n(s, U_s(\varepsilon, n))) \, dB_s.$$

In order to explain the interest of this SDE, we assume for the moment that (3.15) is solvable. Then, applying Itô's formula to the quantity $\zeta_n(\cdot, U.(\varepsilon, n))$, we obtain for every $t \in [0, T]$,

$$\zeta_n(t, U_t(\varepsilon, n)) - \hat{Y}_t(\varepsilon)$$

$$= H_n(U_T(\varepsilon, n)) - H(X_T(\varepsilon)) + \varepsilon \hat{e}(\bar{X}_T(\varepsilon), Y_T(\varepsilon))$$

$$+ \int_t^T [\bar{v}_n(\zeta_n(s, U_s(\varepsilon, n)), \nabla_x \zeta_n(s, U_s(\varepsilon, n)))$$

$$- v(\bar{X}_s(\varepsilon), Y_s(\varepsilon), \hat{Z}_s(\varepsilon))] \, ds$$

$$+ \int_t^T \nabla_x \zeta_n(s, U_s(\varepsilon, n)) [\bar{u}(\zeta_n(s, U_s(\varepsilon, n)), \nabla_x \zeta_n(s, U_s(\varepsilon, n)))$$

$$- u(\bar{X}_s(\varepsilon), Y_s(\varepsilon), \hat{Z}_s(\varepsilon))] \, ds$$

$$(3.16) \qquad + \frac{1}{2} \int_t^T \nabla^2_{x,x} \zeta_n(s, U_s(\varepsilon, n))$$

$$\times \left[ \bar{\alpha}(\zeta_n(s, U_s(\varepsilon, n))) - \frac{p(\bar{X}_s(\varepsilon), \zeta_n(s, U_s(\varepsilon, n)))}{p(\bar{X}_s(\varepsilon), Y_s(\varepsilon))} \right.$$

$$\times \alpha(\bar{X}_s(\varepsilon), \zeta_n(s, U_s(\varepsilon, n))) \bigg] \, ds$$

$$- \int_t^T \left[ \frac{p^{1/2}(\bar{X}_s(\varepsilon), \zeta_n(s, U_s(\varepsilon, n)))}{p^{1/2}(\bar{X}_s(\varepsilon), Y_s(\varepsilon))} \nabla_x \zeta_n(s, U_s(\varepsilon, n)) \right.$$

$$\times \, [(I + \nabla_x \hat{b}) \sigma] (\, \bar{X}_s(\varepsilon), \zeta_n(s, U_s(\varepsilon, n)))$$

$$- \, \hat{Z}_s(\varepsilon) \sigma(\, \bar{X}_s(\varepsilon), Y_s(\varepsilon)) \Big] dB_s$$

$$- \int_t^T dS_s(\varepsilon).$$

Assume that we can apply Theorem 2.4 to the second-order term appearing in (3.16). Then, for every $n \in \mathbb{N}^*$ and for every $t \in [0, T]$,

$$
\begin{aligned}
(3.17) \quad & \int_t^T \nabla_{x,x}^2 \zeta_n(s, U_s(\varepsilon, n)) \\
& \times \Big[ \bar{\alpha}(\zeta_n(s, U_s(\varepsilon, n))) - \frac{p(\bar{X}_s(\varepsilon), \zeta_n(s, U_s(\varepsilon, n)))}{p(\bar{X}_s(\varepsilon), Y_s(\varepsilon))} \\
& \times \alpha(\, \bar{X}_s(\varepsilon), \zeta_n(s, U_s(\varepsilon, n))) \Big] \, ds \underset{\varepsilon \to 0}{\approx} 0.
\end{aligned}
$$

Hence, (3.17) is a negligible term. Applying the strategy given in Section 3.4, we are finally able to face the explosion of the second-order terms $(\nabla_{x,x}^2 \zeta_n)_{n \in \mathbb{N}^*}$.

### 3.6. *Remaining difficulties.* In fact, there remain three main difficulties:

1. First, due to the term $H_n(U_T(\varepsilon, n)) - H(X_T(\varepsilon))$ in (3.16), we also have to estimate the distance between the processes $U(\varepsilon, n)$ and $X(\varepsilon)$.
2. Second, as usual when dealing with FBSDEs, we have to prove first a local version of the global property given in Point 2 of Theorem (PHP) and then complete the proof by induction.
3. The third difficulty comes from the regularity of $p$: since we just know that $y \mapsto p(\cdot, y) \in L^2(\mathbb{T}^P)$ is Lipschitzian, we cannot solve the SDE (3.15) and cannot apply ergodic Theorem 2.4 in (3.17).

Let us explain how to face the second and third difficulties.

### 3.6.1. *Localization.* As explained above, we have to prove first a local version of the homogenization property. To this end, we will first establish the following:

THEOREM 3.1. *There exists a constant $0 < c_{3.1} \leq T$, only depending on $k$, $K$, $\lambda$, $\Lambda$, $P$, $Q$ and $T$, such that*

$$
\begin{aligned}
(3.18) \quad & \lim_{\varepsilon \to 0} \Big[ \mathbf{E} \sup_{T - c_{3.1} \leq t \leq T} |Y_t(\varepsilon) - \theta(t, X_t(\varepsilon))|^2 \\
& + \mathbf{E} \int_{T - c_{3.1}}^T |\hat{Z}_t(\varepsilon) - \nabla_x \theta(t, X_t(\varepsilon))(I + \nabla_x \hat{b})(\, \bar{X}_t(\varepsilon), Y_t(\varepsilon))|^2 \, dt \Big] = 0.
\end{aligned}
$$

To this end, we fix $t_0 \in [0, T]$. We will then prove that (3.18) holds on $[t_0, T]$ as soon as $t_0$ is close enough to $T$. As a consequence, we have to define $U(\varepsilon, n)$ on the interval $[t_0, T]$, and therefore, to choose its new initial condition. It is quite natural to set $U_{t_0}(\varepsilon, n) = X_{t_0}(\varepsilon)$.

3.6.2. *Regularization of $p$ and related auxiliary SDEs.* To face the third difficulty, we have to regularize the function $p$. Hence, we set for every $m \in \mathbb{N}^*$,

$$(3.19) \qquad \forall (x, y) \in \mathbb{R}^P \times \mathbb{R}^Q, \qquad p_m(x, y) = m^Q \int p(x, y') \rho_Q(m(y - y')) \, dy'.$$

Thanks to Theorem 2.1, the sequence of functions $(y \in \mathbb{R}^Q \mapsto \|p_m(\cdot, y) - p(\cdot, y)\|_2)_{m \in \mathbb{N}^*}$ tends to 0 as $m \to +\infty$, uniformly on every compact set. Note also that, for all $m \in \mathbb{N}^*$, $p_m$ is infinitely differentiable with respect to $y$, with bounded derivatives of every order, and also satisfies Theorem 2.1.

Recall that for every $n \in \mathbb{N}^*$, $\zeta_n$ is Lipschitzian with respect to $x$ [see (3.7)]. Hence, following the definition of the auxiliary SDE given in (3.15), we now claim:

PROPOSITION 3.2. *Let $(m(n))_{n \in \mathbb{N}^*} \in \mathbb{N}^{\mathbb{N}^*}$. Then, for all $n \in \mathbb{N}^*$ and $\varepsilon > 0$, the SDE*

$$\forall t \in [t_0, T], \qquad U_t(\varepsilon, n) = X_{t_0}(\varepsilon) + \int_{t_0}^{t} u(\bar{X}_s(\varepsilon), Y_s(\varepsilon), \hat{Z}_s(\varepsilon)) \, ds,$$

$$(3.20) \qquad\qquad + \int_{t_0}^{t} \frac{p_{m(n)}^{1/2}(\bar{X}_s(\varepsilon), \zeta_n(s, U_s(\varepsilon, n)))}{p_{m(n)}^{1/2}(\bar{X}_s(\varepsilon), Y_s(\varepsilon))}$$

$$\qquad\qquad\qquad \times [(I + \nabla_x \hat{b}) \sigma](\bar{X}_s(\varepsilon), \zeta_n(s, U_s(\varepsilon, n))) \, dB_s,$$

*admits a unique solution. We still denote by $U(\varepsilon, n)$ the solution of this SDE [which is different from (3.15)]. Note that the choice of the sequence $(m(n))_{n \in \mathbb{N}^*}$ will be described in a few lines.*

Following the definition of the homogenized coefficients given in Section 1.3, we set for all $n \in \mathbb{N}^*$ and $(x, y, y') \in \mathbb{T}^P \times \mathbb{R}^Q \times \mathbb{R}^Q$,

$$\alpha_n(x, y, y') = \frac{p_{m(n)}(x, y)}{p_{m(n)}(x, y')} \alpha(x, y),$$

$$(3.21) \qquad \bar{\alpha}_n(y, y') = \int_{\mathbb{T}^P} \alpha_n(x, y, y') p(x, y') \, dx$$

$$\qquad\qquad = \int_{\mathbb{T}^P} \alpha(x, y) p_{m(n)}(x, y) \frac{p(x, y')}{p_{m(n)}(x, y')} \, dx.$$

Let us now choose $(m(n))_{n\in\mathbb{N}^*}$: when applying Itô's formula to the quantity $\zeta_n(t, U_t(\varepsilon, n))$, we obtain formula (3.16) with $p$ replaced by $p_{m(n)}$. The second-order term is then given by

$$
\begin{aligned}
(3.22) \quad & \int_t^T \nabla^2_{x,x}\zeta_n(s, U_s(\varepsilon, n)) \\
& \times \Big[ \bar{\alpha}(\zeta_n(s, U_s(\varepsilon, n))) - \frac{p_{m(n)}(\bar{X}_s(\varepsilon), \zeta_n(s, U_s(\varepsilon, n)))}{p_{m(n)}(\bar{X}_s(\varepsilon), Y_s(\varepsilon))} \\
& \qquad\qquad\qquad \times \alpha(\bar{X}_s(\varepsilon), \zeta_n(s, U_s(\varepsilon, n))) \Big] ds.
\end{aligned}
$$

From ergodic Theorem 2.4, this term looks like

$$
\begin{aligned}
(3.23) \quad & \int_t^T \nabla^2_{x,x}\zeta_n(s, U_s(\varepsilon, n)) \\
& \times [\bar{\alpha}(\zeta_n(s, U_s(\varepsilon, n))) - \bar{\alpha}_n(\zeta_n(s, U_s(\varepsilon, n)), Y_s(\varepsilon))] ds.
\end{aligned}
$$

It is then clear that $(m(n))_{n\in\mathbb{N}^*}$ has to be chosen to make (3.23) tend toward 0 as $n \to +\infty$. We let the reader see that the following choice implies that $((t,x,y,y') \mapsto \nabla^2_{x,x}\zeta_n(t,x)(\bar{\alpha}(y) - \bar{\alpha}_n(y,y')))_{n\in\mathbb{N}^*}$ uniformly converges on $[0,T] \times \mathbb{R}^P \times \bar{B}_Q(0, C_{3.1}) \times \bar{B}_Q(0, C_{2.4})$ toward 0 as $n \to +\infty$:

$$
\begin{aligned}
(3.24) \quad \forall n \in \mathbb{N}^*, \qquad m(n) = \inf\Big\{ & m \geq n, \\
& \sup_{[0,T]\times\mathbb{R}^P} |\nabla^2_{x,x}\zeta_n(t,x)| \\
& \times \sup_{|y| \leq C_{2.4}+C_{3.1}} \|p_m(\cdot, y) - p(\cdot, y)\|_2 \leq \frac{1}{n} \Big\}.
\end{aligned}
$$

Note that (3.8) is crucial to define $(m(n))_{n\in\mathbb{N}^*}$ as done above.

## 4. Proof of the homogenization property.

We now apply our strategy to prove the homogenization result [recall that $(\mathcal{H}.8)$ is in force]. To this end, we set for all $\varepsilon > 0$, $n \in \mathbb{N}^*$ and $t \in [t_0, T]$,

$$
\begin{aligned}
(4.1) \quad & V_t(\varepsilon, n) = \zeta_n(t, U_t(\varepsilon, n)), \\
& W_t(\varepsilon, n) = \nabla_x \zeta_n(t, U_t(\varepsilon, n)), \\
& \hat{W}_t(\varepsilon, n) = \nabla_x \zeta_n(t, U_t(\varepsilon, n))(I + \nabla_x \hat{b})(\bar{X}_t(\varepsilon), V_t(\varepsilon, n)).
\end{aligned}
$$

We then proceed in five steps: 1. We estimate the distance between $X(\varepsilon)$ and $U(\varepsilon, n)$. 2. We estimate the distance between $Y(\varepsilon)$ and $V(\varepsilon, n)$. 3. We establish Theorem 3.1 (recall that it is a local version of the homogenization property). 4. We establish point 2 of Theorem (PHP) and then deduce point 3 of Theorem (HP). 5. We establish point 3 of Theorem (PHP).

Before, we give a preliminary estimate of the processes involved in our strategy.

LEMMA 4.1.   $\sup_{\varepsilon>0, n\in\mathbb{N}^*} \mathbf{E}\sup_{t_0\leq t\leq T}[|X_t(\varepsilon)|^2 + |\hat{X}_t(\varepsilon)|^2 + |U_t(\varepsilon,n)|^2] < \infty$.

PROOF.   From Theorems 2.2 and 2.3, $\sup_{\varepsilon>0} \mathbf{E}\sup_{t_0\leq t\leq T}[|X_t(\varepsilon)|^2 + |\hat{X}_t(\varepsilon)|^2] < \infty$. Thanks to (3.20), to the properties of $(p_m)_{m\in\mathbb{N}^*}$ (which satisfy Theorem 2.1), to Section 2.5 (in which we have estimated $u$) and again to Theorem 2.2, $\sup_{\varepsilon>0, n\in\mathbb{N}^*} \mathbf{E}\sup_{t_0\leq t\leq T}[|U_t(\varepsilon,n)|^2] < \infty$.   □

4.1. *Estimate of the distance between $X(\varepsilon)$ and $U(\varepsilon,n)$.*   In this section we first compare $X(\varepsilon)$ with $U(\varepsilon,n)$.

PROPOSITION 4.2.   *There exists a constant $C_{4.1}$, only depending on $k$, $K$, $\lambda$, $\Lambda$, $P$, $Q$ and $T$, such that for all $n\in\mathbb{N}^*$, $\varepsilon>0$ and $t\in[t_0,T]$,*

$$(4.2)\quad \mathbf{E}\sup_{t_0\leq s\leq t} |U_s(\varepsilon,n) - X_s(\varepsilon)|^2 \leq \delta_{4.1}(\varepsilon,n) + C_{4.1}\mathbf{E}\int_{t_0}^t |V_s(\varepsilon,n) - Y_s(\varepsilon)|^2\,ds,$$

*where, for every $n\in\mathbb{N}^*$, $\lim_{\varepsilon\to 0}\delta_{4.1}(\varepsilon,n) = 0$.*

PROOF.   From (3.1) and (3.20), we deduce that for all $\varepsilon>0$, $n\in\mathbb{N}^*$ and $s\in[t_0,T]$,

$$d|U_s(\varepsilon,n) - \hat{X}_s(\varepsilon)|^2$$

$$= d\bigg[\int_{t_0}^\cdot \bigg(\frac{p_{m(n)}^{1/2}(\bar{X}_r(\varepsilon), V_r(\varepsilon,n))}{p_{m(n)}^{1/2}(\bar{X}_r(\varepsilon), Y_r(\varepsilon))}$$

$$\times [(I + \nabla_x\hat{b})\sigma](\bar{X}_r(\varepsilon), V_r(\varepsilon,n))$$

$$- [(I + \nabla_x\hat{b})\sigma](\bar{X}_r(\varepsilon), Y_r(\varepsilon))\bigg)\,dB_r - \int_{t_0}^\cdot R_r(2,\varepsilon)\,dB_r\bigg]_s$$

$$(4.3)\qquad + 2\bigg\langle U_s(\varepsilon,n) - \hat{X}_s(\varepsilon),$$

$$\bigg[\frac{p_{m(n)}^{1/2}(\bar{X}_s(\varepsilon), V_s(\varepsilon,n))}{p_{m(n)}^{1/2}(\bar{X}_s(\varepsilon), Y_s(\varepsilon))}[(I + \nabla_x\hat{b})\sigma](\bar{X}_s(\varepsilon), V_s(\varepsilon,n))$$

$$- [(I + \nabla_x\hat{b})\sigma](\bar{X}_s(\varepsilon), Y_s(\varepsilon))\bigg]\,dB_s\bigg\rangle$$

$$- 2\langle U_s(\varepsilon,n) - \hat{X}_s(\varepsilon), dR_s(\varepsilon)\rangle.$$

According to Section 3.1 and to Lemma 4.1, the influence of the terms $R(\varepsilon)$ and $R(2,\varepsilon)$ is negligible. In particular, from the Burkholder–Davis–Gundy inequalities, we obtain for all $t \in [t_0, T]$,

$$
\begin{aligned}
\mathbf{E} &\sup_{t_0 \leq s \leq t} |U_s(\varepsilon, n) - \hat{X}_s(\varepsilon)|^2 \\
&\leq \delta_{4.1}(\varepsilon, n) + c_{4.1} \mathbf{E} \int_{t_0}^{t} \Big| \frac{p_{m(n)}^{1/2}(\bar{X}_s(\varepsilon), V_s(\varepsilon, n))}{p_{m(n)}^{1/2}(\bar{X}_s(\varepsilon), Y_s(\varepsilon))} \\
&\qquad\qquad\qquad \times [(I + \nabla_x \hat{b})\sigma](\bar{X}_s(\varepsilon), V_s(\varepsilon, n)) \\
&\qquad\qquad\qquad\qquad - [(I + \nabla_x \hat{b})\sigma](\bar{X}_s(\varepsilon), Y_s(\varepsilon)) \Big|^2 ds \\
&= \delta_{4.1}(\varepsilon, n) + c_{4.1} \Delta_{4.1}(t, \varepsilon, n),
\end{aligned} \tag{4.4}
$$

where $c_{4.1}$ is a universal constant, and $(\delta_{4.1}(\varepsilon, n))_{\varepsilon > 0, n \in \mathbb{N}^*}$ is a family of reals whose values may change from one inequality to another and that satisfies, for all $n \in \mathbb{N}^*$, $\lim_{\varepsilon \to 0} \delta_{4.1}(\varepsilon, n) = 0$ [in fact, $\delta_{4.1}(\varepsilon, n)$ could be chosen in (4.4) independently of $n$, but this notation will be used next]. Applying Theorem 2.4 with $G = U(\varepsilon, n)$ [recall that $V.(\varepsilon, n) = \zeta_n(\cdot, U.(\varepsilon, n))$], we have for all $\varepsilon > 0$ and $n \in \mathbb{N}^*$,

$$
\begin{aligned}
\mathbf{E} \sup_{t_0 \leq t \leq T} &\Big| \int_{t_0}^{t} \Big[ \Big| \frac{p_{m(n)}^{1/2}(\bar{X}_s(\varepsilon), V_s(\varepsilon, n))}{p_{m(n)}^{1/2}(\bar{X}_s(\varepsilon), Y_s(\varepsilon))} [(I + \nabla_x \hat{b})\sigma](\bar{X}_s(\varepsilon), V_s(\varepsilon, n)) \\
&\qquad\qquad\qquad\qquad - [(I + \nabla_x \hat{b})\sigma](\bar{X}_s(\varepsilon), Y_s(\varepsilon)) \Big|^2 \\
&\qquad - \int_{\mathbb{T}^P} \Big| \frac{p_{m(n)}^{1/2}(x, V_s(\varepsilon, n))}{p_{m(n)}^{1/2}(x, Y_s(\varepsilon))} [(I + \nabla_x \hat{b})\sigma](x, V_s(\varepsilon, n)) \\
&\qquad\qquad\qquad\qquad - [(I + \nabla_x \hat{b})\sigma](x, Y_s(\varepsilon)) \Big|^2 \\
&\qquad\qquad\qquad\qquad\qquad \times p(x, Y_s(\varepsilon)) \, dx \Big] ds \Big| \\
&\leq \delta_{4.1}(\varepsilon, n).
\end{aligned} \tag{4.5}
$$

From Theorems 2.1 and 2.2, $p_{m(n)}$ and $(I + \nabla_x \hat{b})\sigma$ are bounded and $p_{m(n)}$ is bounded from below by a positive constant. Moreover, again from Theorems 2.1 and 2.2, the functions $y \in \mathbb{R}^Q \mapsto ((I + \nabla_x \hat{b})\sigma)(\cdot, y) \in L^2(\mathbb{T}^P)$ and $y \in \mathbb{R}^Q \mapsto p_{m(n)}(\cdot, y) \in L^2(\mathbb{T}^P)$ are Lipschitzian. Hence, for all $\varepsilon > 0$, $n \in \mathbb{N}^*$ and $t \in [t_0, T]$,

$$
\Delta_{4.1}(t, \varepsilon, n) \leq \delta_{4.1}(\varepsilon, n) + C_{4.1} \mathbf{E} \int_{t_0}^{t} |V_s(\varepsilon, n) - Y_s(\varepsilon)|^2 \, ds. \tag{4.6}
$$

We deduce that (4.2) holds with $X(\varepsilon)$ replaced by $\hat{X}(\varepsilon)$. Since $\hat{X}(\varepsilon) = X(\varepsilon) + \varepsilon \hat{b}$, we complete the proof. $\square$

4.2. *Estimate of the distance between $Y(\varepsilon)$ and $V(\varepsilon, n)$.* We now compare $Y(\varepsilon)$ with $V(\varepsilon, n)$.

PROPOSITION 4.3. *There exists a constant $c_{4.2}$, only depending on $k$, $K$, $\lambda$, $\Lambda$, $P$, $Q$ and $T$, such that for all $t_0 \in [T - c_{4.2}, T]$, $\varepsilon > 0$ and $n \in \mathbb{N}^*$,*

$$
\mathbf{E} \sup_{t_0 \le t \le T} |U_t(n) - X_t(\varepsilon)|^2
$$

$$
(4.7) \qquad + \mathbf{E} \sup_{t_0 \le t \le T} |V_t(\varepsilon, n) - Y_t(\varepsilon)|^2 + \mathbf{E} \int_{t_0}^{T} |\hat{W}_t(\varepsilon, n) - \hat{Z}_t(\varepsilon)|^2 \, dt
$$

$$
\le \delta_{4.2}(\varepsilon, n) + \delta_{4.3}(n),
$$

*where $\lim_{n \to +\infty} \delta_{4.3}(n) = 0$, and for every $n \in \mathbb{N}^*$, $\lim_{\varepsilon \to 0} \delta_{4.2}(\varepsilon, n) = 0$.*

PROOF. Following (3.16), we claim that for all $n \in \mathbb{N}^*$, $\varepsilon > 0$ and $t \in [t_0, T]$,

$$
|V_t(\varepsilon, n) - \hat{Y}_t(\varepsilon)|^2
$$

$$
= |V_T(\varepsilon, n) - \hat{Y}_T(\varepsilon)|^2
$$

$$
+ 2 \int_t^T \langle V_s(\varepsilon, n) - \hat{Y}_s(\varepsilon),
$$

$$
\bar{v}_n(V_s(\varepsilon, n), W_s(\varepsilon, n)) - v(\bar{X}_s(\varepsilon), Y_s(\varepsilon), \hat{Z}_s(\varepsilon)) \rangle \, ds
$$

$$
+ 2 \int_t^T \langle V_s(\varepsilon, n) - \hat{Y}_s(\varepsilon),
$$

$$
W_s(\varepsilon, n)[\bar{u}(V_s(\varepsilon, n), W_s(\varepsilon, n)) - u(\bar{X}_s(\varepsilon), Y_s(\varepsilon), \hat{Z}_s(\varepsilon))] \rangle \, ds
$$

$$
+ \int_t^T \langle V_s(\varepsilon, n) - \hat{Y}_s(\varepsilon), \nabla_{x,x}^2 \zeta_n(s, U_s(\varepsilon, n))
$$

$$
(4.8) \qquad\qquad\qquad \times [\bar{\alpha}(V_s(\varepsilon, n)) - \alpha_n(\bar{X}_s(\varepsilon), V_s(\varepsilon, n), Y_s(\varepsilon))] \rangle \, ds
$$

$$
- \Big[ \int_t^{\cdot} \Big( \frac{p_{m(n)}^{1/2}(\bar{X}_s(\varepsilon), V_s(\varepsilon, n))}{p_{m(n)}^{1/2}(\bar{X}_s(\varepsilon), Y_s(\varepsilon))}
$$

$$
\times \hat{W}_s(\varepsilon, n) \sigma(\bar{X}_s(\varepsilon), V_s(\varepsilon, n)) - \hat{Z}_s(\varepsilon) \sigma(\bar{X}_s(\varepsilon), Y_s(\varepsilon)) \Big) \, dB_s
$$

$$
+ \int_t^{\cdot} S_s(2, \varepsilon) \, dB_s \Big]_T
$$

$$- 2 \int_t^T \Big\langle V_s(\varepsilon, n) - \hat{Y}_s(\varepsilon),$$

$$\Big( \frac{p_{m(n)}^{1/2}(\bar{X}_s(\varepsilon), V_s(\varepsilon, n))}{p_{m(n)}^{1/2}(\bar{X}_s(\varepsilon), Y_s(\varepsilon))} \hat{W}_s(\varepsilon, n) \sigma(\bar{X}_s(\varepsilon), V_s(\varepsilon, n))$$

$$- \hat{Z}_s(\varepsilon) \sigma(\bar{X}_s(\varepsilon), Y_s(\varepsilon)) \Big) dB_s \Big\rangle$$

$$- 2 \int_t^T \langle V_s(\varepsilon, n) - \hat{Y}_s(\varepsilon), dS_s(\varepsilon) \rangle.$$

Thanks to the uniform boundedness of the processes $(\hat{Y}(\varepsilon))_{\varepsilon>0}$ and $(V(\varepsilon, n))_{\varepsilon>0, n\in\mathbb{N}^*}$, to the uniform boundedness of the quantities $\mathbf{E}\int_{t_0}^T (|\hat{Z}_t(\varepsilon)|^2 + |\hat{W}_t(\varepsilon, n)|^2)\, dt$ and to Section 3.1, the influence of the terms $S(\varepsilon)$ and $S(2, \varepsilon)$ is negligible. Therefore, for all $\varepsilon > 0$, $n \in \mathbb{N}^*$ and $t \in [t_0, T]$,

$$\mathbf{E}|V_t(\varepsilon, n) - \hat{Y}_t(\varepsilon)|^2$$

$$+ \mathbf{E} \int_t^T \Big| \frac{p_{m(n)}^{1/2}(\bar{X}_s(\varepsilon), V_s(\varepsilon, n))}{p_{m(n)}^{1/2}(\bar{X}_s(\varepsilon), Y_s(\varepsilon))} \hat{W}_s(\varepsilon, n) \sigma(\bar{X}_s(\varepsilon), V_s(\varepsilon, n))$$

$$- \hat{Z}_s(\varepsilon) \sigma(\bar{X}_s(\varepsilon), Y_s(\varepsilon)) \Big|^2 ds$$

$$\leq \delta_{4.2}(\varepsilon, n) + \mathbf{E}|H_n(U_T(\varepsilon)) - \hat{Y}_T(\varepsilon)|^2$$

$$+ 2\mathbf{E} \int_t^T \langle V_s(\varepsilon, n) - \hat{Y}_s(\varepsilon),$$

$$\bar{v}_n(V_s(\varepsilon, n), W_s(\varepsilon, n)) - v(\bar{X}_s(\varepsilon), Y_s(\varepsilon), \hat{Z}_s(\varepsilon)) \rangle\, ds$$

$$(4.9) \qquad + 2\mathbf{E} \int_t^T \langle V_s(\varepsilon, n) - \hat{Y}_s(\varepsilon), W_s(\varepsilon, n)[\bar{u}(V_s(\varepsilon, n), W_s(\varepsilon, n))$$

$$- u(\bar{X}_s(\varepsilon), Y_s(\varepsilon), \hat{Z}_s(\varepsilon))] \rangle\, ds$$

$$+ \mathbf{E} \int_t^T \langle V_s(\varepsilon, n) - \hat{Y}_s(\varepsilon),$$

$$\nabla^2_{x,x} \zeta_n(s, U_s(\varepsilon, n))$$

$$\times [\bar{\alpha}(V_s(\varepsilon, n)) - \alpha_n(\bar{X}_s(\varepsilon), V_s(\varepsilon, n), Y_s(\varepsilon))] \rangle\, ds$$

$$= \delta_{4.2}(\varepsilon, n) + \mathbf{E}|H_n(U_T(\varepsilon, n)) - \hat{Y}_T(\varepsilon)|^2$$

$$+ 2\Delta_{4.2.1}(t, \varepsilon, n) + 2\Delta_{4.2.2}(t, \varepsilon, n) + \Delta_{4.2.3}(t, \varepsilon, n),$$

where $(\delta_{4.2}(\varepsilon,n))_{\varepsilon>0,n\in\mathbb{N}^*}$ satisfies, for every $n\in\mathbb{N}^*$, $\lim_{\varepsilon\to0}\delta_{4.2}(\varepsilon,n)=0$ (the values of this family may change from one inequality to another). Note that for all $\varepsilon>0$, $n\in\mathbb{N}^*$ and $t\in[t_0,T]$,

$$
\begin{aligned}
&\Delta_{4.2.1}(t,\varepsilon,n) \\
&\quad\leq \mathbf{E}\int_t^T \langle V_s(\varepsilon,n)-\hat{Y}_s(\varepsilon), \\
&\qquad\qquad [\bar{v}_n(V_s(\varepsilon,n),W_s(\varepsilon,n))-\bar{v}(V_s(\varepsilon,n),W_s(\varepsilon,n))] \\
&\qquad\qquad + [\bar{v}(V_s(\varepsilon,n),W_s(\varepsilon,n))-\bar{v}(Y_s(\varepsilon),W_s(\varepsilon,n))] \\
&\qquad\qquad + [\bar{v}(Y_s(\varepsilon),W_s(\varepsilon,n)) \\
&\qquad\qquad\quad -v(\bar{X}_s(\varepsilon),Y_s(\varepsilon),W_s(\varepsilon,n)(I+\nabla_x\hat{b})(\bar{X}_s(\varepsilon),Y_s(\varepsilon)))] \\
&\qquad\qquad + [v(\bar{X}_s(\varepsilon),Y_s(\varepsilon),W_s(\varepsilon,n)(I+\nabla_x\hat{b})(\bar{X}_s(\varepsilon),Y_s(\varepsilon))) \\
&\qquad\qquad\qquad\qquad -v(\bar{X}_s(\varepsilon),Y_s(\varepsilon),\hat{Z}_s(\varepsilon))]\rangle\,ds,
\end{aligned}
\tag{4.10}
$$

$$
\begin{aligned}
&\Delta_{4.2.2}(t,\varepsilon,n) \\
&\quad\leq \mathbf{E}\int_t^T \langle V_s(\varepsilon,n)-\hat{Y}_s(\varepsilon), \\
&\qquad\qquad W_s(\varepsilon,n)\{[\bar{u}(V_s(\varepsilon,n),W_s(\varepsilon,n))-\bar{u}(Y_s(\varepsilon),W_s(\varepsilon,n))] \\
&\qquad\qquad + [\bar{u}(Y_s(\varepsilon),W_s(\varepsilon,n)) \\
&\qquad\qquad\quad -u(\bar{X}_s(\varepsilon),Y_s(\varepsilon), \\
&\qquad\qquad\qquad W_s(\varepsilon,n)(I+\nabla_x\hat{b})(\bar{X}_s(\varepsilon),Y_s(\varepsilon)))] \\
&\qquad\qquad + [u(\bar{X}_s(\varepsilon),Y_s(\varepsilon), \\
&\qquad\qquad\qquad W_s(\varepsilon,n)(I+\nabla_x\hat{b})(\bar{X}_s(\varepsilon),Y_s(\varepsilon))) \\
&\qquad\qquad\qquad\qquad -u(\bar{X}_s(\varepsilon),Y_s(\varepsilon),\hat{Z}_s(\varepsilon))]\}\rangle\,ds,
\end{aligned}
\tag{4.11}
$$

$$
\begin{aligned}
&\Delta_{4.2.3}(t,\varepsilon,n) \\
&\quad\leq \mathbf{E}\int_t^T \langle V_s(\varepsilon,n)-\hat{Y}_s(\varepsilon),\nabla^2_{x,x}\zeta_n(s,U_s(\varepsilon,n)) \\
&\qquad\qquad\quad \times \{[\bar{\alpha}(V_s(\varepsilon,n))-\bar{\alpha}_n(V_s(\varepsilon,n),Y_s(\varepsilon))] \\
&\qquad\qquad\quad + [\bar{\alpha}_n(V_s(\varepsilon,n),Y_s(\varepsilon)) \\
&\qquad\qquad\qquad -\alpha_n(\bar{X}_s(\varepsilon),V_s(\varepsilon,n),Y_s(\varepsilon))]\}\rangle\,ds.
\end{aligned}
\tag{4.12}
$$

Recall the following

1. $(\bar{v}_n-\bar{v})_{n\in\mathbb{N}^*}$ uniformly converges on every compact set toward 0 and $\nabla^2_{x,x}\zeta_n(\bar{\alpha}-\bar{\alpha}_n)$ uniformly converges on $[0,T]\times\mathbb{R}^P\times\bar{B}_Q(0,C_{3.1})\times\bar{B}_Q(0,C_{2.4})$ toward 0 [see (3.24)].

2. $\nabla_x\hat{b}$ is Lipschitzian with respect to $y$ (see Theorem 2.2), $u$ and $\bar{u}$ are locally Lipschitzian in $y$ and $z$, and $v$ and $\bar{v}$ are locally bounded, locally monotonous in $y$ and locally Lipschitzian in $z$ (see Section 2.5).

3. From Theorem 2.3 and Section 3.3, $(\hat{Y}(\varepsilon))_{\varepsilon>0}$, $(V(\varepsilon,n))_{\varepsilon>0,n\in\mathbb{N}^*}$, $(W(\varepsilon, n))_{\varepsilon>0,n\in\mathbb{N}^*}$ and $(\hat{W}(\varepsilon,n))_{\varepsilon>0,n\in\mathbb{N}^*}$ are bounded, uniformly in $\varepsilon$ and $n$.

4. $\forall \varepsilon>0$, $\hat{Y}(\varepsilon) = Y(\varepsilon) - \varepsilon\hat{e}$.

Hence, applying Theorem 2.4 [with $G = U(\varepsilon,n)$], there exists a constant $C_{4.2}$ (whose value may change from one inequality to another), only depending on $k$, $K$, $\lambda$, $\Lambda$, $P$, $Q$ and $T$, such that for all $\varepsilon>0$, $n\in\mathbb{N}^*$ and $t\in[t_0,T]$,

$$
\begin{aligned}
\Delta_{4.2.1}&(t,\varepsilon,n) + \Delta_{4.2.2}(t,\varepsilon,n) + \Delta_{4.2.3}(t,\varepsilon,n) \\
&\leq \delta_{4.2}(\varepsilon,n) + \delta_{4.3}(n) \\
&\quad + C_{4.2}\mathbf{E}\int_t^T |V_s(\varepsilon,n) - Y_s(\varepsilon)| \\
&\qquad\qquad \times [|V_s(\varepsilon,n) - Y_s(\varepsilon)| \\
&\qquad\qquad\qquad + |W_s(\varepsilon,n)(I + \nabla_x\hat{b})(\bar{X}_s(\varepsilon),Y_s(\varepsilon)) - \hat{Z}_s(\varepsilon)|]\,ds \\
&\leq \delta_{4.2}(\varepsilon,n) \\
&\quad + \delta_{4.3}(n) + C_{4.2}\mathbf{E}\int_t^T |V_s(\varepsilon,n) - \hat{Y}_s(\varepsilon)| \\
&\qquad\qquad \times [|V_s(\varepsilon,n) - Y_s(\varepsilon)| + |\hat{W}_s(\varepsilon,n) - \hat{Z}_s(\varepsilon)|]\,ds,
\end{aligned}
\tag{4.13}
$$

where $\lim_{n\to+\infty}\delta_{4.3}(n) = 0$. Note that the values of $(\delta_{4.3}(n))_{n\in\mathbb{N}^*}$ may change from one line to another. Moreover, thanks to the ellipticity of $\sigma$, note that for all $\varepsilon>0$, $n\in\mathbb{N}^*$ and $t\in[t_0,T]$,

$$
\begin{aligned}
\mathbf{E}&\int_t^T |\hat{W}_s(\varepsilon,n) - \hat{Z}_s(\varepsilon)|^2\,ds \\
&\leq C_{4.2}\mathbf{E}\int_t^T \left|\frac{p_{m(n)}^{1/2}(\bar{X}_s(\varepsilon),V_s(\varepsilon,n))}{p_{m(n)}^{1/2}(\bar{X}_s(\varepsilon),Y_s(\varepsilon))}\hat{W}_s(\varepsilon,n)\sigma(\bar{X}_s(\varepsilon),V_s(\varepsilon,n))\right. \\
&\qquad\qquad\qquad\qquad\qquad\qquad \left. - \hat{Z}_s(\varepsilon)\sigma(\bar{X}_s(\varepsilon),Y_s(\varepsilon))\right|^2\,ds \\
&\quad + C_{4.2}\mathbf{E}\int_t^T \left|\frac{p_{m(n)}^{1/2}(\bar{X}_s(\varepsilon),V_s(\varepsilon,n))}{p_{m(n)}^{1/2}(\bar{X}_s(\varepsilon),Y_s(\varepsilon))}\hat{W}_s(\varepsilon,n)\sigma(\bar{X}_s(\varepsilon),V_s(\varepsilon,n))\right. \\
&\qquad\qquad\qquad\qquad\qquad\qquad \left. - \hat{W}_s(\varepsilon,n)\sigma(\bar{X}_s(\varepsilon),Y_s(\varepsilon))\right|^2\,ds.
\end{aligned}
\tag{4.14}
$$

Following (4.4)–(4.6), we deduce from the boundedness of $(\hat{W}(\varepsilon,n))_{\varepsilon>0,n\in\mathbb{N}^*}$ (see Section 3.3) that for all $\varepsilon>0$, $n\in\mathbb{N}^*$ and $t\in[t_0,T]$,

$$
\mathbf{E}\int_t^T \left|\frac{p_{m(n)}^{1/2}(\bar{X}_s(\varepsilon),V_s(\varepsilon,n))}{p_{m(n)}^{1/2}(\bar{X}_s(\varepsilon),Y_s(\varepsilon))}\hat{W}_s(\varepsilon,n)\sigma(\bar{X}_s(\varepsilon),V_s(\varepsilon,n))\right.
$$

$$(4.15) \qquad\qquad\qquad\qquad - \hat{W}_s(\varepsilon, n) \sigma(\bar{X}_s(\varepsilon), Y_s(\varepsilon)) \Big|^2 \, ds$$

$$\leq \delta_{4.2}(\varepsilon, n) + C_{4.2} \mathbf{E} \int_t^T |V_s(\varepsilon, n) - Y_s(\varepsilon)|^2 \, ds.$$

Hence, from (4.14) and (4.15), we deduce that for all $\varepsilon > 0$, $n \in \mathbb{N}^*$ and $t \in [t_0, T]$,

$$\mathbf{E} \int_t^T |\hat{W}_s(\varepsilon, n) - \hat{Z}_s(\varepsilon)|^2 \, ds$$

$$\leq \delta_{4.2}(\varepsilon, n) + C_{4.2} \Big[ \mathbf{E} \int_t^T \Big| \frac{p_{m(n)}^{1/2}(\bar{X}_s(\varepsilon), V_s(\varepsilon, n))}{p_{m(n)}^{1/2}(\bar{X}_s(\varepsilon), Y_s(\varepsilon))}$$

$$(4.16) \qquad\qquad\qquad \times \hat{W}_s(\varepsilon, n) \sigma(\bar{X}_s(\varepsilon), V_s(\varepsilon, n))$$

$$- \hat{Z}_s(\varepsilon) \sigma(\bar{X}_s(\varepsilon), Y_s(\varepsilon)) \Big|^2 \, ds$$

$$+ \mathbf{E} \int_t^T |V_s(\varepsilon, n) - Y_s(\varepsilon)|^2 \, ds \Big].$$

Finally, we deduce from (4.9), (4.13) and (4.16), that for all $\varepsilon > 0$, $n \in \mathbb{N}^*$ and $t \in [t_0, T]$,

$$\mathbf{E}|V_t(\varepsilon, n) - Y_t(\varepsilon)|^2 + \mathbf{E} \int_t^T |\hat{W}_s(\varepsilon, n) - \hat{Z}_s(\varepsilon)|^2 \, ds$$

$$(4.17) \qquad \leq \delta_{4.2}(\varepsilon, n) + \delta_{4.3}(n) + \mathbf{E}|H_n(U_T(\varepsilon, n)) - H(X_T(\varepsilon))|^2$$

$$+ C_{4.2} \mathbf{E} \int_t^T |V_s(\varepsilon, n) - Y_s(\varepsilon)|^2 \, ds.$$

Recall that $H$ is Lipschitzian [see (H.3)] and that $(H_n)_{n \in \mathbb{N}^*}$ uniformly converges toward $H$. Hence, from Proposition 4.2,

$$\sup_{t_0 \leq t \leq T} \mathbf{E}|V_t(\varepsilon, n) - Y_t(\varepsilon)|^2 + \mathbf{E} \int_{t_0}^T |\hat{W}_t(\varepsilon, n) - \hat{Z}_t(\varepsilon)|^2 \, dt$$

$$(4.18) \qquad \leq \delta_{4.2}(\varepsilon, n) + \delta_{4.3}(n) + C_{4.2} \mathbf{E} \int_{t_0}^T |V_t(\varepsilon, n) - Y_t(\varepsilon)|^2 \, dt.$$

In fact, applying the Burkholder–Davis–Gundy inequalities, the same scheme leads to [note that we had to estimate first the term $\mathbf{E} \int_{t_0}^T |\hat{W}_t(\varepsilon, n) - \hat{Z}_t(\varepsilon)|^2 \, dt$ to apply the Burkholder–Davis–Gundy inequalities]

$$\mathbf{E} \sup_{t_0 \leq t \leq T} |V_t(\varepsilon, n) - Y_t(\varepsilon)|^2 + \mathbf{E} \int_{t_0}^T |\hat{W}_t(\varepsilon, n) - \hat{Z}_t(\varepsilon)|^2 \, dt$$

$$(4.19)$$

$$\leq \delta_{4.2}(\varepsilon, n) + \delta_{4.3}(n) + C_{4.2} \mathbf{E} \int_{t_0}^{T} |V_t(\varepsilon, n) - Y_t(\varepsilon)|^2 \, dt.$$

Thanks to Proposition 4.2, we complete the proof. $\square$

4.3. *Proof of Theorem* 3.1. We now establish Theorem 3.1. From the definitions of $(V(\varepsilon, n))_{\varepsilon > 0, n \in \mathbb{N}^*}$ and $(W(\varepsilon, n))_{\varepsilon > 0, n \in \mathbb{N}^*}$, we deduce from Proposition 4.3 that for $T - t_0 \leq c_{4.2}$ and for all $\varepsilon > 0$ and $n \in \mathbb{N}^*$,

$$
\begin{aligned}
&\mathbf{E} \sup_{t_0 \leq t \leq T} |U_t(\varepsilon, n) - X_t(\varepsilon)|^2 + \mathbf{E} \sup_{t_0 \leq t \leq T} |\zeta_n(t, U_t(\varepsilon, n)) - Y_t(\varepsilon)|^2 \\
(4.20) \quad &+ \mathbf{E} \int_{t_0}^{T} |\nabla_x \zeta_n(t, U_t(\varepsilon, n))(I + \nabla_x \hat{b})(\bar{X}_t(\varepsilon), \zeta_n(t, U_t(\varepsilon, n))) - \hat{Z}_t(\varepsilon)|^2 \, dt \\
&\leq \delta_{4.2}(\varepsilon, n) + \delta_{4.3}(n).
\end{aligned}
$$

Recall from Section 2.5 that $\theta$ is $C_{2.6}$-Lipschitzian with respect to the variable $x$. Hence, we claim that for all $\varepsilon > 0$, $n \in \mathbb{N}^*$ and $M > 0$,

$$
\begin{aligned}
&\mathbf{E} \sup_{t_0 \leq t \leq T} |\theta(t, X_t(\varepsilon)) - Y_t(\varepsilon)|^2 \\
&\quad \leq 3 \mathbf{E} \sup_{t_0 \leq t \leq T} \big[ |\theta(t, X_t(\varepsilon)) - \theta(t, U_t(\varepsilon, n))|^2 \\
&\qquad\qquad\qquad + |\theta(t, U_t(\varepsilon, n)) - \zeta_n(t, U_t(\varepsilon, n))|^2 \\
&\qquad\qquad\qquad + |\zeta_n(t, U_t(\varepsilon, n)) - Y_t(\varepsilon)|^2 \big] \\
&\quad \leq 3 C_{2.6}^2 \mathbf{E} \sup_{t_0 \leq t \leq T} |X_t(\varepsilon) - U_t(\varepsilon, n)|^2 \\
(4.21) \quad &\qquad + 3 \sup_{[t_0, T] \times \bar{B}_P(0, M)} |(\theta - \zeta_n)(t, x)|^2 \\
&\qquad + 3 \|\theta - \zeta_n\|_\infty^2 \mathbf{P} \Big\{ \sup_{t_0 \leq t \leq T} |U_t(\varepsilon, n)| > M \Big\} + \delta_{4.2}(\varepsilon, n) + \delta_{4.3}(n) \\
&\quad \leq 3 \sup_{[t_0, T] \times \bar{B}_P(0, M)} |(\theta - \zeta_n)(t, x)|^2 \\
&\qquad + 3 \|\theta - \zeta_n\|_\infty^2 \mathbf{P} \Big\{ \sup_{t_0 \leq t \leq T} |U_t(\varepsilon, n)| > M \Big\} + \delta_{4.2}(\varepsilon, n) + \delta_{4.3}(n).
\end{aligned}
$$

Thanks to Section 3.3, $(\zeta_n)_{n \in \mathbb{N}^*}$ is uniformly bounded on $[0, T] \times \mathbb{R}^P$ and uniformly converges on every compact subset of $[0, T] \times \mathbb{R}^P$ toward $\theta$. From Lemma 4.1, $\sup_{\varepsilon > 0, n \in \mathbb{N}^*} \mathbf{E} \sup_{t_0 \leq t \leq T} |U_t(\varepsilon, n)|^2$ is finite. Therefore,

$$(4.22) \qquad \lim_{\varepsilon \to 0} \mathbf{E} \sup_{t_0 \leq t \leq T} |\theta(t, X_t(\varepsilon)) - Y_t(\varepsilon)|^2 = 0.$$

Let us now deal with the gradient term appearing in (3.18). For all $\varepsilon > 0$, $n \in \mathbb{N}^*$ and $0 < \eta < T - t_0$,

$$\mathbf{E} \int_{t_0}^{T} |\nabla_x \theta(t, X_t(\varepsilon))(I + \nabla_x \hat{b})(\bar{X}_t(\varepsilon), Y_t(\varepsilon)) - \hat{Z}_t(\varepsilon)|^2 \, dt$$

$$\leq 3\mathbf{E}\int_{t_0}^{T-\eta}|\nabla_x\theta(t,X_t(\varepsilon))(I+\nabla_x\hat{b})(\bar{X}_t(\varepsilon),Y_t(\varepsilon))$$
$$-\nabla_x\theta(t,U_t(\varepsilon,n))(I+\nabla_x\hat{b})(\bar{X}_t(\varepsilon),\theta(t,U_t(\varepsilon,n)))|^2\,dt$$
$$+3\mathbf{E}\int_{t_0}^{T-\eta}|\nabla_x\theta(t,U_t(\varepsilon,n))(I+\nabla_x\hat{b})(\bar{X}_t(\varepsilon),\theta(t,U_t(\varepsilon,n)))$$
$$\text{(4.23)}\qquad\qquad -\nabla_x\zeta_n(t,U_t(\varepsilon,n))(I+\nabla_x\hat{b})(\bar{X}_t(\varepsilon),\zeta_n(t,U_t(\varepsilon,n)))|^2\,dt$$
$$+3\mathbf{E}\int_{t_0}^{T-\eta}|\nabla_x\zeta_n(t,U_t(\varepsilon,n))(I+\nabla_x\hat{b})(\bar{X}_t(\varepsilon),\zeta_n(t,U_t(\varepsilon,n)))$$
$$-\hat{Z}_t(\varepsilon)|^2\,dt$$
$$+\mathbf{E}\int_{T-\eta}^{T}|\nabla_x\theta(t,X_t(\varepsilon))(I+\nabla_x\hat{b})(\bar{X}_t(\varepsilon),Y_t(\varepsilon))-\hat{Z}_t(\varepsilon)|^2\,dt$$
$$=\Delta_{4.3.1}(\varepsilon,n,\eta)+\Delta_{4.3.2}(\varepsilon,n,\eta)+\Delta_{4.3.3}(\varepsilon,n,\eta)+\Delta_{4.3.4}(\varepsilon,n,\eta).$$

Thanks to Theorem 2.2, $\nabla_x\hat{b}$ is bounded and Lipschitzian with respect to $y$. From Section 2.5, $\nabla_x\theta$ is bounded on $[0,T[\times\mathbb{R}^P$ and Hölderian on $[t_0,T-\eta]\times\mathbb{R}^P$ for every $0<\eta<T-t_0$. Hence, thanks to (4.3.1),

$$\text{(4.24)}\quad\forall\varepsilon>0,\ \forall n\in\mathbb{N}^*,\ \forall\eta>0,\qquad\Delta_{4.3.1}(\varepsilon,n,\eta)\leq\delta_{4.4}(n,\eta)+\delta_{4.5}(\varepsilon,n,\eta),$$

where for all $0<\eta<T-t_0$, $\lim_{n\to+\infty}\delta_{4.4}(n,\eta)=0$, and for all $n\in\mathbb{N}^*$, $\lim_{\varepsilon\to0}\delta_{4.5}(\varepsilon,n,\eta)=0$.

Recall now from Section 3.3 that the sequences of functions $(\zeta_n)_{n\in\mathbb{N}^*}$ and $(\nabla_x\zeta_n)_{n\in\mathbb{N}^*}$ are bounded on $[0,T]\times\mathbb{R}^P$ by $C_{3.1}$ and uniformly converge on every compact subset of $[0,T[\times\mathbb{R}^P$ toward $\theta$ and $\nabla_x\theta$. Hence, using Lemma 4.1 [follow (4.21)],

$$\text{(4.25)}\qquad\forall\varepsilon>0,\ \forall n\in\mathbb{N}^*,\ \forall\eta>0,\qquad\Delta_{4.3.2}(\varepsilon,n,\eta)\leq\delta_{4.4}(n,\eta).$$

Note now that $\Delta_{4.3.4}(\varepsilon,n,\eta)\leq\delta_{4.6}(\eta)$, with $\lim_{\eta\to0}\delta_{4.6}(\eta)=0$. Therefore, thanks to (4.20), (4.23)–(4.25),

$$\mathbf{E}\int_{t_0}^{T}|\nabla_x\theta(t,X_t(\varepsilon))(I+\nabla_x\hat{b})(\bar{X}_t(\varepsilon),\theta(t,X_t(\varepsilon)))-\hat{Z}_t(\varepsilon)|^2\,dt$$
$$\text{(4.26)}\qquad\qquad\leq\delta_{4.4}(\varepsilon,\eta,n)+\delta_{4.5}(\eta,n)+\delta_{4.6}(\eta).$$

Thanks to (4.22) and (4.26), we deduce that Theorem 3.1 holds with $c_{3.1}=c_{4.2}$.

4.4. *Proof of the homogenization property.* We first establish point 2 of Theorem (PHP). To this end, we wish to extend our local convergence result to the whole interval $[0,T]$. For this purpose, we apply the scheme developed on the interval $[t_0,T]$ to a new interval $[t_0-\delta,t_0]$ with $\delta>0$, and then conclude by induction. Actually, the main difference between these two cases comes from the new final condition $\theta(t_0,\cdot)$.

Indeed, in the former sections, the final conditions $(Y_T(\varepsilon))_{\varepsilon>0}$ were given by the relationships $\forall \varepsilon > 0$, $Y_T(\varepsilon) = H(X_T(\varepsilon))$, where $H$ is $k$-Lipschitzian, but, from Theorem 3.1, we just know that the new final conditions $(Y_{t_0}(\varepsilon))_{\varepsilon>0}$ satisfy

$$(4.27) \qquad \mathbf{E}|Y_{t_0}(\varepsilon) - \theta(t_0, X_{t_0}(\varepsilon))|^2 \to 0 \qquad \text{as } \varepsilon \to 0,$$

where $\theta(t_0, \cdot)$ is, thanks to (2.12), $C_{2.6}$-Lipschitzian.

Actually, property (4.27) is sufficient to apply the strategy developed in Sections 4.2 and 4.3. In particular, modifying $c_{4.2}$ if necessary, we can prove that (3.18) holds on the interval $[t_0 - c_{4.2}, t_0]$, and therefore, on the interval $[t_0 - c_{4.2}, T]$. Following the induction method introduced in Delarue (2002b), we then prove that it holds on $[0, T]$. In particular, since we control the Lipschitz constant of $\theta$ with respect to the variable $x$ [see again (2.12)], we know that the "length $c_{4.2}$" is bounded from below by a nonnegative real during the induction. This completes the proof of point 2 of Theorem (PHP).

We directly deduce $\lim_{\varepsilon \to 0} \theta_\varepsilon(0, x_0) = \theta(0, x_0)$. This proves point 3 of Theorem (HP).

4.5. *Proof of point 3 of Theorem* (PHP). In this final section we prove point 3 of Theorem (PHP). Basically, it can be seen as a direct consequence of point 2. We first set for all $\varepsilon > 0$ and $t \in [0, T]$,

$$(4.28) \qquad \begin{aligned} N_t(\varepsilon) &= \int_0^t [(I + \nabla_x \hat{b})\sigma](\bar{X}_s(\varepsilon), Y_s(\varepsilon)) \, dB_s, \\ M_t(\varepsilon) &= \int_0^t \hat{Z}_s(\varepsilon)\sigma(\bar{X}_s(\varepsilon), Y_s(\varepsilon)) \, dB_s. \end{aligned}$$

Note that for every $\varepsilon > 0$ and for every $t \in [0, T]$,

$$(4.29) \qquad \begin{aligned} [N(\varepsilon)]_t &= \int_0^t \alpha(\bar{X}_s(\varepsilon), Y_s(\varepsilon)) \, dB_s, \\ [M(\varepsilon)]_t &= \int_0^t \hat{Z}_s(\varepsilon)a(\bar{X}_s(\varepsilon), Y_s(\varepsilon))\hat{Z}_s^*(\varepsilon) \, ds, \\ [M(\varepsilon), N(\varepsilon)]_t &= \int_0^t \hat{Z}_s(\varepsilon)a(\bar{X}_s(\varepsilon), Y_s(\varepsilon))(I + \nabla_x \hat{b})^*(\bar{X}_s(\varepsilon), Y_s(\varepsilon)) \, ds. \end{aligned}$$

Note from Theorem 2.2 that $(I + \nabla_x \hat{b})\sigma$ is bounded and, therefore, that

$$\sup_{\varepsilon > 0} \mathbf{E} \sup_{0 \le t \le T} |N_t(\varepsilon)|^2 < \infty.$$

Moreover, note from Theorem 2.3 that

$$\sup_{\varepsilon > 0} \mathbf{E} \sup_{0 \le t \le T} |M_t(\varepsilon)|^2 < \infty.$$

Thanks to this notation, recall from Section 3.1 that for all $\varepsilon > 0$ and $t \in [0, T]$,

$$
\begin{aligned}
X_t(\varepsilon) &+ \varepsilon \hat{b}\big( \bar{X}_t(\varepsilon), Y_t(\varepsilon) \big) \\
&= x_0 + \varepsilon \hat{b}\Big( \frac{x_0}{\varepsilon}, Y_0(\varepsilon) \Big) \\
&\quad + \int_0^t u\big( \bar{X}_s(\varepsilon), Y_s(\varepsilon), \hat{Z}_s(\varepsilon) \big) \, ds + N_t(\varepsilon) + R_t(\varepsilon),
\end{aligned}
$$

(4.30)
$$
\begin{aligned}
Y_t(\varepsilon) &- \varepsilon \hat{e}\big( \bar{X}_t(\varepsilon), Y_t(\varepsilon) \big) \\
&= H(X_T(\varepsilon)) - \varepsilon \hat{e}\big( \bar{X}_T(\varepsilon), Y_T(\varepsilon) \big) \\
&\quad + \int_t^T v\big( \bar{X}_s(\varepsilon), Y_s(\varepsilon), \hat{Z}_s(\varepsilon) \big) \, ds \\
&\qquad - (M_T(\varepsilon) - M_t(\varepsilon)) + S_T(\varepsilon) - S_t(\varepsilon).
\end{aligned}
$$

We now claim:

LEMMA 4.4. *The processes $(X(\varepsilon), M(\varepsilon))_{\varepsilon > 0}$ are tight in the space $(\mathcal{C}([0, T], \mathbb{R}^P \times \mathbb{R}^Q), \| \cdot \|_\infty)$.*

PROOF. From (4.30), Section 3.1 [which permits us to estimate $(R(\varepsilon))_{\varepsilon > 0}$], Theorem 2.3 [which provides estimates of $(X(\varepsilon))_{\varepsilon > 0}$ and $(\hat{Z}(\varepsilon))_{\varepsilon > 0}$] and Kolmogorov's criterium, there exists a family of continuous and $\mathbb{R}^P$-valued processes $(r(\varepsilon))_{\varepsilon > 0}$ satisfying $\sup_{0 \le t \le T} |r_t(\varepsilon)| \to 0$ in probability as $\varepsilon \to 0$ such that the family $(X(\varepsilon) + r(\varepsilon))_{\varepsilon > 0}$ is tight in the space $(\mathcal{C}([0, T], \mathbb{R}^P), \| \cdot \|_\infty)$.

Moreover, it is readily seen that for every $\varepsilon_0 > 0$, the family of processes $(X(\varepsilon))_{\varepsilon \ge \varepsilon_0}$ is tight in the same space as above. We easily deduce that $(X(\varepsilon))_{\varepsilon > 0}$ is tight.

In particular, the family $(\theta(\cdot, X.(\varepsilon)))_{\varepsilon > 0}$ is tight in the space $(\mathcal{C}([0, T], \mathbb{R}^Q), \| \cdot \|_\infty)$. Moreover, referring to Delarue (2002b), for every $\varepsilon_0 > 0$, the family $(Y(\varepsilon))_{\varepsilon > \varepsilon_0}$ is tight in the same space as above (roughly speaking, we know that for every $\varepsilon > \varepsilon_0$, for every $t \in [0, T]$, $Y_t(\varepsilon) = \theta_\varepsilon(t, X_t(\varepsilon))$, where the family $(\theta_\varepsilon : (t, x) \in [0, T] \times \mathbb{R}^P \mapsto Y_t(\varepsilon, t, x))_{\varepsilon > \varepsilon_0}$ is equicontinuous). Hence, we deduce from point 2 of Theorem (PHP) that the family $(Y(\varepsilon))_{\varepsilon > 0}$ is tight.

Finally, using again the same argument, we deduce from (4.30), Section 3.1 [which permits estimate $(S(\varepsilon))_{\varepsilon > 0}$], Theorem 2.3 [which permits to estimate $(\hat{Z}(\varepsilon))_{\varepsilon > 0}$] and Kolmogorov's citerium [which permits us to prove that the family of processes $(\int_0^\cdot v(\bar{X}_s(\varepsilon), Y_s(\varepsilon), \hat{Z}_s(\varepsilon)) \, ds)_{\varepsilon > 0}$ is tight] that the family of processes $(M(\varepsilon))_{\varepsilon > 0}$ is tight in the space $(\mathcal{C}([0, T], \mathbb{R}^Q), \| \cdot \|_\infty)$. This completes                              the                              proof.     □

Let us now complete the proof of point 3 of Theorem (PHP).

PROOF OF POINT 3 OF THEOREM (PHP). Thanks to Lemma 4.5, there exist a subsequence still indexed by $\varepsilon > 0$ and a continuous and $\mathbb{R}^P \times \mathbb{R}^Q$ valued process, denoted by $(X_t, M_t)_{0 \le t \le T}$, such that,

$$(4.31) \quad (X_t(\varepsilon), M_t(\varepsilon))_{0 \le t \le T} \implies (X_t, M_t)_{0 \le t \le T} \qquad \text{as } \varepsilon \to 0,$$

where $\implies$ denotes the convergence in law on the space $\mathcal{C}([0, T], \mathbb{R}^P \times \mathbb{R}^Q)$ endowed with the topology of the uniform convergence.

Applying estimates given in Sections 2.3 and 2.5 and using point 2 of Theorem (PHP), note that

$$
\begin{aligned}
(4.32) \quad \sup_{0 \le t \le T} \Bigg| &\int_0^t u(\bar{X}_s(\varepsilon), Y_s(\varepsilon), \hat{Z}_s(\varepsilon))\, ds \\
&- \int_0^t u(\bar{X}_s(\varepsilon), Y_s(\varepsilon), \\
&\qquad \nabla_x \theta(s, X_s(\varepsilon))(I + \nabla_x \hat{b})(\bar{X}_s(\varepsilon), Y_s(\varepsilon)))\, ds \Bigg| \xrightarrow{\mathbf{P}} 0 \\
&\hspace{8cm} \text{as } \varepsilon \to 0.
\end{aligned}
$$

Hence, from point 2 of Theorem 2.4 [to treat the term $\nabla_x \theta(\cdot, X.(\varepsilon))$], recall that $\nabla_x \theta$ is bounded and locally Hölder continuous, approximate $X(\varepsilon)$ by $\hat{X}(\varepsilon)$ and set $G = \hat{X}(\varepsilon)$ in Theorem 2.4. We deduce that

$$
\begin{aligned}
(4.33) \quad \sup_{0 \le t \le T} \Bigg| &\int_0^t u(\bar{X}_s(\varepsilon), Y_s(\varepsilon), \hat{Z}_s(\varepsilon))\, ds \\
&- \int_0^t \bar{u}(Y_s(\varepsilon), \nabla_x \theta(s, X_s(\varepsilon)))\, ds \Bigg| \xrightarrow{\mathbf{P}} 0 \qquad \text{as } \varepsilon \to 0.
\end{aligned}
$$

Therefore, using once again estimates given in Sections 2.3 and 2.5 and applying point 2 of Theorem (PHP), we claim that

$$
\begin{aligned}
(4.34) \quad \sup_{0 \le t \le T} \Bigg| &\int_0^t u(\bar{X}_s(\varepsilon), Y_s(\varepsilon), \hat{Z}_s(\varepsilon))\, ds \\
&- \int_0^t \bar{u}(\theta(s, X_s(\varepsilon)), \nabla_x \theta(s, X_s(\varepsilon)))\, ds \Bigg| \xrightarrow{\mathbf{P}} 0 \qquad \text{as } \varepsilon \to 0.
\end{aligned}
$$

In the same way, we deduce that

$$
\begin{aligned}
(4.35) \quad \sup_{0 \le t \le T} \Bigg| &\int_0^t v(\bar{X}_s(\varepsilon), Y_s(\varepsilon), \hat{Z}_s(\varepsilon))\, ds \\
&- \int_0^t \bar{v}(\theta(s, X_s(\varepsilon)), \nabla_x \theta(s, X_s(\varepsilon)))\, ds \Bigg| \xrightarrow{\mathbf{P}} 0 \qquad \text{as } \varepsilon \to 0,
\end{aligned}
$$

and from Theorem 2.3 [in particular, from the estimate of $(\hat{Z}(\varepsilon))_{\varepsilon>0}$] that

$$\sup_{0 \le t \le T} \left| [N(\varepsilon)]_t - \int_0^t \bar{\alpha}(\theta(\theta(s, X_s(\varepsilon)))) \, ds \right| \xrightarrow{\mathbf{P}} 0 \qquad \text{as } \varepsilon \to 0,$$

$$\sup_{0 \le t \le T} \left| [M(\varepsilon)]_t - \int_0^t \nabla_x \theta(s, X_s(\varepsilon)) \right.$$

$$(4.36) \qquad\qquad \left. \times \bar{\alpha}(\theta(s, X_s(\varepsilon)))(\nabla_x \theta(s, X_s(\varepsilon)))^* \, ds \right| \xrightarrow{\mathbf{P}} 0 \qquad \text{as } \varepsilon \to 0,$$

$$\sup_{0 \le t \le T} \left| [M(\varepsilon), N(\varepsilon)]_t - \int_0^t \nabla_x \theta(s, X_s(\varepsilon)) \bar{\alpha}(\theta(s, X_s(\varepsilon))) \, ds \right| \xrightarrow{\mathbf{P}} 0 \qquad \text{as } \varepsilon \to 0.$$

For the sake of simplicity, we set for every $t \in [0, T]$, $Y_t = \theta(t, X_t)$ and for every $t \in [0, T[$, $Z_t = \nabla_x \theta(t, X_t)$.

From (4.34)–(4.36), we deduce that

$$\left( \left( X_t(\varepsilon), Y_t(\varepsilon), M_t(\varepsilon), \int_0^t u(\bar{X}_s(\varepsilon), Y_s(\varepsilon), \hat{Z}_s(\varepsilon)) \, ds, \right. \right.$$

$$\int_0^t v(\bar{X}_s(\varepsilon), Y_s(\varepsilon), \hat{Z}_s(\varepsilon)) \, ds,$$

$$(4.37) \qquad\qquad \left. [N(\varepsilon)]_t, [M(\varepsilon)]_t, [M(\varepsilon), N(\varepsilon)]_t \right)_{0 \le t \le T} \Bigg)_{\varepsilon > 0}$$

$$\implies \left( X_t, Y_t, M_t, \int_0^t \bar{u}(Y_s, Z_s) \, ds, \int_0^t \bar{v}(Y_s, Z_s) \, ds, \right.$$

$$\left. \int_0^t \bar{\alpha}(Y_s) \, ds, \int_0^t Z_s \bar{\alpha}(Y_s) Z_s^* \, ds, \int_0^t Z_s \bar{\alpha}(Y_s) \, ds \right)_{0 \le t \le T},$$

where $\implies$ denotes the convergence in law on $(\mathcal{C}([0, T], \mathbb{R}^D), \|\cdot\|_\infty)$ ($D$ denotes an appropriate integer). Referring to Jacod and Shiryaev (1987), Chapter IX, Theorem 2.4, we deduce from Section 3.1, (4.29) [from which we know that $\sup_{\varepsilon>0} \mathbf{E} \sup_{0 \le t \le T}(|N_t(\varepsilon)|^2 + |M_t(\varepsilon)|^2) < \infty$], (4.30) and (4.37) that for every $0 \le t \le T$,

$$X_t = x_0 + \int_0^t \bar{u}(Y_s, Z_s) \, ds + N_t,$$

$$(4.38)$$

$$Y_t = H(X_T) + \int_t^T \bar{v}(Y_s, Z_s) \, ds - (M_T - M_t),$$

where $(N_t)_{0 \le t \le T}$ is an $\{\mathcal{F}_t^X\}_{0 \le t \le T}$-continuous square integrable martingale (recall that $\{\mathcal{F}_t^X\}_{0 \le t \le T}$ denotes the usual augmentation of the natural filtration of $X$) whose quadratic variation is given by

$$(4.39) \qquad \forall t \in [0, T], \qquad [N]_t = \int_0^t \bar{\alpha}(Y_s) \, ds.$$

Moreover, $(M_t)_{0\le t\le T}$ is also an $\{\mathcal{F}_t^X\}_{0\le t\le T}$-continuous square integrable martingale [note from (4.38) that $(M_t)_{0\le t\le T}$ is clearly $\{\mathcal{F}_t^X\}_{0\le t\le T}$-adapted] such that

$$(4.40) \qquad \forall t \in [0, T], \qquad [M]_t = \int_0^t Z_s \bar{\alpha}(Y_s) Z_s^* \, ds,$$

$$[M, N]_t = \int_0^t Z_s \bar{\alpha}(Y_s) \, ds.$$

In particular, $\forall 0 \le t \le T$, $[M - \int_0^\cdot Z_s \, dN_s]_t = 0$. Therefore, for every $t \in [0, T]$, $M_t = \int_0^t Z_s \, dN_s$. We now set $\forall t \in [0, T]$, $\bar{B}_t = \int_0^t \bar{\alpha}^{-1/2}(Y_s) \, dN_s$.

It is clear that $(\bar{B}_t)_{0\le t\le T}$ is an $\{\mathcal{F}_t^X\}_{0\le t\le T}$-Brownian motion. Moreover, we deduce from (4.38) that $(X, Y, Z)$ satisfies the FBSDE

$$(4.41) \qquad X_t = x_0 + \int_0^t \bar{u}(Y_s, Z_s) \, ds + \int_0^t \bar{\alpha}^{1/2}(Y_s) \, d\bar{B}_s,$$

$$Y_t = H(X_T) + \int_t^T \bar{v}(Y_s, Z_s) \, ds - \int_t^T Z_s \bar{\alpha}^{1/2}(Y_s) \, d\bar{B}_s.$$

Thanks to Delarue (2002b) (see Remarks 1.6 and 2.7), $(X, Y, Z)$ and $(X(0, x_0)$, $Y(0, x_0), Z(0, x_0))$ have the same law on $(\mathcal{C}([0, T], \mathbb{R}^P), \|\cdot\|_\infty) \times (\mathcal{C}([0, T], \mathbb{R}^Q)$, $\|\cdot\|_\infty) \times (L^2([0, T], \mathbb{R}^{Q\times P}), \|\cdot\|_2)$.

Note now that for every function $\varphi$ satisfying the assumptions of point 3 of Theorem (PHP), the function $x \in \mathcal{C}([0, T], \mathbb{R}^P) \mapsto (t \in [0, T] \mapsto \int_0^t \bar{\varphi}(s, x_s, \theta(s, x_s)$, $\nabla_x \theta(s, x_s)) \, ds)$ is continuous for the topology of the uniform convergence (recall that $\nabla_x \theta$ is bounded and locally Hölder continuous). Thanks to Theorem 2.4 and to point 2 of Theorem (PHP), we easily deduce the asymptotic behavior of the processes $(\int_0^\cdot \varphi(s, \bar{X}_s(\varepsilon), X_s(\varepsilon), Y_s(\varepsilon), \hat{Z}_s(\varepsilon)) \, ds)_{\varepsilon > 0}$. This completes the proof. □

## APPENDIX A: PROOFS OF THEOREMS 2.1 AND 2.2

Here is our strategy to prove Theorems 2.1 and 2.2:

1. Concerning point 1 of Theorem 2.1, we refer the reader to Paragraphs 3.2 and 3.3 of Chapter III of Bensoussan, Lions and Papanicolaou (1978) and to Pardoux (1999b) for detailed proofs.
2. We detail the sketch of the proof of Theorem 2.2.
3. We finally prove point 2 of Theorem 2.1.

A.1. *Proof of Theorem 2.2.* In order to solve the auxiliary problems $(\mathcal{A}ux(y))_{y\in\mathbb{R}^Q}$, we study in this section the following kind of Poisson equations on $\mathbb{T}^P$:

$$(A.1) \quad L_y \hat{\varphi}(\cdot, y, u) + \varphi(\cdot, y, u) = 0, \qquad \text{with } \int_{\mathbb{T}^P} \varphi(x, y, u) p(x, y) \, dx = 0,$$

where $u$ is a parameter lying in $\mathbb{R}^d$, $d \in \mathbb{N}^*$. In particular, we aim to study the regularity of the function $\hat{\varphi}$ with respect to the variables $y$ and $u$ [see also Pardoux and Veretennikov (2003) for similar results].

Actually, such equations are also involved in the proof of the ergodic theorem given in Section 2.4. This is the reason why we treat the case where $\varphi$ also depends on a second parameter $u$ (note by the way that $u$, at the opposite of $y$, does not appear in the definition of $L_y$).

Here is the main proposition of this section. [We refer the reader to Paragraphs 3.2 and 3.3 of Chapter III of Bensoussan, Lions and Papanicolaou (1978) for a detailed proof of this result and to Chapters VIII and IX of Gilbarg and Trudinger (1983) for usual estimates of the solutions of elliptic PDEs. Related results are also given in Pardoux and Veretennikov (2002)]:

PROPOSITION A.1. *Let us consider* $q \geq 2$ *and a function* $\varphi : \mathbb{T}^P \times \mathbb{R}^Q \times \mathbb{R}^d \to \mathbb{R}$ *such that*

(A.2)
$$\forall (y, u) \in \mathbb{R}^Q \times \mathbb{R}^d,$$
$$\|\varphi(\cdot, y, u)\|_q < \infty \quad and \quad \int_{\mathbb{T}^P} \varphi(x, y, u) p(x, y) \, dx = 0,$$

*then,* $\forall (y, u) \in \mathbb{R}^Q \times \mathbb{R}^d$, *the Poisson equation on* $\mathbb{T}^P$

(A.3)
$$L_y \hat{\varphi}(\cdot, y, u) + \varphi(\cdot, y, u) = 0, \qquad \int_{\mathbb{T}^P} \hat{\varphi}(x, y, u) \, dx = 0,$$

*admits a unique solution in* $W^{2,q}(\mathbb{T}^P)$. *Moreover, there exists a constant* $C_{A.1}^{(q)}$, *only depending on* $k$, $\lambda$, $\Lambda$, $P$ *and* $q$, *such that*

(A.4)
$$\forall (y, u) \in \mathbb{R}^Q \times \mathbb{R}^d, \qquad \|\hat{\varphi}(\cdot, y, u)\|_{2,q} \leq C_{A.1}^{(q)} \|\varphi(\cdot, y, u)\|_q.$$

*In particular, for* $q > \frac{P}{2}$, *there exists a constant* $C_{A.2}^{(q)}$, *only depending on* $k$, $\lambda$, $\Lambda$, $P$ *and* $q$, *such that:*

(A.5)
$$\forall (y, u) \in \mathbb{R}^Q \times \mathbb{R}^d, \qquad \sup_{x \in \mathbb{T}^P} |\hat{\varphi}(x, y, u)| \leq C_{A.2}^{(q)} \|\varphi(\cdot, y, u)\|_q.$$

*In the following corollaries, we investigate the regularity of* $\hat{\varphi}(\cdot, y, u)$ *with respect to the parameters* $y$ *and* $u$.

COROLLARY A.2. *Let* $\varphi : \mathbb{T}^P \times \mathbb{R}^Q \times \mathbb{R}^d \to \mathbb{R}$ *be such that, for every* $q \geq 2$, *the function* $(y, u) \mapsto \varphi(\cdot, y, u)$ *is continuous from* $\mathbb{R}^Q \times \mathbb{R}^d$ *into* $L^q(\mathbb{T}^P)$ *and satisfies*

(A.6)
$$\forall (y, u) \in \mathbb{R}^Q \times \mathbb{R}^d, \qquad \int_{\mathbb{T}^P} \varphi(x, y, u) p(x, y) \, dx = 0.$$

*Then, for every* $q \geq 2$, *the function* $(y, u) \mapsto \hat{\varphi}(\cdot, y, u)$ *is continuous from* $\mathbb{R}^Q \times \mathbb{R}^d$ *into* $W^{2,q}(\mathbb{T}^P)$, *and, in particular, the function* $\hat{\varphi}$ *is continuous and continuously differentiable with respect to* $x$ *on* $\mathbb{T}^P \times \mathbb{R}^Q \times \mathbb{R}^d$.

Proof.   Let us consider $(y, \delta_1) \in (\mathbb{R}^Q)^2$ and $(u, \delta_2) \in (\mathbb{R}^d)^2$. We just have to note that

$$
\begin{aligned}
\text{(A.7)} \quad & L_y(\hat{\varphi}(\cdot, y + \delta_1, u + \delta_2) - \hat{\varphi}(\cdot, y, u)) + [\varphi(\cdot, y + \delta_1, u + \delta_2) - \varphi(\cdot, y, u)] \\
& \quad + [(L_{y+\delta_1} - L_y)\hat{\varphi}(\cdot, y + \delta_1, u + \delta_2)] = 0.
\end{aligned}
$$

Using Proposition A.1, as well as the assumptions on the coefficients, we complete the proof (note that continuity and $x$-continuous differentiability of $\hat{\varphi}$ on $\mathbb{T}^P \times \mathbb{R}^Q \times \mathbb{R}^d$ follow from the Sobolev imbedding theorems).   $\square$

The following corollary deals with the differentiability of the solution of the Poisson equations with respect to the parameters $y$ and $u$. The proof relies on the estimates given in Proposition A.1 and on the formula (A.7) given in the proof of Corollary A.2. It is left to the reader.

Corollary A.3.   Let $\varphi : \mathbb{T}^P \times \mathbb{R}^Q \times \mathbb{R}^d \to \mathbb{R}$ be such that, for every $q \geq 2$, the function $(y, u) \mapsto \varphi(\cdot, y, u)$ is continuously differentiable from $\mathbb{R}^Q \times \mathbb{R}^d$ into $L^q(\mathbb{T}^P)$ [the partial derivatives are denoted by $(\frac{\partial \varphi}{\partial y_i})_{1 \leq i \leq Q}$ and $(\frac{\partial \varphi}{\partial u_j})_{1 \leq j \leq d}$] and satisfies (A.6). Then, for every $q \geq 2$, the function $(y, u) \mapsto \hat{\varphi}(\cdot, y, u)$ is continuously differentiable from $\mathbb{R}^Q \times \mathbb{R}^d$ into $W^{2,q}(\mathbb{T}^P)$. In particular, from the Sobolev imbedding theorems, the function $\hat{\varphi}$ belongs to $\mathcal{C}^{1,1,1}(\mathbb{T}^P \times \mathbb{R}^Q \times \mathbb{R}^d, \mathbb{R})$ and the function $\nabla_x \hat{\varphi}$ to $\mathcal{C}^{0,1,1}(\mathbb{T}^P \times \mathbb{R}^Q \times \mathbb{R}^d, \mathbb{R}^P)$. Moreover, for all $q \geq 2$, $(y, u) \in \mathbb{R}^Q \times \mathbb{R}^d$ and $(i, j) \in \{1, \ldots, Q\} \times \{1, \ldots, d\}$, $\frac{\partial \hat{\varphi}}{\partial y_i}(\cdot, y, u)$ and $\frac{\partial \hat{\varphi}}{\partial u_j}(\cdot, y, u)$ are the unique solutions in $W^{2,q}(\mathbb{T}^P)$ of the following equations on $\mathbb{T}^P$:

$$
\begin{aligned}
\text{(A.8)} \quad & L_y \frac{\partial \hat{\varphi}}{\partial y_i}(\cdot, y, u) + \frac{\partial \varphi}{\partial y_i}(\cdot, y, u) + \frac{\partial L_y}{\partial y_i} \hat{\varphi}(\cdot, y, u) = 0, \\
& \int_{\mathbb{T}^P} \frac{\partial \hat{\varphi}}{\partial y_i}(x, y, u) p(x, y) \, dx = 0, \\
& L_y \frac{\partial \hat{\varphi}}{\partial u_j}(\cdot, y, u) + \frac{\partial \varphi}{\partial u_j}(\cdot, y, u) = 0, \\
& \int_{\mathbb{T}^P} \frac{\partial \hat{\varphi}}{\partial u_j}(x, y, u) p(x, y) \, dx = 0.
\end{aligned}
$$

Proof of Theorem 2.2.   Let us now turn to the proof of Theorem 2.2. From Proposition A.1, we are able to solve $(\mathcal{A}u x(y))_{y \in \mathbb{R}^Q}$. Therefore, point 1 of Theorem 2.2 is easily established. Moreover, thanks to Assumption $(\mathcal{H})$, we also deduce from Proposition A.1 that, for every $q \geq 2$, there exists a constant $C_{A.3}^{(q)}$ (whose value may change from one inequality to another), only depending on $k$, $K$, $\lambda$, $\Lambda$, $P$, $q$ and $Q$, such that

$$
\text{(A.9)} \qquad \|\hat{b}(\cdot, y)\|_{2,q} + \|\hat{e}(\cdot, y)\|_{2,q} \leq C_{A.3}^{(q)}.
$$

Hence, we just have to study the regularity of $\hat{b}$ and $\hat{e}$ with respect to the variable $y$. Since $b$ and $e$ satisfy the same kind of assumptions, it is sufficient to focus on the case of $b$. In particular, thanks to $(\mathcal{H}.6)$ and $(\mathcal{H}.7)$, we know that, for every $1 \leq \ell \leq P$, the function $(x,y) \in \mathbb{T}^P \times \mathbb{R}^Q \mapsto b_\ell(x,y)$ satisfies the hypotheses of Corollary A.3. Therefore, we deduce that, for every $q \geq 2$, the function $y \in \mathbb{R}^Q \mapsto \hat{b}_\ell(\cdot, y) \in W^{2,q}(\mathbb{T}^P)$ is continuously differentiable. In particular, the function $(x,y) \mapsto \hat{b}_\ell(x,y)$ belongs to $\mathcal{C}^{1,1}(\mathbb{T}^P \times \mathbb{R}^Q, \mathbb{R})$ and the function $(x,y) \mapsto \nabla_x \hat{b}_\ell(x,y)$ to $\mathcal{C}^{0,1}(\mathbb{T}^P \times \mathbb{R}^Q, \mathbb{R}^P)$. Moreover, for all $y \in \mathbb{R}^Q$, the function $\nabla_y \hat{b}(\cdot, y)$, which belongs to $\bigcap_{q \geq 2} W^{2,q}(\mathbb{T}^P)$, satisfies on $\mathbb{T}^P$

(A.10)
$$\forall (\ell, i) \in \{1, \ldots, P\} \times \{1, \ldots, Q\},$$
$$L_y \frac{\partial \hat{b}_\ell}{\partial y_i}(\cdot, y) + \frac{\partial b_\ell}{\partial y_i}(\cdot, y) + \frac{\partial L_y}{\partial y_i} \hat{b}_\ell(\cdot, y) = 0.$$

In particular, we deduce from Assumption $(\mathcal{H})$ and Proposition A.1 that, for every $q \geq 2$,

(A.11)
$$\|\nabla_y \hat{b}(\cdot, y)\|_{2,q} \leq C_{A.3}^{(q)}.$$

Repeating the same argument for the second-order derivatives, we finally establish point 2 of Theorem 2.2. From the Sobolev imbedding theorems, we then deduce point 3 of Theorem 2.2. □

A.2. *Proof of point 2 of Theorem 2.1* In this second section we investigate the regularity of the densities $p(\cdot, y)$ upon the paramater $y$.

LEMMA A.4. *There exists a constant $C_{A.4}$ only depending on $k$, $K$, $\lambda$, $\Lambda$ and $P$, such that*

(A.12)
$$\forall (y, h) \in (\mathbb{R}^Q)^2, \qquad \|p(\cdot, y + h) - p(\cdot, y)\|_2 \leq C_{A.4} |h|.$$

PROOF. Let us consider $\varphi \in L^2(\mathbb{T}^P)$. We first define

(A.13)
$$\forall y \in \mathbb{R}^Q, \qquad \bar{\varphi}(y) = \int_{\mathbb{T}^P} \varphi(x) p(x,y) \, dx.$$

Let us now fix $(y, h) \in (\mathbb{R}^Q)^2$. From the equality

$$\int_{\mathbb{T}^P} p(x,y) \, dx = \int_{\mathbb{T}^P} p(x, y + h) \, dx = 1,$$

we deduce

(A.14)
$$\bar{\varphi}(y + h) - \bar{\varphi}(y) = \int_{\mathbb{T}^P} (\varphi(x) - \bar{\varphi}(y))(p(x, y + h) - p(x, y)) \, dx.$$

Moreover, thanks to Proposition A.1, there exists $\widehat{(\varphi - \bar{\varphi})}(\cdot, y) \in W^{2,2}(\mathbb{T}^P)$ such that

$$(A.15) \qquad L_y \widehat{(\varphi - \bar{\varphi})}(\cdot, y) + (\varphi(\cdot) - \bar{\varphi}(y)) = 0.$$

Hence,

$$(A.16) \qquad \begin{aligned} &\bar{\varphi}(y + h) - \bar{\varphi}(y) \\ &= -\int_{\mathbb{T}^P} L_y \widehat{(\varphi - \bar{\varphi})}(x, y)(p(x, y + h) - p(x, y)) \, dx \\ &= \int_{\mathbb{T}^P} (L_{y+h} - L_y) \widehat{(\varphi - \bar{\varphi})}(x, y) p(x, y + h) \, dx. \end{aligned}$$

Therefore, thanks to Assumption $(\mathcal{H})$, to point 1 of Theorem 2.1 and to Theorem 2.2, there exists a constant $C_{A.4}$ (whose value may change from one inequality to another), only depending on $k$, $K$, $\lambda$, $\Lambda$, $P$, $Q$ and $T$, such that

$$(A.17) \qquad \begin{aligned} |\bar{\varphi}(y + h) - \bar{\varphi}(y)| &\leq C_{A.4} \| \widehat{\varphi - \bar{\varphi}} \|_{2,2} |h| \\ &\leq C_{A.4} \|\varphi - \bar{\varphi}\|_2 |h| \leq C_{A.4} \|\varphi\|_2 |h|. \end{aligned}$$

This completes the proof. $\quad\square$

LEMMA A.5. *The function* $p : y \in \mathbb{R}^Q \mapsto p(\cdot, y) \in L^2(\mathbb{T}^P)$ *is continuously differentiable [the partial derivatives are denoted by* $(\frac{\partial p}{\partial y_j})_{1 \leq j \leq Q}$*]. Moreover, there exist two constants* $C_{A.5}$ *and* $C_{A.6}$*, only depending on* $k$*,* $K$*,* $\lambda$*,* $\Lambda$*,* $P$ *and* $Q$*, such that for all* $j \in \{1, \ldots, Q\}$ *and* $(y, h) \in (\mathbb{R}^Q)^2$*,*

$$(A.18) \qquad \left\| \frac{\partial p}{\partial y_j}(\cdot, y) \right\|_2 \leq C_{A.5},$$

$$(A.19) \qquad \left\| \frac{\partial p}{\partial y_j}(\cdot, y + h) - \frac{\partial p}{\partial y_j}(\cdot, y) \right\|_2 \leq C_{A.6} |h|.$$

PROOF. We consider $\varphi \in L^2(\mathbb{T}^P)$, $j \in \{1, \ldots, Q\}$ and $(y, h) \in \mathbb{R}^Q \times \mathbb{R}^*$, and we set $y + h = y + h e_j$. Thanks to Lemma A.4 [from which we deduce that $y \mapsto p(\cdot, y) \in L^2(\mathbb{T}^P)$ is continuous] and to (A.16), we deduce that $\bar{\varphi}$ is differentiable with respect to $y_j$ and that

$$(A.20) \qquad \frac{\partial \bar{\varphi}}{\partial y_j}(y) = \int_{\mathbb{T}^P} \frac{\partial L_y}{\partial y_j} \widehat{(\varphi - \bar{\varphi})}(x, y) p(x, y) \, dx.$$

Let us now prove that $y \mapsto p(\cdot, y)$ is continuously differentiable. For this purpose, set

$$(A.21) \qquad H_y = \left\{ \psi \in L^2(\mathbb{T}^P), \ \int_{\mathbb{T}^P} \psi(x) p(x, y) \, dx = 0 \right\}.$$

Thanks to point 1 of Theorem 2.1 and to Theorem 2.2, it is readily seen that there exists a constant $C_{A.5}$, only depending on $k$, $K$, $\lambda$, $\Lambda$, $P$ and $Q$, such that

$$(A.22) \qquad \forall \psi \in H_y, \qquad \left| \int_{\mathbb{T}^P} \frac{\partial L_y}{\partial y_j} \hat{\psi}(x,y) p(x,y)\, dx \right| \leq C_{A.5} \|\psi\|_2.$$

Hence, there exists $u_j(\cdot, y) \in H_y$ such that

$$(A.23) \quad \forall \psi \in H_y, \qquad \int_{\mathbb{T}^P} \psi(x) u_j(x,y)\, dx = \int_{\mathbb{T}^P} \frac{\partial L_y}{\partial y_j} \hat{\psi}(x,y) p(x,y)\, dx.$$

Set

$$v_j(\cdot, y) = u_j(\cdot, y) - \left( \int_{\mathbb{T}^P} u_j(x,y)\, dx \right) p(\cdot, y).$$

Hence, for every $\psi \in L^2(\mathbb{T}^P)$,

$$\int_{\mathbb{T}^P} \psi(x) v_j(x,y)\, dx$$

$$(A.24)$$
$$= \int_{\mathbb{T}^P} (\psi(x) - \bar{\psi}(y)) v_j(x,y)\, dx + \bar{\psi}(y) \overbrace{\int_{\mathbb{T}^P} v_j(x,y)\, dx}^{0}$$

$$= \int_{\mathbb{T}^P} (\psi(x) - \bar{\psi}(y)) u_j(x,y)\, dx$$

$$- \int_{\mathbb{T}^P} u_j(x,y)\, dx \overbrace{\int_{\mathbb{T}^P} (\psi(x) - \bar{\psi}(y)) p(x,y)\, dx}^{0}$$

$$= \int_{\mathbb{T}^P} \frac{\partial L_y}{\partial y_j} \widehat{(\psi - \bar{\psi})}(x,y) p(x,y)\, dx.$$

Modifying $C_{A.5}$ if necessary, we deduce that (A.18) is satisfied with $\frac{\partial p}{\partial y_j}$ replaced by $v_j$.

Let us prove that (A.19) holds with $\frac{\partial p}{\partial y_j}$ replaced by $v_j$. To this end, fix once again $\varphi \in L^2(\mathbb{T}^P)$ and $(y, h) \in (\mathbb{R}^Q)^2$.

Following the proof of Corollary A.2, we deduce from Proposition A.1 and Lemma A.4 that there exists a constant $C_{A.6}$ only depending on $k$, $K$, $\lambda$, $\Lambda$, $P$ and $Q$, such that

$$(A.25) \qquad \|\widehat{\varphi - \bar{\varphi}}(\cdot, y + h) - \widehat{\varphi - \bar{\varphi}}(\cdot, y)\|_{2,2} \leq C_{A.6} |h| \|\varphi\|_2.$$

Hence, thanks to Lemma A.4, to (A.24) and to Assumption $(\mathcal{H})$, we deduce (modifying $C_{A.6}$ if necessary)

$$\left| \int_{\mathbb{T}^P} \varphi(x) (v_j(x, y+h) - v_j(x,y))\, dx \right|$$

$$
\text{(A.26)} \quad
\begin{aligned}
&= \left| \int_{\mathbb{T}^P} \frac{\partial L_{y+h}}{\partial y_j}(\widehat{\varphi - \bar{\varphi}})(x, y+h) p(x, y+h)\, dx \right.\\
&\qquad \left. - \int_{\mathbb{T}^P} \frac{\partial L_y}{\partial y_j}(\widehat{\varphi - \bar{\varphi}})(x, y) p(x, y)\, dx \right|\\
&\leq C_{A.6} \|\varphi\|_2 |h|.
\end{aligned}
$$

Therefore, (A.19) holds with $\frac{\partial p}{\partial y_j}$ replaced by $v_j$. Moreover, thanks to (A.20), (A.24) and (A.26), we deduce from the mean value theorem that

$$
\text{(A.27)} \quad
\begin{aligned}
&\left| \frac{1}{h} \int_{\mathbb{T}^P} \varphi(x)(p(x, y+h) - p(x, y))\, dx \right.\\
&\qquad \left. - \int_{\mathbb{T}^P} \varphi(x) v_j(x, y)\, dx \right| \leq C_{A.6} \|\varphi\|_2 |h|.
\end{aligned}
$$

Therefore, $p: y \in \mathbb{R}^P \mapsto p(\cdot, y) \in L^2(\mathbb{T}^P)$ is continuously differentiable, and, for every $j \in \{1, \ldots, Q\}$, $\frac{\partial p}{\partial y_j} = v_j$. $\square$

PROOF OF THEOREM 2.1. We already know that point 1 of Theorem 2.1 holds. According to Lemmas A.4 and A.5, we know that the function $y \in \mathbb{R}^Q \mapsto p(\cdot, y) \in L^2(\mathbb{T}^P)$ is continuously differentiable with respect to $y$ and that there exists a constant $C_{A.5}$, given by Lemma A.5, such that

$$
\text{(A.28)} \quad \forall i \in \{1, \ldots, Q\}, \ \forall y \in \mathbb{R}^Q, \qquad \left\| \frac{\partial p}{\partial y_i}(\cdot, y) \right\|_2 \leq C_{A.5}.
$$

Actually, repeating the arguments of Lemma A.5 [and using in particular (A.18) and (A.19)], we can study in the same way the differentiability of order 2 of $p$, and then achieve the proof of Theorem 2.2. $\square$

## APPENDIX B: PROOF OF THEOREM 2.3

This section is devoted to the proof of Theorem 2.3. Following Pardoux (1999a, b), we start with the following:

LEMMA B.1. *There exists a nondecreasing function $C_{B.1}: \mathbb{R}_+ \to \mathbb{R}_+$, only depending on $\lambda$ and $\Lambda$, such that*

$$
\forall \varepsilon > 0, \ \forall t \in [0, T],
$$

$$
\text{(B.1)} \quad
\begin{aligned}
&\mathbf{E} \int_t^T (1 + |Y_s(\varepsilon)|) |Z_s(\varepsilon)|^2\, ds \\
&\qquad \leq C_{B.1}(T) \left[ 1 + \varepsilon^{-1} \mathbf{E} \int_t^T (1 + |Y_s(\varepsilon)|^2)\, ds \right].
\end{aligned}
$$

PROOF.    Consider the function $\Psi : y \in \mathbb{R}^Q \to [1 + |y|^2]^{3/2}$. From Itô's formula, $\forall \varepsilon > 0$,

$$
\begin{aligned}
(\text{B.2}) \quad d\Psi(Y_t(\varepsilon)) &= -3[1 + |Y_t(\varepsilon)|^2]^{1/2} \langle Y_t(\varepsilon), (\varepsilon^{-1}e + f)(\bar{X}_t(\varepsilon), Y_t(\varepsilon), Z_t(\varepsilon)) \rangle \, dt \\
&\quad + \tfrac{3}{2}[1 + |Y_t(\varepsilon)|^2]^{1/2} |Z_t(\varepsilon)\sigma(\bar{X}_t(\varepsilon), Y_t(\varepsilon))|^2 \, dt \\
&\quad + 3[1 + |Y_t(\varepsilon)|^2]^{1/2} \langle Y_t(\varepsilon), Z_t(\varepsilon)\sigma(\bar{X}_t(\varepsilon), Y_t(\varepsilon)) \, dB_t \rangle \\
&\quad + \tfrac{3}{2}[1 + |Y_t(\varepsilon)|^2]^{-1/2} |Y_t^*(\varepsilon)Z_t(\varepsilon)\sigma(\bar{X}_t(\varepsilon), Y_t(\varepsilon))|^2 \, dt.
\end{aligned}
$$

Recall from Section 2.1 that, for every $\varepsilon > 0$, $\theta_\varepsilon$ is bounded. In particular, for all $\varepsilon > 0$, $\exists \Gamma(\varepsilon) > 0$ s.t. $\mathbf{P}\{\sup_{0 \leq t \leq T} |Y_t(\varepsilon)| \leq \Gamma(\varepsilon)\} = 1$. Hence, $\forall \varepsilon > 0$, $\forall t \in [0, T]$,

$$
\begin{aligned}
(\text{B.3}) \quad & \mathbf{E}\Psi(Y_t(\varepsilon)) + \tfrac{3}{2}\mathbf{E}\int_t^T [1 + |Y_s(\varepsilon)|^2]^{1/2} |Z_s(\varepsilon)\sigma(\bar{X}_s(\varepsilon), Y_s(\varepsilon))|^2 \, ds \\
&+ \tfrac{3}{2}\mathbf{E}\int_t^T [1 + |Y_s(\varepsilon)|^2]^{-1/2} |Y_s^*(\varepsilon)Z_s(\varepsilon)\sigma(\bar{X}_s(\varepsilon), Y_s(\varepsilon))|^2 \, ds \\
&= \mathbf{E}\Psi(Y_T(\varepsilon)) + 3\mathbf{E}\int_t^T [1 + |Y_s(\varepsilon)|^2]^{1/2} \\
&\qquad\qquad\qquad \times \langle Y_s(\varepsilon), (\varepsilon^{-1}e + f)(\bar{X}_s(\varepsilon), Y_s(\varepsilon), Z_s(\varepsilon)) \rangle \, ds.
\end{aligned}
$$

Therefore, from the boundedness of $H$ [and then of $Y_T(\varepsilon)$] and the properties of $e$, $f$ and $\sigma$, there exists a constant $c_{B.1}$, only depending on $\lambda$ and $\Lambda$, such that

$$
\begin{aligned}
(\text{B.4}) \quad & \mathbf{E}\Psi(Y_t(\varepsilon)) + \mathbf{E}\int_t^T [1 + |Y_s(\varepsilon)|^2]^{1/2} |Z_s(\varepsilon)|^2 \, ds \\
&+ \mathbf{E}\int_t^T [1 + |Y_s(\varepsilon)|^2]^{-1/2} |Y_s^*(\varepsilon)Z_s(\varepsilon)|^2 \, ds \\
&\leq c_{B.1} + \varepsilon^{-1} c_{B.1} \mathbf{E}\int_t^T [1 + |Y_s(\varepsilon)|^2] \, ds + c_{B.1} \mathbf{E}\int_t^T \Psi(Y_s(\varepsilon)) \, ds.
\end{aligned}
$$

In particular, from Gronwall's lemma, there exists a nondecreasing function $C_{B.1} : \mathbb{R}_+ \to \mathbb{R}_+$, only depending on $\lambda$ and $\Lambda$, such that

$$
\begin{aligned}
(\text{B.5}) \quad & \forall \varepsilon > 0, \ \forall t \in [0, T], \\
& \mathbf{E}\Psi(Y_t(\varepsilon)) \leq C_{B.1}(T)\left[1 + \varepsilon^{-1}\mathbf{E}\int_t^T [1 + |Y_s(\varepsilon)|^2] \, ds\right].
\end{aligned}
$$

Therefore, modifying $C_{B.1}$ if necessary, we deduce from (B.4) and (B.5) that for all $\varepsilon > 0$ and $t \in [0, T]$,

$$
\mathbf{E}\int_t^T [1 + |Y_s(\varepsilon)|^2]^{1/2} |Z_s(\varepsilon)|^2 \, ds
$$

(B.6)
$$\leq C_{B.1}(T)\left[1 + \varepsilon^{-1}\mathbf{E}\int_t^T [1 + |Y_s(\varepsilon)|^2]\,ds\right],$$

from which the desired result follows easily. $\square$

PROPOSITION B.2. *There exists a nondecreasing function* $C_{B.2}\colon\mathbb{R}_+ \to \mathbb{R}_+$, *only depending on* $k$, $K$, $\lambda$, $\Lambda$, $P$ *and* $Q$, *such that*

(B.7) $\qquad \forall \varepsilon > 0, \qquad \sup_{0 \leq t \leq T} \mathbf{E}|Y_t(\varepsilon)|^2 + \mathbf{E}\int_0^T |Z_s(\varepsilon)|^2\,ds \leq C_{B.2}(T).$

PROOF. From (2.5), we deduce

$$
\begin{aligned}
d|\hat{Y}_t(\varepsilon)|^2 = {} & -2\langle \hat{Y}_t(\varepsilon), v(\bar{X}_t(\varepsilon), Y_t(\varepsilon), \hat{Z}_t(\varepsilon))\rangle\,dt \\
& + |\hat{Z}_t(\varepsilon)\sigma(\bar{X}_t(\varepsilon), Y_t(\varepsilon)) - \varepsilon\nabla_y \hat{e}(\bar{X}_t(\varepsilon), Y_t(\varepsilon)) \\
& \qquad\qquad \times Z_t(\varepsilon)\sigma(\bar{X}_t(\varepsilon), Y_t(\varepsilon))|^2\,dt
\end{aligned}
$$

(B.8)
$$
\begin{aligned}
& - \varepsilon\langle \hat{Y}_t(\varepsilon), \nabla^2_{y,y}\hat{e}(\bar{X}_t(\varepsilon), Y_t(\varepsilon))[Z_t(\varepsilon)a(\bar{X}_t(\varepsilon), Y_t(\varepsilon))Z_t^*(\varepsilon)]\rangle\,dt \\
& + 2\varepsilon\langle \hat{Y}_t(\varepsilon), (\nabla_y\hat{e}f)(\bar{X}_t(\varepsilon), Y_t(\varepsilon), Z_t(\varepsilon))\rangle\,dt \\
& + 2\langle \hat{Y}_t(\varepsilon), \hat{Z}_t(\varepsilon)\sigma(\bar{X}_t(\varepsilon), Y_t(\varepsilon))\,dB_t\rangle \\
& - 2\varepsilon\langle \hat{Y}_t(\varepsilon), \nabla_y\hat{e}(\bar{X}_t(\varepsilon), Y_t(\varepsilon))Z_t(\varepsilon)\sigma(\bar{X}_t(\varepsilon), Y_t(\varepsilon))\,dB_t\rangle.
\end{aligned}
$$

We obtain

(B.9)
$$
\begin{aligned}
\mathbf{E}|\hat{Y}_t(\varepsilon)|^2 + \mathbf{E}\int_t^T & |\hat{Z}_s(\varepsilon)\sigma(\bar{X}_s(\varepsilon), Y_s(\varepsilon)) \\
& - \varepsilon\nabla_y\hat{e}(\bar{X}_s(\varepsilon), Y_s(\varepsilon))Z_s(\varepsilon)\sigma(\bar{X}_s(\varepsilon), Y_s(\varepsilon))|^2\,ds \\
= \mathbf{E}|\hat{Y}_T(\varepsilon)|^2 + 2\mathbf{E}\int_t^T & \langle \hat{Y}_s(\varepsilon), v(\bar{X}_s(\varepsilon), Y_s(\varepsilon), \hat{Z}_s(\varepsilon))\rangle\,ds \\
+ \varepsilon\mathbf{E}\int_t^T & \langle \hat{Y}_s(\varepsilon), \nabla^2_{y,y}\hat{e}(\bar{X}_s(\varepsilon), Y_s(\varepsilon)) \\
& \qquad \times [Z_s(\varepsilon)a(\bar{X}_s(\varepsilon), Y_s(\varepsilon))Z_s^*(\varepsilon)]\rangle\,ds \\
- 2\varepsilon\mathbf{E}\int_t^T & \langle \hat{Y}_s(\varepsilon), (\nabla_y\hat{e}f)(\bar{X}_s(\varepsilon), Y_s(\varepsilon), Z_s(\varepsilon))\rangle\,ds.
\end{aligned}
$$

Recall that $\hat{Z}.(\varepsilon) = Z.(\varepsilon) - \nabla_x\hat{e}(\bar{X}.(\varepsilon), Y.(\varepsilon))$. Therefore, thanks to Assumption $(\mathcal{H})$, to (2.5.i) and to Theorem 2.2, there exists a nondecreasing function $c_{B.2}\colon\mathbb{R}_+ \to \mathbb{R}_+$, only depending on $k$, $K$, $\lambda$, $\Lambda$, $P$ and $Q$, such that

$\forall\, 0 < \varepsilon < 1,\ \forall\, t \in [0, T],$

$$\mathbf{E}|Y_t(\varepsilon)|^2 + \mathbf{E}\int_t^T |Z_s(\varepsilon)|^2\, ds$$

(B.10)
$$\leq c_{B.2}(T)\bigg[1 + \mathbf{E}\int_t^T |Y_s(\varepsilon)|(1 + |Y_s(\varepsilon)| + |Z_s(\varepsilon)|)\, ds$$

$$+ \varepsilon \mathbf{E}\int_t^T |Y_s(\varepsilon)||Z_s(\varepsilon)|^2\, ds + \varepsilon^2 \mathbf{E}\int_t^T |Z_s(\varepsilon)|^2\, ds\bigg].$$

Therefore, from Lemma B.1, there exist a nondecreasing function $C_{B.2}\colon \mathbb{R}_+ \to \mathbb{R}_+$, only depending on $k$, $K$, $\lambda$, $\Lambda$, $P$ and $Q$, and a nonincreasing function $\eta_{B.1}\colon \mathbb{R}_+ \to \mathbb{R}_+$, only depending on $k$, $K$, $\lambda$, $\Lambda$, $P$ and $Q$, such that $\forall\, 0 < \varepsilon < \eta_{B.1}(T),\ \forall\, t \in [0, T],$

(B.11) $\quad \mathbf{E}|Y_t(\varepsilon)|^2 + \mathbf{E}\int_t^T |Z_s(\varepsilon)|^2\, ds \leq C_{B.2}(T)\bigg[1 + \mathbf{E}\int_t^T |Y_s(\varepsilon)|^2\, ds\bigg].$

Hence, applying Gronwall's lemma, and modifying $C_{B.2}$ if necessary, we have $\forall\, 0 < \varepsilon < \eta_{B.1}(T),$

(B.12) $$\sup_{0 \leq t \leq T} \mathbf{E}|Y_t(\varepsilon)|^2 + \mathbf{E}\int_0^T |Z_s(\varepsilon)|^2\, ds \leq C_{B.2}(T).$$

For every $\varepsilon \geq \eta_{B.1}(T)$, we can apply usual estimates of solutions of FBSDEs to the system $(\mathrm{E}(\varepsilon, 0, x_0))$ [see, e.g., Theorem 3.2 in Delarue (2002b)]. This completes the proof of Proposition B.2. $\quad\square$

PROPOSITION B.3. *There exists a constant $C_{B.3}$, only depending on $k$, $K$, $\lambda$, $\Lambda$, $P$, $Q$ and $T$, such that*

(B.13) $$\forall\, \varepsilon > 0, \qquad \mathbf{P}\bigg\{\sup_{t \in [0, T]} |Y_t(\varepsilon)| \leq C_{B.3}\bigg\} = 1.$$

PROOF. Recall that for all $\varepsilon > 0$, $t \in [0, T]$ and $x \in \mathbb{R}^P$, $(X_s(\varepsilon, t, x), Y_s(\varepsilon, t, x), Z_s(\varepsilon, t, x))_{t \leq s \leq T}$ denotes the solution of the FBSDE $\mathrm{E}(\varepsilon, t, x)$. From Proposition B.2, we then deduce

$\forall\, x \in \mathbb{R}^P,\ \forall\, \varepsilon > 0,\ \forall\, t \in [0, T],$

(B.14)
$$\sup_{t \leq s \leq T} \mathbf{E}|Y_s(\varepsilon, t, x)|^2 \leq C_{B.2}(T - t) \leq C_{B.2}(T).$$

Hence, from Theorems 3.1 and 3.2 of Delarue (2002b), we complete the proof of Proposition B.3. $\quad\square$

COROLLARY B.4. *There exists a constant $C_{B.4}$, only depending on $k$, $K$, $\lambda$, $\Lambda$, $P$, $Q$ and $T$, such that*

$$\text{(B.15)} \qquad \forall \varepsilon > 0, \qquad \mathbf{E}\left(\int_0^T |Z_s(\varepsilon)|^2 \, ds\right)^2 \leq C_{B.4}.$$

PROOF. We know from (B.8) that $\forall \varepsilon > 0$, $\forall t \in [0, T]$,

$$
\begin{aligned}
|\hat{Y}_t(\varepsilon)|^2 &+ \int_t^T |\hat{Z}_s(\varepsilon)\sigma(\bar{X}_s(\varepsilon), Y_s(\varepsilon)) \\
&\qquad - \varepsilon\nabla_y\hat{e}(\bar{X}_s(\varepsilon), Y_s(\varepsilon))Z_s(\varepsilon)\sigma(\bar{X}_s(\varepsilon), Y_s(\varepsilon))|^2 \, ds \\
&= |\hat{Y}_T(\varepsilon)|^2 + 2\int_t^T \langle \hat{Y}_s(\varepsilon), v(\bar{X}_s(\varepsilon), Y_s(\varepsilon), \hat{Z}_s(\varepsilon)) \rangle \, ds \\
&\quad + \varepsilon\int_t^T \langle \hat{Y}_s(\varepsilon), \nabla^2_{y,y}\hat{e}(\bar{X}_s(\varepsilon), Y_s(\varepsilon)) \\
&\qquad\qquad\qquad \times [Z_s(\varepsilon)a(\bar{X}_s(\varepsilon), Y_s(\varepsilon))Z_s^*(\varepsilon)] \rangle \, ds \\
&\quad - 2\varepsilon\int_t^T \langle \hat{Y}_s(\varepsilon), (\nabla_y\hat{e}f)(\bar{X}_s(\varepsilon), Y_s(\varepsilon), Z_s(\varepsilon)) \rangle \, ds \\
&\quad - 2\int_t^T \langle \hat{Y}_s(\varepsilon), \hat{Z}_s(\varepsilon)\sigma(\bar{X}_s(\varepsilon), Y_s(\varepsilon)) \, dB_s \rangle \\
&\quad + 2\varepsilon\int_t^T \langle \hat{Y}_s(\varepsilon), \nabla_y\hat{e}(\bar{X}_s(\varepsilon), Y_s(\varepsilon))Z_s(\varepsilon)\sigma(\bar{X}_s(\varepsilon), Y_s(\varepsilon)) \, dB_s \rangle.
\end{aligned}
\tag{B.16}
$$

Hence, from Proposition B.3 and (2.5.i), there exists a constant $c_{B.3}$ (whose value may change from one inequality to another), only depending on $k$, $K$, $\lambda$, $\Lambda$, $P$, $Q$ and $T$, such that $\forall\, 0 < \varepsilon < 1$,

$$
\begin{aligned}
\int_0^T &|Z_s(\varepsilon)|^2 \, ds \\
&\leq c_{B.3} + c_{B.3}\int_0^T |Z_s(\varepsilon)| \, ds + c_{B.3}\varepsilon\int_0^T |Z_s(\varepsilon)|^2 \, ds \\
&\quad - 2\int_0^T \langle \hat{Y}_s(\varepsilon), \hat{Z}_s(\varepsilon)\sigma(\bar{X}_s(\varepsilon), Y_s(\varepsilon)) \, dB_s \rangle \\
&\quad + 2\varepsilon\int_0^T \langle \hat{Y}_s(\varepsilon), \nabla_y\hat{e}(\bar{X}_s(\varepsilon), Y_s(\varepsilon))Z_s(\varepsilon)\sigma(\bar{X}_s(\varepsilon), Y_s(\varepsilon)) \, dB_s \rangle.
\end{aligned}
\tag{B.17}
$$

Therefore, thanks to Theorem 2.2 (which provides an estimate of $\nabla_x\hat{e}$), we have $\forall\, 0 < \varepsilon < 1$,

$$\mathbf{E}\left(\int_0^T |Z_s(\varepsilon)|^2 \, ds\right)^2$$

$$
\text{(B.18)} \qquad \leq c_{B.3}\Big[1 + \mathbf{E}\Big(\int_0^T |Z_s(\varepsilon)|\,ds\Big)^2
$$

$$
+ \varepsilon^2 \mathbf{E}\Big(\int_0^T |Z_s(\varepsilon)|^2\,ds\Big)^2 + \mathbf{E}\int_0^T |Z_s(\varepsilon)|^2\,ds\Big].
$$

In particular, from Proposition B.2, we have $\forall\, 0 < \varepsilon < 1$,

$$
\text{(B.19)} \qquad \mathbf{E}\Big(\int_0^T |Z_s(\varepsilon)|^2\,ds\Big)^2 \leq c_{B.3} + c_{B.3}\varepsilon^2 \mathbf{E}\Big(\int_0^T |Z_s(\varepsilon)|^2\,ds\Big)^2.
$$

Moreover [see, e.g., Delarue (2002b)],

$$
\forall\, \varepsilon > 0, \qquad \mathbf{E}\Big(\int_0^T |Z_s(\varepsilon)|^2\,ds\Big)^2 < \infty.
$$

We deduce that there exist two constants $C_{B.4}$ and $\eta_{B.2}$, only depending on $k$, $K$, $\lambda$, $\Lambda$, $P$, $Q$ and $T$, such that

$$
\text{(B.20)} \qquad \sup_{0 < \varepsilon < \eta_{B.2}} \mathbf{E}\Big(\int_0^T |Z_s(\varepsilon)|^2\,ds\Big)^2 \leq C_{B.4}.
$$

For every $\varepsilon \geq \eta_{B.2}$, we can apply usual estimates of solutions of FBSDEs to the system $(\mathrm{E}(\varepsilon, 0, x_0))$ [see once again Delarue (2002b)]. This completes the proof. $\qquad \square$

From Propositions B.2 and B.3, Corollary B.4 and (2.4), we deduce the following:

COROLLARY B.5.   *There exists a constant $C_{B.5}$, only depending on $k$, $K$, $\lambda$, $\Lambda$, $P$, $Q$ and $T$, such that*

$$
\text{(B.21)} \qquad \forall\, \varepsilon > 0, \qquad \mathbf{E}\sup_{t \in [0,T]} |X_t(\varepsilon)|^2 \leq C_{B.5}.
$$

*Note that Theorem* 2.3 *easily follows from Proposition* B.3 *and Corollaries* B.4 *and* B.5.

## APPENDIX C: PROOF OF THEOREM 2.4

Let us now turn to the proof of Theorem 2.4. To this end, keep the notation introduced in Section 2.4 and set

$$
\text{(C.1)} \quad c_{C.1} = \sup_{0 < \varepsilon < 1}\Big(\mathbf{E}|G_0(\varepsilon)| + \mathbf{E}\int_0^T (|g_t(1,\varepsilon)| + |g_t(2,\varepsilon)|^2)\,dt\Big) < \infty.
$$

The strategy is quite clear: if the function $\varphi$ is regular with respect to $t$, $y$ and $u$ and null outside a compact set, the solution, denoted by $\hat{\varphi}$, of the family of equations on $\mathbb{T}^P$

$$\forall\,(t,y,g)\in[0,T]\times\mathbb{R}^Q\times\mathbb{R}^d,$$

(C.2)
$$L_y\hat\varphi(t,\cdot,y,g)+\varphi(t,\cdot,y,g)-\bar\varphi(t,y,g)=0,$$

$$\int_{\mathbb{T}^P}\hat\varphi(t,x,y,g)\,dx=0,$$

is also regular with respect to $t$, $y$ and $u$ and null outside a compact set. Hence, thanks to Itô's formula, we have for all $\varepsilon>0$ and $t\in[0,T]$,

(C.3)
$$\begin{aligned}
&\varepsilon^2\hat\varphi(t,\bar X_t(\varepsilon),Y_t(\varepsilon),G_t(\varepsilon))\\
&\quad=\varepsilon^2\hat\varphi\Big(0,\frac{x_0}{\varepsilon},Y_0(\varepsilon),G_0(\varepsilon)\Big)\\
&\quad\quad+\int_0^t[\bar\varphi-\varphi](s,\bar X_s(\varepsilon),Y_s(\varepsilon),G_s(\varepsilon))\,ds+\varepsilon r_t(\varepsilon),
\end{aligned}$$

where $(r(\varepsilon))_{\varepsilon>0}$ is a family of semi-martingales satisfying $\mathbf{E}\sup_{t\in[0,T]}|r_t(\varepsilon)|<\infty$. The result follows then easily.

Since $\varphi$ is not regular, we apply this strategy to a regularization sequence of $\varphi$. To this end, we first extend $\varphi$ to the whole set $\mathbb{R}\times\mathbb{R}^P\times\mathbb{R}^Q\times\mathbb{R}^d$ by setting for every $(t,x,y,g)$, $\varphi(t,x,y,g)=\varphi(0,x,y,g)$ if $t\le 0$ and $\varphi(t,x,y,g)=\varphi(T,x,y,g)$ if $t\ge T$. Following the notation given in the Introduction, we define for all $n\in\mathbb{N}^*$ and $(t,x,y,g)\in\mathbb{R}\times\mathbb{R}^P\times\mathbb{R}^Q\times\mathbb{R}^d$,

(C.4)
$$\begin{aligned}
\varphi_n(t,x,y,g)=n^{Q+d+1}\int\varphi(t',x,y',g')\rho_{Q+d+1}\\
\times(n(t-t',y-y',g-g'))\,dt'\,dy'\,dg'.
\end{aligned}$$

Note from the assumptions made on $\varphi$ that $(\varphi_n)_{n\in\mathbb{N}^*}$ uniformly converges on every compact subset of $\mathbb{R}\times\mathbb{R}^P\times\mathbb{R}^Q\times\mathbb{R}^d$ toward $\varphi$. We set

(C.5)
$$\forall\,n\in\mathbb{N}^*,\ \forall\,(t,y,g)\in\mathbb{R}\times\mathbb{R}^Q\times\mathbb{R}^d,$$

$$\bar\varphi_n(t,y,g)=\int_{\mathbb{T}^P}\varphi_n(t,x,y,g)p(x,y)\,dx.$$

Thanks to Theorem 2.1, $\bar\varphi_n$ is twice continuously differentiable with respect to $t$, $y$ and $g$. Moreover, $(\bar\varphi_n)_{n\in\mathbb{N}^*}$ uniformly converges on every compact set toward $\bar\varphi$ as $n\to+\infty$.

Thanks to Proposition A.1, we can find a sequence of functions $(\hat\varphi_n)_{n\in\mathbb{N}^*}$, such that for every $n\in\mathbb{N}^*$,

$$\forall\,(t,y,g)\in\mathbb{R}\times\mathbb{R}^Q\times\mathbb{R}^d,$$

(C.6)
$$L_y\hat\varphi_n(t,\cdot,y,g)+[\varphi_n-\bar\varphi_n](t,\cdot,y,g)=0,$$

$$\int_{\mathbb{T}^P}\hat\varphi_n(t,x,y,g)\,dx=0.$$

Hence, thanks to Proposition A.1, for every $q \geq 2$, there exists a constant $C_{C.1}^{(q)}$, only depending on $\|\varphi\|_\infty$, $k$, $\lambda$, $\Lambda$, $P$ and $q$, such that

$$(C.7) \quad \forall n \in \mathbb{N}^*, \ \forall (t, y, g) \in \mathbb{R} \times \mathbb{R}^Q \times \mathbb{R}^d, \qquad \|\hat{\varphi}_n(t, \cdot, y, g)\|_{2,q} \leq C_{C.1}^{(q)}.$$

Thanks to Corollary A.3, we know that for all $n \in \mathbb{N}^*$ and $q \geq 2$, the function $(t, y, g) \mapsto \hat{\varphi}_n(t, \cdot, y, g)$ is twice continuously differentiable from $\mathbb{R} \times \mathbb{R}^Q \times \mathbb{R}^d$ into $W^{2,q}(\mathbb{T}^P)$ (in particular, $\hat{\varphi}_n$ and $\nabla_x \hat{\varphi}_n$ are twice continuously differentiable with respect to $t$, $y$ and $g$). Moreover, for all $q \geq 2$ and $n \in \mathbb{N}^*$,

$$
\begin{aligned}
\sup_{(t,y,g) \in \mathbb{R} \times \mathbb{R}^Q \times \mathbb{R}^d} & [\|\nabla_t \hat{\varphi}_n(t, \cdot, y, g)\|_{2,q} + \|\nabla_y \hat{\varphi}_n(t, \cdot, y, g)\|_{2,q} \\
(C.8) \qquad\qquad & + \|\nabla_g \hat{\varphi}_n(t, \cdot, y, g)\|_{2,q} + \|\nabla_{y,g}^2 \hat{\varphi}_n(t, \cdot, y, g)\|_{2,q} \\
& + \|\nabla_{y,y}^2 \hat{\varphi}_n(t, \cdot, y, g)\|_{2,q} + \|\nabla_{g,g}^2 \hat{\varphi}_n(t, \cdot, y, g)\|_{2,q}] < \infty.
\end{aligned}
$$

Hence, from the Itô–Krylov formula [we refer once again to Krylov (1980) and to Pardoux and Veretennikov (2002)] we have for all $n \in \mathbb{N}^*$, $\varepsilon > 0$ and $t \in [0, T]$,

$$
\begin{aligned}
d\hat{\varphi}_n & (t, \bar{X}_t(\varepsilon), Y_t(\varepsilon), G_t(\varepsilon)) \\
(C.9) \qquad & = \frac{1}{\varepsilon^2} [L_{Y_t(\varepsilon)} \hat{\varphi}_n](t, \bar{X}_t(\varepsilon), Y_t(\varepsilon), G_t(\varepsilon)) \, dt + \frac{1}{\varepsilon} \, dr_t(n, \varepsilon) \\
& = \frac{1}{\varepsilon^2} [\bar{\varphi}_n - \varphi_n](t, \bar{X}_t(\varepsilon), Y_t(\varepsilon), G_t(\varepsilon)) \, dt + \frac{1}{\varepsilon} \, dr_t(n, \varepsilon),
\end{aligned}
$$

where $(r(n, \varepsilon))_{n \in \mathbb{N}^*, \varepsilon > 0}$ is, thanks to Theorem 2.3, to (C.1) and to (C.8), a family of semimartingales satisfying

$$(C.10) \qquad \forall n \in \mathbb{N}^*, \qquad c_{C.1}(n) = \sup_{0 < \varepsilon < 1} \mathbf{E} \sup_{0 \leq t \leq T} |r_t(n, \varepsilon)| < \infty.$$

Hence, from (C.7) and (C.10), there exists an $\mathbb{R}_+$ valued sequence, denoted by $(\gamma_n)_{n \in \mathbb{N}^*}$, such that

$$(C.11) \ \forall 0 < \varepsilon < 1, \qquad \mathbf{E} \sup_{0 \leq t \leq T} \left| \int_0^t [\varphi_n - \bar{\varphi}_n](s, \bar{X}_s(\varepsilon), Y_s(\varepsilon), G_s(\varepsilon)) \, ds \right| \leq \varepsilon \gamma_n.$$

Therefore, from Theorem 2.3, there exists a constant $C_{C.2}$, only depending on $\|\varphi\|_\infty$ and $T$, such that for all $(n, N) \in (\mathbb{N}^*)^2$ and $0 < \varepsilon < 1$,

$$
\begin{aligned}
\mathbf{E} \sup_{0 \leq t \leq T} & \left| \int_0^t [\varphi - \bar{\varphi}](s, \bar{X}_s(\varepsilon), Y_s(\varepsilon), G_s(\varepsilon)) \, ds \right| \\
(C.12) \qquad & \leq \varepsilon \gamma_n + C_{C.2} \Bigg( \mathbf{P} \Big\{ \sup_{0 \leq t \leq T} |G_t(\varepsilon)| \geq N \Big\} \\
& \qquad\qquad + \sup_{\Delta_N} \{ (|\varphi - \varphi_n| + |\bar{\varphi} - \bar{\varphi}_n|)(t, x, y, g) \} \Bigg),
\end{aligned}
$$

where $\Delta_N = [0, T] \times \mathbb{T}^P \times \{y \in \mathbb{R}^Q, \; |y| \leq C_{2.4}\} \times \{g \in \mathbb{R}^d, \; |g| \leq N\}$.

From (C.1), there exists a constant $C_{C.3}$, only depending on $c_{C.1}$, $\|\varphi\|_\infty$ and $T$, such that for all $(n, N) \in (\mathbb{N}^*)^2$ and $0 < \varepsilon < 1$,

(C.13)
$$\mathbf{E} \sup_{0 \leq t \leq T} \left| \int_0^t [\varphi - \bar{\varphi}](s, \bar{X}_s(\varepsilon), Y_s(\varepsilon), G_s(\varepsilon)) \, ds \right|$$
$$\leq \varepsilon \gamma_n + C_{C.3} \left( \frac{1}{N} + \sup_{\Delta_N} \{ (|\varphi - \varphi_n| + |\bar{\varphi} - \bar{\varphi}_n|)(t, x, y, g) \} \right).$$

Point 1 of Theorem 2.4 easily follows. Point 2 is a direct consequence of point 1.

## APPENDIX D: PROOF OF SOLVABILITY PROPERTIES

We now establish the solvability results mentioned in point 1 of Theorem (HP), in Section 2.5 and in Section 3.3. In particular, we aim to do the following:

1. To prove that the systems of PDEs $(\mathcal{E}(\varepsilon))_{\varepsilon > 0}$, $\mathcal{E}(\lim)$ and $(\mathcal{E}_{\text{reg}}(n))_{n \in \mathbb{N}^*}$ [recall that systems $(\mathcal{E}_{\text{reg}}(n))_{n \in \mathbb{N}^*}$ are given in Section 3.3] are uniquely solvable in the space $\mathcal{V}$ and to establish that the solutions of $(\mathcal{E}(\varepsilon))_{\varepsilon > 0}$ and $\mathcal{E}(\lim)$ coincide with $(\theta_\varepsilon)_{\varepsilon > 0}$ and $\theta$ [refer to Sections 2.1 and 2.5 for the definitions of $(\theta_\varepsilon)_{\varepsilon > 0}$ and $\theta$].
2. To estimate the solutions of these systems and, in particular, to establish (2.13) and (3.8).

**Solvability of $(\mathcal{E}(\varepsilon))_{\varepsilon > 0}$, $\mathcal{E}(\lim)$ and $(\mathcal{E}_{\text{reg}}(n))_{n \in \mathbb{N}^*}$.** Since other cases are similar, we just expose the proof of the unique solvability of $\mathcal{E}(1)$. To this end, set for the sake of simplicity,

(D.1)
$$\forall (x, y, z) \in \mathbb{R}^P \times \mathbb{R}^Q \times \mathbb{R}^{Q \times P}, \qquad C(x, y, z) = (b + c)(x, y, z),$$
$$F(x, y, z) = (e + f)(x, y, z).$$

We first prove that $\mathcal{E}(1)$ admits a solution in the space $\mathcal{V}$. For this purpose, we consider a sequence of functions, denoted by $(C_m, F_m, \sigma_m, h_m)_{m \in \mathbb{N}}$, satisfying [see, e.g., Delarue (2002b) for the construction of such a sequence].

($\mathcal{H}$.m1). There exist three constants $\tilde{k}$, $\tilde{\lambda}$ and $\tilde{\Lambda}$ and a nonincreasing function $\widetilde{K}$, only depending on $k$, $K$, $\lambda$ and $\Lambda$, such that $C_m$, $F_m$, $h_m$ and $\sigma_m$ are smooth, with bounded derivatives of any order and satisfy properties ($\mathcal{H}$.1)–($\mathcal{H}$.5) with respect to $\tilde{k}$, $\widetilde{K}$, $\tilde{\lambda}$ and $\tilde{\Lambda}$.

($\mathcal{H}$.m2). $(C_m, F_m, \sigma_m, h_m)_{m \in \mathbb{N}}$ uniformly converges on every compact set toward $(C, F, \sigma, H)$.

Then, we know from Ma, Protter and Yong (1994) [see also Corollary B.7 given in Delarue (2002b)] that, for every $m \in \mathbb{N}$, the system of PDEs:

For $(t, x) \in [0, T[ \times \mathbb{R}^P$ and $\ell \in \{1, \ldots, Q\}$,

$$\frac{\partial(\varphi_m)_\ell}{\partial t}(t, x) + \frac{1}{2}(a_m)_{i,j}(x, \varphi_m(t, x))\frac{\partial^2(\varphi_m)_\ell}{\partial x_i \partial x_j}(t, x)$$

$\mathcal{S}(m)$ $\qquad + (C_m)_i(x, \varphi_m(t, x), \nabla_x \varphi_m(t, x))\dfrac{\partial(\varphi_m)_\ell}{\partial x_i}(t, x)$

$$+ (F_m)_\ell(\varphi_m(t, x), \nabla_x \varphi_m(t, x)) = 0,$$

$$\varphi_m(T, x) = h_m(x),$$

with $a_m = \sigma_m(\sigma_m)^*$, admits a unique bounded classical solution $\varphi_m \in \mathcal{C}^{1,2}([0, T] \times \mathbb{R}^P, \mathbb{R}^Q)$. Moreover, thanks to Lemma 2.1 of Delarue (2002b), we know that

(D.2) $$\sup_{m \in \mathbb{N}} \|\varphi_m\|_{\mathcal{C}^{(1/2,1)}([0,T] \times \mathbb{R}^P, \mathbb{R}^Q)} < \infty.$$

Thanks to Theorem 5.1 of Chapter VII of Ladyzenskaja, Solonnikov and Ural'ceva (1968), we deduce that there exists $0 < \gamma \le 1$, only depending on $\bar{\lambda}$, $\bar{\Lambda}$ and $P$, such that

(D.3) $$\forall 0 < \delta < T, \qquad \sup_{m \in \mathbb{N}} \|\nabla_x \varphi_m\|_{\mathcal{C}^{(\gamma/2, \gamma)}([0, T-\delta] \times \mathbb{R}^P, \mathbb{R}^Q)} < \infty.$$

Finally, thanks to Theorem 9.1 of Chapter IV of Ladyzenskaja, Solonnikov and Ural'ceva (1968), we prove that

(D.4)
$$\forall q \ge 2, \ \forall 0 < \delta < T, \ \forall R > 0,$$
$$\sup_{m \in \mathbb{N}} \|\varphi_m\|_{W^{1,2,q}(]0, T-\delta[ \times B_P(0,R), \mathbb{R}^Q)} < \infty.$$

Hence, using a compactness argument, we deduce the existence of a solution $\varphi$ to the system $\mathcal{E}(1)$ in the space $\mathcal{V}$.

Let us now prove that $\varphi = \theta_1$ (see Section 2.1 for the definition of $\theta_1$). For this purpose, fix $(t, x) \in [0, T[ \times \mathbb{R}^P$. Since $\varphi$ and $\nabla_x \varphi$ are bounded, there exists, thanks to Theorem 6.5.1 of Stroock and Varadhan (1979), a triple $((U, W), (\widetilde{\Omega}, \widetilde{\mathcal{F}}, \widetilde{\mathbf{P}}))$, which is a weak solution of the SDE

(D.5)
$$\forall s \in [t, T], \qquad U_s = U_t + \int_t^s C(U_r, \varphi(r, U_r), \nabla_x \varphi(r, U_r)) \, dr$$
$$+ \int_t^s \sigma(U_r, \varphi(r, U_r)) \, dW_r.$$

Let us define for every $n \in \mathbb{N}$,

(D.6) $$\tau(n) = \inf\{t \le s \le T, \ |X_s - x| \ge n\}, \qquad (\inf \varnothing = T).$$

Thanks to the Itô–Krylov formula [Theorem 1, Section 10, Chapter II of Krylov (1980)] and to the system of PDEs $\mathcal{E}(1)$, we know that for all $n \in \mathbb{N}$ and $0 < \delta < T - t$,

$$
\begin{aligned}
\forall s \in [t, \tau(n) \wedge (T-\delta)], \\
\varphi(s, U_s) = \varphi(\tau(n) &\wedge (T-\delta), U_{\tau(n) \wedge (T-\delta)}) \\
&+ \int_s^{\tau(n) \wedge (T-\delta)} F(U_r, \varphi(r, U_r), \nabla_x \varphi(r, U_r))\, dr \\
&- \int_s^{\tau(n) \wedge (T-\delta)} \nabla_x \varphi(r, U_r) \sigma(U_r, \varphi(r, U_r))\, dW_r.
\end{aligned}
\tag{D.7}
$$

Letting $\delta \to 0$ and $n \to +\infty$, we deduce from the boundedness of $\varphi$ and $\nabla_x \varphi$ and the continuity of $\varphi$ that

$$
\begin{aligned}
\forall s \in [t, T], \qquad \varphi(s, U_s) = H(U_T) &+ \int_s^T F(U_r, \varphi(r, U_r), \nabla_x \varphi(r, U_r))\, dr \\
&- \int_s^T \nabla_x \varphi(r, U_r) \sigma(U_r, \varphi(r, U_r))\, dW_r.
\end{aligned}
\tag{D.8}
$$

Hence, $(U_s, \varphi(s, U_s), \nabla_x \varphi(s, U_s))_{t \leq s \leq T}$ is a "weak solution" of the FBSDE $\mathrm{E}(1, t, x)$. From Delarue [(2002b), Remark 2.7] we deduce that $\varphi = \theta_1$.

**Regularity of $(\nabla_x \theta_\varepsilon)_{\varepsilon > 0}$ and $\nabla_x \theta$.** We first investigate the regularity of $\nabla_x \theta_1$: thanks to (D.3), we deduce that $\theta_1$ is continuously differentiable with respect to $x$ on $[0, T[ \times \mathbb{R}^P$, and that $\nabla_x \theta_1$ is Hölderian on every set $[0, T - \eta] \times \mathbb{R}^P$, with $0 < \eta < T$. Actually, this is also true for the functions $(\theta(\varepsilon))_{\varepsilon > 0}$ and $\theta$: in particular, there exists a constant $0 < \beta \leq 1$, only depending on $\bar{\lambda}$, $\bar{\Lambda}$ and $P$, such that $\nabla_x \theta$ satisfies (2.13) for every $\eta > 0$.

**Regularity of $(\zeta_n)_{n \in \mathbb{N}^*}$.** Let us finally investigate the special case of $\zeta_n$ for a given $n \in \mathbb{N}^*$ (see Section 3.3 for the definition of $\zeta_n$): keep the notation introduced in the section "Solvability of $(\mathcal{E}(\varepsilon))_{\varepsilon > 0}$, $\mathcal{E}(\lim)$ and $(\mathcal{E}_{\mathrm{reg}}(n))_{n \in \mathbb{N}^*}$" and assume in addition to $(\mathcal{H}.1)$–$(\mathcal{H}.5)$ that the coefficient $F$ satisfies $\forall (x, y, y', z) \in \mathbb{R}^P \times \mathbb{R}^Q \times \mathbb{R}^Q \times \mathbb{R}^{Q \times P}$,

$$
|F(x, y, z) - F(x, y', z)| \leq K(|y| + |y'| + |z|)|y - y'|,
\tag{D.9}
$$

and that $H$ is smooth, with bounded derivatives of every order. Then, from Theorem 5.1, Chapter VII of Ladyzenskaja, Solonnikov and Ural'ceva (1968), we can assume without loss of generality that

$$
\forall 0 < \eta < 1, \qquad \sup_{m \in \mathbb{N}} \|\varphi_m\|_{\mathcal{C}^{(1+\eta/2, 2+\eta)}([0,T] \times \mathbb{R}^P, \mathbb{R}^Q)} < \infty.
\tag{D.10}
$$

Applying this regularization procedure to $\mathcal{E}_{\mathrm{reg}}(n)$, we deduce that, for every $0 < \eta < 1$, $\zeta_n \in \mathcal{C}^{(1+\eta/2, 2+\eta)}([0, T] \times \mathbb{R}^P, \mathbb{R}^Q)$ and that (3.8) holds.

**Convergence of** $(\zeta_n)_{n \in \mathbb{N}^*}$. It is readily seen that (D.2)–(D.4) hold with $\varphi_m$ replaced by $\zeta_m$. Hence, from a compactness argument, $(\zeta_n)_{n \in \mathbb{N}^*}$ uniformly converges on every compact subset of $[0, T] \times \mathbb{R}^P$ toward $\theta$ and $(\nabla_x \zeta_n)_{n \in \mathbb{N}^*}$ uniformly converges on every compact subset of $[0, T[ \times \mathbb{R}^P$ toward $\nabla_x \theta$.

U.F.R. DE MATHEMATIQUES
UNIVERSITÉ PARIS VII
CASE 7012
2 PLACE JUSSIEU
75 251 PARIS CEDEX 05
FRANCE
E-MAIL: delarue@math.jussieu.fr